\documentclass[a4paper,fleqn,numbers]{cas-sc}
\usepackage{amsmath,amssymb,amsfonts,amsthm}
\usepackage{bm}

\hypersetup{
	colorlinks=true,
	linkcolor=blue,
	citecolor=blue,
	urlcolor=blue
}

\usepackage[nameinlink,capitalise,noabbrev]{cleveref}

\usepackage{graphicx}
\usepackage{enumitem}
\usepackage{booktabs}
\usepackage{mathtools}
\let\mathscr\mathcal 
\usepackage{subcaption}
\usepackage{algorithm}
\usepackage{algpseudocode}
\usepackage{multirow}
\usepackage{tikz}
\usepackage[authoryear]{natbib}
\usepackage{float}
\newtheorem{theorem}{Theorem}
\newtheorem{proposition}{Proposition}
\newtheorem{lemma}{Lemma}
\newtheorem{corollary}{Corollary}

\newtheorem{remark}{Remark}
\usepackage{float}     
\usepackage{placeins}  
\newcommand{\R}{\mathbb{R}}

\newcommand{\LD}{L_{\mathrm{ref}}}
\newcommand{\LK}{L_{\mathrm{KNN}}}
\graphicspath{{figures/v5_20261114_V2/}}
\setcitestyle{numbers,square,sort&compress}
\usepackage[mathlines]{lineno}
\begin{document}
	\let\WriteBookmarks\relax
	\shorttitle{Structure-Preserving PHS RBF-FD Methods}
	\shortauthors{G. Yao and Z. Diefendorf}
	\title[mode=title]{Structure-Preserving PHS RBF-FD Methods for Schr\"odinger Equations}
	
	\author[1]{Guangming Yao}[orcid = 0000-0002-4819-0101]
	\cormark[1]
	\author[1]{Zachary Diefendorf}
	
	\affiliation[1]{organization={Department of Mathematics, Clarkson University},
		city={Potsdam},state={New York},country={USA}}
	\cortext[cor1]{Corresponding author: gyao@clarkson.edu}
	
	\begin{abstract}
		We develop and analyze a structure-preserving polyharmonic-spline radial-basis-function finite-difference (PHS RBF-FD) method for the two-dimensional Schrödinger equation with homogeneous Dirichlet boundary conditions. The method combines a 29-point complete-shell stencil with odd boundary reflection to construct a symmetric discrete Laplacian. Unlike conventional nearest-neighbor RBF-FD discretizations, whose independently selected stencils do not generally preserve global symmetry, the proposed construction admits an exact discrete sine eigenbasis. We derive a Fourier-symbol representation and establish a sufficient criterion for strict negative definiteness, verified for the \(r^7+\mathbb P_4\) and \(r^{11}+\mathbb P_5\) formulations using validated coefficient calculations. The resulting symmetric spatial operator preserves discrete mass and energy under exact semidiscrete evolution, while fully discrete conservation depends on the time integrator. Numerical experiments examine kernel selection, dispersion, convergence, symmetry defects, and computational cost. Comparisons with conventional RBF-FD and fourth-order finite differences include exact Fourier modes and a variable-potential Schrödinger problem with an independently refined spectral reference. The \(r^{11}+\mathbb P_5\) formulation provides favorable accuracy in the tested problems while maintaining the desired spatial structure. The analysis establishes a connection between boundary treatment, global operator symmetry, and structure preservation in localized PHS discretizations.
	\end{abstract}
	
	\maketitle

	
	\section{Introduction}
	
	Radial basis function (RBF) methods provide flexible
	meshfree approximations for partial differential
	equations (PDEs), particularly on scattered nodes
	and geometrically complex domains
	\cite{FornbergFlyer2015,Wendland2005,Fasshauer2007}.
	Among these methods, radial basis function-generated
	finite differences (RBF-FD) combine local RBF
	interpolation with sparse differentiation matrices.
	Polynomial augmentation improves approximation
	accuracy and provides exact reproduction of
	polynomials up to a prescribed degree
	\cite{FlyerEtAl2016,BayonaEtAl2017}.
	Polyharmonic splines (PHS) are especially attractive
	because they do not require a shape parameter,
	thereby avoiding the associated parameter-selection
	problem.
	
	\paragraph{Related developments in localized RBF methods.}
	Localized RBF approximations and implicit meshless
	methods have been investigated extensively for
	the numerical solution of PDEs.
	\cite{ChenFanWen2012}
	developed the method of approximate particular
	solutions for solving PDEs.
	Related localized formulations demonstrated how
	local RBF approximations could be assembled into
	sparse global systems.
	\cite{YaoDuoChenShen2015}
	developed implicit local RBF interpolation methods,
	while subsequent work employed polyharmonic splines
	in modified particular-solution formulations
	\cite{YaoChenZheng2017}.
	More recently, \cite{RubasingheYaoNiuTsogtgerel2023}
	investigated PHS interpolation on scattered nodes.
	These developments demonstrate the effectiveness
	of localization and polynomially augmented RBF
	approximations. However, the present work addresses
	a different question: whether the assembled
	differentiation operator preserves the algebraic
	structure required by a Hamiltonian evolution
	equation.
	
	\paragraph{Structure preservation for
		Schr\"odinger equations.}
	The Schr\"odinger equation provides an important
	setting in which the structure of a spatial
	discretization directly influences numerical
	stability and conservation properties.
	Its continuous Hamiltonian operator is
	self-adjoint, and the corresponding evolution
	is unitary.
	Consequently, the exact solution preserves
	mass and energy.
	Structure-preserving discretizations of
	Schr\"odinger equations have been studied
	extensively, including conservative numerical
	schemes \cite{SanzSerna1984},
	splitting methods \cite{Lubich2008},
	and time-splitting spectral approximations
	\cite{BaoJinMarkowich2002,BaoJinMarkowich2003}.
	The distinction between semidiscrete invariants
	and their preservation under time integration
	is also central to geometric numerical
	integration \cite{HairerLubichWanner2006}.
	Loss of global symmetry can alter the
	spectral properties of the discrete
	Schr\"odinger generator and may introduce
	artificial growth or decay.
	Thus, local approximation accuracy alone
	does not guarantee preservation of the
	underlying Hamiltonian structure.
	
	\paragraph{Related RBF-FD research and
		the present approach.}
	For RBF-FD discretizations, however,
	the algebraic structure of the assembled
	spatial operator requires particular attention.
	Although the polynomially augmented RBF
	interpolation matrix is symmetric on each
	local stencil, independently computed
	differentiation weights do not generally
	produce a symmetric global matrix.
	To see this distinction, let $w_{ij}$
	denote the differentiation weight assigned
	to node $j$ by the stencil centered at node $i$.
	Global symmetry requires
	\(
	w_{ij}=w_{ji}
	\)
	for every pair of interior nodes.
	Standard nearest-neighbor (KNN) stencil
	selection does not enforce this reciprocity.
	The difficulty is especially pronounced
	near boundaries, where one-sided stencil
	geometries and distance-tie choices
	can destroy the reflection symmetries
	present in interior stencils. 
	The influence of polynomial augmentation
	and boundary geometry on RBF-FD accuracy
	has been investigated in
	\cite{FlyerEtAl2016,BayonaEtAl2017,
		BayonaEtAl2019}.
	Stabilization techniques for RBF-FD
	differentiation operators have also
	been developed for advection--diffusion
	problems
	\cite{ShankarFogelson2018,
		ShankarWrightNarayan2020}.

	In this paper, we consider a different
	approach based on the geometry of the
	local stencil and the treatment of
	homogeneous Dirichlet boundaries.
	On a uniform Cartesian grid, we employ
	a 29-point complete-shell stencil whose
	geometry is invariant under coordinate
	reflections and rotations.
	The same reference differentiation
	weights are used at every interior node.
	Stencil points outside the physical
	domain are eliminated through odd
	reflection across the boundary.
	This construction produces a symmetric
	global discrete Laplacian without
	postprocessing the assembled matrix
	by replacing $L$ with $(L+L^T)/2$.
	Moreover, the resulting operator admits
	an exact discrete sine eigenbasis,
	allowing its spectral properties to
	be analyzed through a Fourier symbol.
	Our analysis focuses on the PHS
	formulations
	$r^7+\mathbb P_4$ and
	$r^{11}+\mathbb P_5$.
	For these formulations, we establish
	a sufficient coefficient-based criterion
	for strict negative definiteness of
	the reflected discrete Laplacian.
	The coefficient conditions are verified
	using validated numerical calculations.
	
	\paragraph{Contributions and scope.}
	The principal contributions of this work
	are summarized as follows.
	
	\begin{enumerate}[label=(\roman*),leftmargin=*]\itemsep -0.15cm
		
		\item We construct a reflected PHS RBF-FD
		Laplacian using a complete-shell stencil
		and odd boundary extension, and prove
		symmetry of the assembled interior matrix.
		
		\item We derive the discrete Fourier symbol
		and establish a sufficient criterion
		for strict negative definiteness,
		applying it to the
		$r^7+\mathbb P_4$ and
		$r^{11}+\mathbb P_5$ formulations.
		
		\item We analyze the dispersion properties
		of the reflected discretizations,
		including their fixed-mode asymptotic
		behavior and the influence of kernel
		and polynomial selection.
		
		\item We investigate numerical accuracy,
		symmetry defects, and conservation
		properties through comparisons with
		standard KNN RBF-FD and fourth-order
		finite differences.
		
		\item We assess the proposed method on
		both exact Fourier-mode solutions
		and a variable-potential
		Schr\"odinger benchmark, while
		examining computational cost and
		the limitations of the construction.
		
	\end{enumerate}
	
	The present analysis is restricted to
	uniform Cartesian grids on a square
	domain with homogeneous Dirichlet
	boundary conditions.
	The resulting symmetry and spectral
	properties rely on this geometry and
	should not be assumed for general
	scattered-node discretizations.
	Nevertheless, the results provide
	a foundation for investigating
	structure-preserving localized RBF
	methods in more general settings.
	The variable-potential Schr\"odinger benchmark
	further evaluates the method beyond the
	free equation, demonstrating convergence
	and favorable accuracy relative to the
	comparison discretizations.

	\section{Mathematical model and PHS RBF-FD discretization}
	
	We consider the linear Schr\"odinger equation on a square with homogeneous
	Dirichlet boundary conditions. This setting permits exact Fourier-mode
	benchmarks and, on a uniform grid, a direct examination of how local
	stencil geometry and boundary treatment affect the assembled Laplacian.
	We first state the continuous model and its invariants, then introduce
	the polynomially augmented PHS RBF-FD weights and the two spatial
	constructions used in the numerical comparisons.

	Let $\Omega=(0,\pi)^2$, and write the complex
	wave function as $\psi=u+iv$, where $u$ and $v$
	are real-valued functions. We consider the free
	Schr\"odinger equation
	\begin{equation}
		i\psi_t=-\Delta\psi,
		\qquad (x,y)\in\Omega,\quad t>0,
		\label{eq:complexSchrodinger}
	\end{equation}
	subject to homogeneous Dirichlet boundary conditions
	\begin{equation}
		\psi|_{\partial\Omega}=0,
		\label{eq:complexDirichlet}
	\end{equation}
	and smooth, boundary-compatible initial data
	\begin{equation}
		\psi(x,y,0)=\psi_0(x,y)
		=u_0(x,y)+iv_0(x,y).
		\label{eq:complexInitial}
	\end{equation}
	
	Separating the real and imaginary parts of
	\eqref{eq:complexSchrodinger} yields the
	equivalent coupled real-valued system
	\begin{align}
		u_t&=-\Delta v,\label{eq:cont1}\\
		v_t&=\phantom{-}\Delta u,\label{eq:cont2}
	\end{align}
	with boundary conditions
	\begin{equation}
		u|_{\partial\Omega}=0,\qquad
		v|_{\partial\Omega}=0,
		\label{eq:dirichlet}
	\end{equation}
	and initial conditions
	$u(\cdot,0)=u_0$ and $v(\cdot,0)=v_0$.

	\subsection{Fourier solution and exact-mode benchmarks}
	\label{subsec:exactModes}
	
	On $\Omega=(0,\pi)^2$, the Dirichlet Laplacian
	has the sine eigenfunctions
	\[
	\phi_{k\ell}(x,y)=\sin(kx)\sin(\ell y),
	\qquad
	-\Delta\phi_{k\ell}
	=(k^2+\ell^2)\phi_{k\ell},
	\quad k,\ell\ge1.
	\]
	We expand the real and imaginary parts of the
	solution in this eigenbasis:
	\[
	u(x,y,t)=\sum_{k,\ell\ge1}
	a_{k\ell}(t)\phi_{k\ell}(x,y),
	\qquad
	v(x,y,t)=\sum_{k,\ell\ge1}
	b_{k\ell}(t)\phi_{k\ell}(x,y).
	\]
	Substituting these expansions into
	$u_t=-\Delta v$ and $v_t=\Delta u$,
	and using the orthogonality of the sine
	eigenfunctions, gives
	\[
	a'_{k\ell}(t)=\omega_{k\ell}b_{k\ell}(t),
	\qquad
	b'_{k\ell}(t)=-\omega_{k\ell}a_{k\ell}(t),
	\qquad
	\omega_{k\ell}=k^2+\ell^2.
	\]
	Each pair of Fourier coefficients therefore
	satisfies a harmonic oscillator system,
	whose solution is
	\begin{align*}
		a_{k\ell}(t)
		&=A_{k\ell}\cos(\omega_{k\ell}t)
		+B_{k\ell}\sin(\omega_{k\ell}t),\\
		b_{k\ell}(t)
		&=-A_{k\ell}\sin(\omega_{k\ell}t)
		+B_{k\ell}\cos(\omega_{k\ell}t).
	\end{align*}
	Here $A_{k\ell}$ and $B_{k\ell}$ are the
	Fourier sine coefficients of the initial data
	$u_0$ and $v_0$, respectively:
	\begin{align*}
		A_{k\ell}
		&=\frac{4}{\pi^2}
		\int_0^\pi\!\int_0^\pi
		u_0(x,y)\phi_{k\ell}(x,y)\,dx\,dy,\\
		B_{k\ell}
		&=\frac{4}{\pi^2}
		\int_0^\pi\!\int_0^\pi
		v_0(x,y)\phi_{k\ell}(x,y)\,dx\,dy.
	\end{align*}
	Consequently, the exact solution in terms
	of its real and imaginary components is
	\begin{align}
		u(x,y,t)
		&=\sum_{k,\ell\ge1}
		\bigl[
		A_{k\ell}\cos(\omega_{k\ell}t)
		+B_{k\ell}\sin(\omega_{k\ell}t)
		\bigr]\phi_{k\ell}(x,y),
		\label{eq:fourierU}\\
		v(x,y,t)
		&=\sum_{k,\ell\ge1}
		\bigl[
		-A_{k\ell}\sin(\omega_{k\ell}t)
		+B_{k\ell}\cos(\omega_{k\ell}t)
		\bigr]\phi_{k\ell}(x,y).
		\label{eq:fourierV}
	\end{align}
	Equivalently, writing $\psi=u+iv$,
	the complex solution takes the form
	\[
	\psi(x,y,t)
	=\sum_{k,\ell\ge1}
	(A_{k\ell}+iB_{k\ell})
	e^{-i\omega_{k\ell}t}
	\phi_{k\ell}(x,y).
	\]
	In particular, for the single-mode initial
	condition
	\(
	u_0(x,y)=\sin(kx)\sin(\ell y),
	v_0(x,y)=0,
	\)
	only one Fourier coefficient is nonzero,
	with $A_{k\ell}=1$ and $B_{k\ell}=0$.
	The exact solution reduces to
	\begin{align}
		u(x,y,t)
		&=\cos\bigl((k^2+\ell^2)t\bigr)
		\sin(kx)\sin(\ell y),
		\label{eq:exactu}\\
		v(x,y,t)
		&=-\sin\bigl((k^2+\ell^2)t\bigr)
		\sin(kx)\sin(\ell y).
		\label{eq:exactv}
	\end{align}
	These exact modes provide benchmarks for
	evaluating spatial consistency, numerical
	dispersion, and accumulated phase error.

	\subsection{Continuous mass and energy conservation}
	\label{subsec:continuousInvariants}
	
	The continuous mass and Dirichlet energy are
	\begin{align}
		M(t)&=\int_\Omega\bigl(u^2+v^2\bigr)\,dx\,dy,
		\label{eq:continuousMass}\\
		E(t)&=\frac12\int_\Omega
		\bigl(|\nabla u|^2+|\nabla v|^2\bigr)\,dx\,dy.
		\label{eq:energy}
	\end{align}
	
	\begin{proposition}[Continuous mass and energy conservation]
		Assume that $u$ and $v$ are sufficiently smooth solutions of
		\eqref{eq:cont1}--\eqref{eq:dirichlet}. Then
		$M(t)=M(0)$ and $E(t)=E(0)$ for all times for which
		the solution exists.
	\end{proposition}
	\begin{proof}
		Since $\psi_t=i\Delta\psi$, integration by parts and the
		homogeneous boundary conditions give
		\[
		\frac{dM}{dt}=2\operatorname{Re}\int_\Omega
		\overline\psi\,i\Delta\psi\,dx\,dy
		=2\operatorname{Re}\left(-i\int_\Omega
		|\nabla\psi|^2\,dx\,dy\right)=0.
		\]
		Similarly, differentiation of the energy and integration by parts yield
		\begin{align*}
			\frac{dE}{dt}
			&=-\int_\Omega (\Delta u)u_t\,dx\,dy
			-\int_\Omega (\Delta v)v_t\,dx\,dy=\int_\Omega (\Delta u)(\Delta v)\,dx\,dy
			-\int_\Omega (\Delta v)(\Delta u)\,dx\,dy=0.
		\end{align*}
	\end{proof}

	These conservation laws provide benchmarks
	for assessing the spatial discretizations
	and time integrators developed below.
	We denote by $M_h(t)$ and $E_h(t)$
	the corresponding discrete mass and
	Dirichlet energy, respectively, where
	$h$ is the spatial grid spacing.
	Their precise definitions are introduced
	in Section~\ref{sec:semidiscrete}.
	We then establish conditions under which
	these discrete quantities are conserved
	exactly, preserving the corresponding
	conservation laws of the continuous
	Schr\"odinger equation.

	\subsection{PHS RBF-FD approximation}
	\label{sec:localPHS}
	
	This section describes the construction of local
	polyharmonic spline (PHS) radial basis function-generated
	finite difference (RBF-FD) approximations to the
	Laplacian. We first introduce the PHS kernels and
	their Laplacians, then formulate the polynomially
	augmented local differentiation system. Finally,
	we discuss polynomial unisolvency, kernel selection,
	and the distinction between local and global symmetry.
	
	\subsubsection{PHS kernels and their Laplacians}
	\label{subsec:phsFamilies}
	
	We consider two standard families of PHS kernels:
	the odd-power family
	\begin{equation}
		\phi(r)=r^{2n+1},
		\label{eq:oddPHS}
	\end{equation}
	and the even logarithmic family
	\begin{equation}
		\phi(r)=r^{2n}\log r,
		\label{eq:logPHS}
	\end{equation}
	where $n$ is a positive integer and $r$ denotes
	Euclidean distance.
	The principal numerical experiments use the
	odd-power kernels $r^7$ and $r^{11}$.
	The logarithmic family is included to describe
	the more general PHS RBF-FD framework.
	
	For a radial function $f(r)$ in two dimensions,
	the Laplacian takes the form
	\begin{equation}
		\Delta f(r)
		=f''(r)+\frac{1}{r}f'(r),
		\qquad r>0.
		\label{eq:radialLaplacian}
	\end{equation}
	
	\paragraph{Odd-power kernels.}
	
	Let $f(r)=r^q$, where $q=2n+1$.
	Differentiation gives
	\[
	f'(r)=qr^{q-1},
	\qquad
	f''(r)=q(q-1)r^{q-2}.
	\]
	Substitution into \eqref{eq:radialLaplacian}
	yields
	\begin{align*}
		\Delta(r^q)
		&=q(q-1)r^{q-2}
		+qr^{q-2} =q^2r^{q-2}.
	\end{align*}
	Thus,
	\begin{equation}
		\Delta(r^q)=q^2r^{q-2},
		\qquad r>0.
		\label{eq:lapPHS}
	\end{equation}
	In particular,
	\[
	\Delta(r^7)=49r^5,
	\qquad
	\Delta(r^{11})=121r^9.
	\]
	
	\paragraph{Even logarithmic kernels.}
	
	For $f(r)=r^p\log r$, where $p=2n$,
	we have
	\begin{align*}
		f'(r)
		&=pr^{p-1}\log r+r^{p-1},\\
		f''(r)
		&=p(p-1)r^{p-2}\log r
		+(2p-1)r^{p-2}.
	\end{align*}
	Applying \eqref{eq:radialLaplacian},
	\begin{align*}
		\Delta(r^p\log r)
		&=p(p-1)r^{p-2}\log r  +(2p-1)r^{p-2} +pr^{p-2}\log r+r^{p-2} =r^{p-2}
		\left(p^2\log r+2p\right).
	\end{align*}
	Therefore,
	\begin{equation}
		\Delta(r^p\log r)
		=r^{p-2}\left(p^2\log r+2p\right),
		\qquad r>0.
		\label{eq:logPHSLaplacian}
	\end{equation}
	For example,
	\[
	\Delta(r^4\log r)
	=r^2(16\log r+8),
	\]
	and
	\[
	\Delta(r^6\log r)
	=r^4(36\log r+12).
	\]
	
	For the odd powers $q\ge3$ and even
	logarithmic powers $p\ge4$ considered here,
	both the kernels and their Laplacians extend
	continuously to $r=0$, where the corresponding
	Laplacian values are zero.
	
	\subsubsection{Polynomial augmentation and local weights}
	\label{subsec:localWeights}
	
	Let $\Pi_m(\R^2)$ denote the space of
	bivariate polynomials of total degree at most
	$m$. Its dimension is
	\begin{equation}
		p_m=\dim\Pi_m(\R^2)
		=\frac{(m+1)(m+2)}{2}.
		\label{eq:polynomialDimension}
	\end{equation}
	In particular,
	\(
	\dim\Pi_4(\R^2)=15,
	\dim\Pi_5(\R^2)=21.
	\)
	For a stencil
	\(
	\mathcal N(\bm x_0)
	=\{\bm x_1,\ldots,\bm x_{n_s}\}
	\)
	centered at $\bm x_0$, we approximate
	the Laplacian by
	\begin{equation}
		(\Delta f)(\bm x_0)
		\approx
		\sum_{j=1}^{n_s}
		w_j(\bm x_0)f(\bm x_j).
		\label{eq:lapapprox}
	\end{equation}
	Let
	$\{p_\alpha\}_{\alpha=1}^{p_m}$
	be a basis of $\Pi_m(\R^2)$.
	Define the RBF interpolation matrix
	$A_1\in\R^{n_s\times n_s}$ and
	polynomial matrix
	$P\in\R^{n_s\times p_m}$ by
	\[
	(A_1)_{ij}
	=\phi(\|\bm x_i-\bm x_j\|),
	\qquad
	P_{i\alpha}=p_\alpha(\bm x_i).
	\]
	The local differentiation weights
	$\bm w=(w_1,\ldots,w_{n_s})^T$
	are obtained from the augmented system
	\begin{equation}
		\begin{pmatrix}
			A_1&P\\
			P^T&0
		\end{pmatrix}
		\begin{pmatrix}
			\bm w\\
			\bm\lambda
		\end{pmatrix}
		=
		\begin{pmatrix}
			\bm b\\
			\bm c
		\end{pmatrix},
		\label{eq:blocksys}
	\end{equation}
	where
	\(
	b_j=
	\left.
	\Delta_{\bm x}
	\phi(\|\bm x-\bm x_j\|)
	\right|_{\bm x=\bm x_0},
	c_\alpha=
	(\Delta p_\alpha)(\bm x_0).
	\)
	The first block row enforces the
	RBF differentiation conditions, while
	the second imposes exact polynomial
	reproduction:
	\begin{equation}
		\sum_{j=1}^{n_s}
		w_jp(\bm x_j)
		=(\Delta p)(\bm x_0),
		\qquad
		p\in\Pi_m(\R^2).
		\label{eq:polynomialReproduction}
	\end{equation}
	Consequently, the resulting differentiation
	formula is exact for every polynomial
	in the augmentation space.
	
	\paragraph{Unisolvency and nonsingularity.}
	
	The polynomial matrix $P$ must have full
	column rank for the stencil to be
	$\Pi_m$-unisolvent. In particular,
	\[
	n_s\ge\dim\Pi_m(\R^2)
	\]
	is necessary but not sufficient:
	the geometric arrangement of the stencil
	nodes must also support independent
	polynomial reproduction constraints.
	For the PHS kernels considered here,
	the standard conditional-definiteness
	theory guarantees nonsingularity of
	the augmented system for distinct,
	polynomially unisolvent stencil nodes,
	provided the polynomial space contains
	the degree required by the kernel.
	For the principal combinations used
	in this work, $r^7+\Pi_4$ and
	$r^{11}+\Pi_5$, these requirements
	are satisfied whenever the corresponding
	polynomial matrices have full column rank.
	
	Near physical boundaries, a naive
	fixed-size KNN stencil may fail this
	rank condition. The matched KNN
	comparisons therefore use rank-safe
	stencil selection to ensure that
	the polynomial reproduction constraints
	are well-defined before the local
	weights are computed.
	
	\subsubsection{Kernel selection and symmetry considerations}
	\label{subsec:kernelSelection}
	
	The choice of PHS kernel, polynomial
	augmentation, and stencil geometry
	jointly determines the resulting
	differentiation weights. Increasing
	the kernel exponent or polynomial
	degree does not automatically improve
	accuracy or stability; the augmented
	system must remain nonsingular, and
	the assembled spatial operator must
	be examined separately.
	
	\begin{figure}[pos=htbp]
		\centering
		
		\begin{subfigure}[H]{0.46\textwidth}
			\centering
			\begin{tikzpicture}[
				scale=0.52,
				>=stealth,
				gridpoint/.style={
					circle, fill=gray!35, inner sep=1.8pt
				},
				stencilpoint/.style={
					circle, fill=blue!75, inner sep=2.6pt
				},
				centerpoint/.style={
					circle, fill=red!85, inner sep=3.5pt
				}
				]
				
				\draw[step=1,gray!25,thin] (-3.5,-3.5)
				grid (3.5,3.5);
				
				\draw[blue!55,dashed,thick]
				(0,0) circle (3);
				
				\foreach \i in {-3,...,3}{
					\foreach \j in {-3,...,3}{
						\node[gridpoint] at (\i,\j) {};
					}
				}
				
				\foreach \i in {-3,...,3}{
					\foreach \j in {-3,...,3}{
						\pgfmathtruncatemacro{\radsq}{\i*\i+\j*\j}
						\ifnum\radsq<10
						\node[stencilpoint] at (\i,\j) {};
						\fi
					}
				}
				
				\node[centerpoint] at (0,0) {};
				
				\draw[->,thick] (-3.5,0) -- (3.65,0)
				node[right] {$a$};
				\draw[->,thick] (0,-3.5) -- (0,3.65)
				node[above] {$b$};
				
				\node[blue!70!black,fill=white,inner sep=1pt]
				at (2.1,2.45) {$a^2+b^2\leq 9$};
				
			\end{tikzpicture}
			\caption{Complete-shell reference stencil $S_{29}$.}
			\label{fig:stencil29}
		\end{subfigure}
		\hfill
		\begin{subfigure}[H]{0.46\textwidth}
			\centering
			\begin{tikzpicture}[
				scale=0.52,
				>=stealth,
				physical/.style={
					circle, fill=blue!75, inner sep=2.6pt
				},
				exterior/.style={
					circle, draw=orange!90!black,
					fill=orange!25, inner sep=2.6pt
				},
				centerpoint/.style={
					circle, fill=red!85, inner sep=3.5pt
				}
				]
				
				\fill[blue!4] (0,-3.5) rectangle (4.5,3.5);
				
				\draw[step=1,gray!25,thin]
				(-2.5,-3.5) grid (4.5,3.5);
				
				\draw[very thick,black]
				(0,-3.5) -- (0,3.5);
				\node[rotate=90,above]
				at (0,2.5) {$x=0,\quad \psi=0$};
				
				\foreach \a in {-3,...,3}{
					\foreach \b in {-3,...,3}{
						\pgfmathtruncatemacro{\radsq}{\a*\a+\b*\b}
						\ifnum\radsq<10
						\pgfmathtruncatemacro{\xx}{1+\a}
						\ifnum\xx<0
						\node[exterior] at (\xx,\b) {};
						\else
						\node[physical] at (\xx,\b) {};
						\fi
						\fi
					}
				}
				
				\node[centerpoint] at (1,0) {};
				
				\draw[->,orange!90!black,thick,bend left=30]
				(-1,1) to (1,1);
				
				\node[orange!90!black,fill=white,inner sep=1pt]
				at (-2.5,1.5) {\scriptsize odd reflection};
				
				\node[fill=white,align=center,inner sep=1.5pt]
				at (-1,-2.4)
				{\scriptsize $\psi(-x,y)=-\psi(x,y)$};
				
				\node[orange!90!black] at (-1.35,-3.15)
				{\scriptsize exterior};
				\node[blue!70!black] at (2.75,3.15)
				{\scriptsize physical};
				
			\end{tikzpicture}
			\caption{Odd reflection at a Dirichlet boundary.}
			\label{fig:stencil-reflection}
		\end{subfigure}
		\caption{Geometry of the structure-preserving
			reflected RBF-FD construction.
			(a) The complete-shell reference stencil
			$S_{29}$ contains 29 lattice points and is
			invariant under coordinate reflections,
			as defined in \eqref{eq:S29}.
			(b) Near the left homogeneous Dirichlet
			boundary, the same reference stencil is
			retained, including stencil positions outside
			the physical domain. These exterior values
			are eliminated by odd reflection,
			$\psi(-x,y)=-\psi(x,y)$, while the boundary
			values vanish. Consequently, contributions
			from exterior stencil positions are transferred
			to their reflected interior counterparts
			with opposite signs. The reflected boundary
			construction is described in detail in
			Section~\ref{reflectedbond}.}
		\label{fig:reflected-stencil-geometry}
	\end{figure}
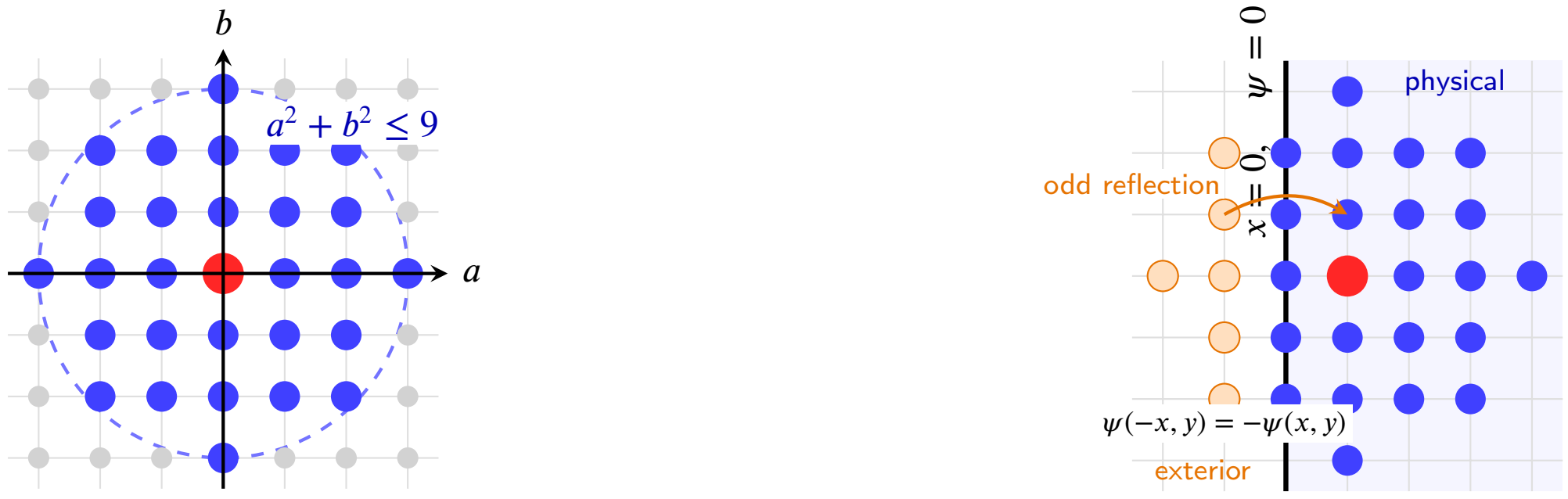
	Our objective is to obtain an accurate
	discrete Laplacian while preserving
	the symmetry of the continuous Laplacian.
	For this purpose, we use a complete-shell
	stencil whose geometric symmetry, together
	with odd boundary reflection, allows the
	assembled discrete operator to remain
	symmetric.
	Our principal reflected construction uses the
	odd-power kernel $r^{11}$ with degree-five
	polynomial augmentation $\Pi_5$ on the
	29-point complete lattice stencil
	\begin{equation}
		S_{29}
		=
		\{(a,b)\in\mathbb Z^2:
		a^2+b^2\le 9\}.
		\label{eq:S29}
	\end{equation}
	Here, $(a,b)$ denotes an integer offset
	relative to the center of the stencil,
	rather than the absolute coordinates
	of a grid point. For a stencil centered
	at $(x_i,y_j)$ with grid spacing $h$,
	the offset $(a,b)$ corresponds to the
	physical point
	\[
	(x_i+ah,\;y_j+bh).
	\]
	Since the stencil is defined by
	$a^2+b^2\le 9$, it follows that
	\[
	(a,b)\in S_{29}
	\quad\Longrightarrow\quad
	(-a,-b)\in S_{29}.
	\]
	The corresponding physical points
	$(x_i+ah,y_j+bh)$ and
	$(x_i-ah,y_j-bh)$ are symmetric
	about the stencil center $(x_i,y_j)$.
	Thus, $S_{29}$ is centrally symmetric
	in local stencil coordinates, see~\ref{fig:reflected-stencil-geometry}.
	Panel~\subref{fig:stencil29}. 
	
	This combination is chosen because the
	degree-five polynomial augmentation provides
	the moment conditions needed for fourth-order
	Laplacian consistency, while the complete-shell
	geometry preserves the reflection symmetries
	required for a symmetric global operator. 
	The radius-three condition selects exactly
	29 integer lattice offsets, providing
	a compact stencil with sufficient nodes
	for the 21-dimensional polynomial space
	$\Pi_5$. The stencil is invariant under
	coordinate reflections and rotations
	by $90^\circ$, properties that are
	central to the symmetry and spectral
	analysis developed below.
	The combination $r^7+\Pi_4$ is also
	examined to assess how kernel and
	polynomial selection affect accuracy
	and dispersion.

	\paragraph{Centrosymmetry and odd polynomial moments.}
	
	The reference stencil $S_{29}$ is centrosymmetric:
	if $(a,b)\in S_{29}$, then $(-a,-b)\in S_{29}$.
	For a radial kernel and the Laplacian evaluated
	at the stencil center, the local differentiation
	problem is invariant under central reflection.
	Provided the augmented system is nonsingular,
	uniqueness implies
	\begin{equation}
		w_{a,b}=w_{-a,-b}.
		\label{eq:centrosymmetricWeights}
	\end{equation}
	Consequently, for every polynomial $p(a,b)$
	satisfying $p(-a,-b)=-p(a,b)$,
	\begin{equation}
		\sum_{(a,b)\in S_{29}}
		w_{a,b}p(a,b)=0.
		\label{eq:oddMoments}
	\end{equation}
	In particular, all monomials of total degree
	five satisfy this identity.
	This observation has a useful consequence for
	polynomial augmentation. On the reference stencil,
	the unique $\Pi_4$-augmented Laplacian weights
	automatically satisfy the additional degree-five
	reproduction conditions, whose Laplacians vanish
	at the stencil center. Thus, provided both
	augmented systems are nonsingular, the $\Pi_4$
	and $\Pi_5$ formulations produce identical
	Laplacian weights in exact arithmetic.
	This equivalence is a consequence of the
	centrosymmetric stencil geometry and cannot
	be assumed for arbitrary KNN stencils.
	We therefore retain $\Pi_5$ in the principal
	$r^{11}$ RBF-FD formulation.

	\paragraph{Assembly of the discrete Laplacian.}
	
	The RBF-FD procedure computes differentiation
	weights locally, one stencil at a time.
	Let $\bm{x}_1,\ldots,\bm{x}_M$ denote the
	interior grid nodes. For each interior
	center $\bm{x}_i$, the local RBF-FD weights
	approximate the Laplacian by
	\[
	\Delta u(\bm{x}_i)
	\approx
	\sum_{\bm{x}_j\in\mathcal S_i}
	w_{ij}u(\bm{x}_j),
	\]
	where $\mathcal S_i$ is the stencil centered
	at $\bm{x}_i$. The weights associated with
	interior unknowns are inserted into row $i$
	of the global matrix $L$. Thus,
	\[
	(L\bm{u})_i
	=
	\sum_{j=1}^{M}L_{ij}u(\bm{x}_j),
	\]
	where $L_{ij}=w_{ij}$ when $\bm{x}_j$ is an
	interior node in the stencil of $\bm{x}_i$,
	and $L_{ij}=0$ otherwise. Boundary
	contributions are treated according to
	the prescribed boundary conditions.
	
	Because the weights are computed separately
	at different stencil centers, the matrix
	entry $L_{ij}$ and its counterpart $L_{ji}$
	need not be equal. This motivates the
	following distinction between local
	algebraic symmetry and global operator
	symmetry.

	\begin{remark}[Local versus global symmetry]
		Two distinct notions of symmetry arise in the
		RBF-FD construction:
		\begin{enumerate}[label=(\roman*),leftmargin=*]
			
			\item \textbf{Local algebraic symmetry.}
			The augmented matrix in \eqref{eq:blocksys}
			is symmetric because $A_1=A_1^T$ and the
			polynomial blocks are $P$ and $P^T$.
			This property holds for each local
			weight-computation system.

			\item \textbf{Global operator symmetry.}
			The assembled Laplacian $L$ is symmetric
			if $L_{ij}=L_{ji}$ for all interior nodes
			$i,j$. This property does not follow
			from symmetry of the local augmented
			systems. It must be established for the
			assembled operator.

		\end{enumerate}
		The two spatial discretization strategies
		introduced below differ fundamentally in
		whether this global symmetry is preserved.
	\end{remark}

	
	\subsection{Rank-safe KNN discretization}
	
	A common approach in RBF-FD methods is to
	construct local stencils using the
	$k$-nearest-neighbor (KNN) rule, which
	selects nodes according to their Euclidean
	distance from each stencil center.
	This approach is simple and flexible,
	particularly for scattered or irregularly
	distributed nodes.
	
	In the present implementation, we use a
	rank-safe modification of KNN to ensure
	that the selected nodes support the
	degree-five polynomial augmentation.
	For each interior center, candidate nodes
	are first ordered by increasing squared
	Euclidean distance. Nodes at equal
	distances are ordered lexicographically
	by their grid indices, providing a
	deterministic tie-breaking rule.
	The algorithm then selects candidate
	nodes whose polynomial evaluation
	vectors increase the numerical rank
	until the full 21-dimensional space
	$\Pi_5$ is represented. The remaining
	stencil positions are filled with the
	nearest unused candidates until
	$n_s=29$ nodes have been selected.
	Thus, the resulting stencil need not
	consist of the 29 geometrically nearest
	nodes.
	
	Let a uniform $N\times N$ grid be placed
	on $[0,\pi]^2$, and split the nodes into
	interior and boundary sets $I$ and $B$.
	Homogeneous Dirichlet values are imposed
	strongly:
	\[
	u_i(t)=v_i(t)=0,\qquad i\in B.
	\]
	At each interior center $\bm x_i$,
	the rank-safe KNN procedure selects
	29 nodes from the full grid, including
	boundary nodes when selected.
	The augmented RBF-FD system
	\eqref{eq:blocksys} is solved to obtain
	the local Laplacian weights.
	The global discrete Laplacian is assembled
	row by row, retaining only weights
	associated with interior unknowns.
	Because the prescribed boundary values
	are zero, their contributions vanish.
	This produces
	\begin{equation}
		(\Delta u)|_I\approx \LK\bm u,
		\qquad
		(\Delta v)|_I\approx \LK\bm v.
		\label{eq:knnL}
	\end{equation}
	Since the local stencils are selected
	independently, the resulting differentiation
	weights need not satisfy reciprocal
	coupling conditions. Consequently,
	the assembled matrix $\LK$ is not
	necessarily symmetric.
	This contrasts with the reflected
	complete-shell construction developed
	below, for which global symmetry is
	established analytically.
	
	\subsubsection{Rank-safe KNN RBF-FD algorithm}
	
	Algorithm~\ref{alg:knn} summarizes the rank-safe
	KNN construction used in the numerical experiments.
	The method first selects polynomially independent
	nodes and then completes each stencil according
	to distance ordering. Local differentiation weights
	are assembled into the global interior Laplacian
	with homogeneous Dirichlet conditions imposed
	strongly.
	
	\begin{algorithm}[ht]
		\caption{Rank-safe KNN PHS RBF-FD method}
		\label{alg:knn}
		\begin{algorithmic}[1]
			
			\Require Grid size $N$, stencil size $n_s=29$,
			PHS kernel $\phi(r)=r^{11}$,
			polynomial space $\mathbb P_5$,
			final time $T$, and time step $\Delta t$.
			
			\Ensure Interior Laplacian matrix $\LK$,
			numerical solutions $\bm u(t)$ and $\bm v(t)$,
			and diagnostic quantities.
			
			\State Generate the uniform $N\times N$ grid on
			$[0,\pi]^2$ with spacing $h=\pi/(N-1)$.
			
			\State Partition the grid nodes into interior
			set $I$ and boundary set $B$.
			
			\State Initialize the sparse interior matrix
			$\LK\in\mathbb R^{|I|\times |I|}$.
			
			\For{each interior node $\bm x_i$, $i\in I$}
			
			\State Generate candidate physical grid nodes
			within ten grid spacings of $\bm x_i$,
			retaining offsets $(a,b)$ satisfying
			$a^2+b^2\le 100$.
			
			\State Sort candidates by increasing
			squared Euclidean distance, breaking ties
			lexicographically by grid indices.
			
			\State Initialize the selected stencil
			$\mathcal S_i=\varnothing$ and an empty
			orthonormal polynomial basis $Q$.
			
			\For{each candidate in the sorted list}
			
			\State Evaluate the 21 monomials spanning
			$\mathbb P_5$ at its local offset,
			using scaled coordinates.
			
			\State Orthogonalize the polynomial
			evaluation vector against $Q$.
			
			\If{the residual norm exceeds $10^{-10}$}
			\State Add the candidate to
			$\mathcal S_i$ and normalize
			the residual into $Q$.
			\EndIf
			
			\If{$|Q|=21$}
			\State \textbf{break}
			\EndIf
			
			\EndFor
			
			\If{$|Q|<21$}
			\State Report polynomial rank failure.
			\EndIf
			
			\State Traverse the distance-sorted
			candidates again, adding unused nodes
			until $|\mathcal S_i|=n_s=29$.
			
			\State Construct the local PHS matrix
			\(
			(A_1)_{jk}
			=
			\phi(\|\bm x_j-\bm x_k\|).
			\)
			
			\State Construct the polynomial matrix
			$P$ using a basis of $\mathbb P_5$.
			
			\State Form the Laplacian right-hand side
			\(
			b_j=
			\left.
			\Delta_{\bm x}
			\phi(\|\bm x-\bm x_j\|)
			\right|_{\bm x=\bm x_i},
			c_\alpha=\Delta p_\alpha(\bm x_i).
			\)
			
			\State Solve the augmented system
			\(
			\begin{pmatrix}
				A_1 & P\\
				P^T & 0
			\end{pmatrix}
			\begin{pmatrix}
				\bm w\\
				\bm\lambda
			\end{pmatrix}
			=
			\begin{pmatrix}
				\bm b\\
				\bm c
			\end{pmatrix}.
			\)
			
			\For{each selected stencil node $\bm x_j$}
			\If{$j\in I$}
			\State Insert the corresponding
			weight into $(\LK)_{ij}$.
			\EndIf
			\EndFor
			
			\EndFor
			
			\State Enforce homogeneous Dirichlet conditions
			by setting $u_i=v_i=0$ for all $i\in B$.
			
			\State Form the semidiscrete system
			\(
			\frac{d}{dt}
			\begin{pmatrix}
				\bm u\\
				\bm v
			\end{pmatrix}
			=
			\begin{pmatrix}
				0 & -\LK\\
				\LK & 0
			\end{pmatrix}
			\begin{pmatrix}
				\bm u\\
				\bm v
			\end{pmatrix}.
			\)
			
			\State Integrate the system using the
			time-stepping method specified in the experiment.
			
			\State Compute solution errors, symmetry defect,
			spectral diagnostics, discrete mass, and energy.
			
		\end{algorithmic}
	\end{algorithm}
	
	
	\subsubsection{Nonsymmetry of the assembled KNN operator}
	
	The rank-safe KNN construction generally produces
	a nonsymmetric assembled Laplacian because the
	local stencils and differentiation weights are
	determined independently at each interior node.
	Unlike the reflected RBF-FD construction, no
	reciprocity condition is imposed on the resulting
	matrix entries. Two features of the KNN stencil
	construction contribute to this nonsymmetry.
	
	\begin{enumerate}[leftmargin=2.2em]
		
		\item \textbf{Independent stencil selection and distance ties.}
		Each interior node selects its stencil independently.
		On a Cartesian grid, several candidate nodes may
		have the same distance from a stencil center.
		The deterministic lexicographic tie-breaking rule
		resolves these ambiguities but does not necessarily
		preserve reflection or rotational symmetry.
		Furthermore, the polynomial rank-selection step
		may skip some nearby nodes in favor of more distant
		ones. Consequently, the selected stencils need not
		be geometrically symmetric, even away from the
		boundary.
		
		\item \textbf{Boundary-induced stencil asymmetry.}
		Near $\partial\Omega$, the available physical nodes
		form a truncated neighborhood. The resulting KNN
		stencils generally differ from translations of
		an interior reference stencil. Because the local
		differentiation weights are computed independently,
		the weight assigned to node $\bm x_j$ by the stencil
		centered at $\bm x_i$ need not equal the weight
		assigned to $\bm x_i$ by the stencil centered at
		$\bm x_j$.
		
		Boundary nodes participate in the local weight
		calculations, but their contributions are omitted
		from the global interior matrix because homogeneous
		Dirichlet values are prescribed. This elimination
		does not itself cause the asymmetry; rather, the
		asymmetry arises from the independently constructed
		stencils and their associated differentiation weights.
	\end{enumerate}
	Thus, symmetry of the local augmented systems
	does not generally imply
	\[
	\LK=\LK^T.
	\]
	We quantify the global symmetry defect by
	\begin{equation}
		\eta_{\mathrm{sym}}(L)
		=
		\frac{\|L-L^T\|_F}{\|L\|_F}.
		\label{eq:symdefect}
	\end{equation}
	The rank-safe KNN method is retained throughout
	the paper as an accuracy baseline. No assumption
	is made that a more symmetric stencil necessarily
	produces a smaller approximation error.
	\subsection{Reflected boundary discretization}
	\label{reflectedbond}
	The reflected-boundary construction uses
	the complete-shell reference stencil
	$S_{29}$ defined in \eqref{eq:S29}.
	Its geometry and boundary treatment
	are illustrated in
	Figure~\ref{fig:reflected-stencil-geometry}.
	Panel~\subref{fig:stencil29} shows the
	29-point reference stencil, whose
	reflection symmetry allows a single
	set of reference weights to be used
	throughout the uniform grid.
	Panel~\subref{fig:stencil-reflection}
	illustrates the treatment of stencil
	points lying outside the physical
	domain. Homogeneous Dirichlet conditions
	are imposed through odd reflection,
	so that exterior contributions are
	transferred to reflected interior nodes
	with the appropriate sign.

	
	\subsubsection{Uniform grid and complete reference stencil}
	
	Let
	\[
	h=\frac{\pi}{N-1},\qquad
	x_i=ih,\qquad y_j=jh,
	\qquad 0\le i,j\le N-1.
	\]
	The interior indices satisfy
	$1\le i,j\le N-2$.
	Unlike standard KNN RBF-FD, which constructs
	a separate stencil at each interior node,
	the reflected method uses the fixed,
	reflection-symmetric reference stencil
	$S_{29}$ defined in \eqref{eq:S29}.
	A single set of dimensionless reference
	Laplacian weights $\widehat w_{ab}$,
	$(a,b)\in S_{29}$, is computed using
	the $r^{11}+\Pi_5$ formulation.
	The corresponding physical weights are
	\begin{equation}
		w_{ab}=h^{-2}\widehat w_{ab}.
		\label{eq:weightScaling}
	\end{equation}
	Because the kernel is radial, the polynomial
	space is invariant under coordinate
	reflections, and the Laplacian is reflection
	invariant, uniqueness of the augmented
	system implies
	\begin{equation}
		w_{a,b}
		=w_{-a,b}
		=w_{a,-b}
		=w_{-a,-b}.
		\label{eq:weightReflection}
	\end{equation}
	Coordinate-interchange symmetry also holds:
	\[
	w_{a,b}=w_{b,a}.
	\]
	
	
	\subsubsection{Odd extension for homogeneous Dirichlet data}
	
	Let $U_{ij}$ denote the value of a grid
	function $U$ at the grid point $(x_i,y_j)$,
	with $U_{ij}=0$ at boundary nodes.
	To evaluate the centered stencil near the
	boundary, we define an odd extension
	across $x=0$ and $x=\pi$ by
	\[
	U_{-i,j}=-U_{i,j},
	\qquad
	U_{N-1+i,j}=-U_{N-1-i,j},
	\]
	and analogously across $y=0$ and $y=\pi$. Thus, values at exterior grid points are
	determined by reflection of interior values
	with a change of sign, introducing no
	additional unknowns. A reflection across one boundary contributes a
	factor $-1$; two reflections contribute $+1$.  Points that reflect onto a physical
	boundary have value zero.
	The reflected discrete Laplacian is then defined at every physical interior node by
	\begin{equation}
		(\LD U)_{ij}
		=\sum_{(a,b)\in\mathcal S_R}
		w_{ab}\,\widetilde U_{i+a,j+b},
		\qquad 1\le i,j\le N-2.
		\label{eq:reflectedL}
	\end{equation}
	Exterior stencil points are therefore not additional unknowns.  Their contributions are
	mapped back to physical interior unknowns with the sign dictated by the odd extension.
	For example, at the first interior layer $i=1$, the offset $a=-2$ reaches index $-1$ and
	contributes
	\[
	w_{-2,b}\widetilde U_{-1,j+b}
	=-w_{-2,b}U_{1,j+b}.
	\]
	
	\subsubsection{Matrix interpretation}
	Let $E$ denote the linear odd-extension operator that maps the physical interior vector to
	values on the extended stencil grid, let $\widetilde L$ denote the translation-invariant
	reference-stencil convolution on the extended grid, and let $R$ restrict the result to
	physical interior nodes.  Then
	\begin{equation}
		\LD=R\widetilde L E.
		\label{eq:extensionMatrix}
	\end{equation}
	This identity is useful conceptually, but the implementation need not form ghost unknowns
	or the matrices $E$ and $\widetilde L$ explicitly.  Each exterior stencil contribution is
	reflected and accumulated directly into the appropriate column of $\LD$.

	\begin{proposition}[Symmetry of the reflected Dirichlet operator]
		\label{prop:reflectedSymmetry}
		For an integer radius $R\ge 1$, let
		\(
		M(R)
		=
		1+4R+
		4\sum_{a=1}^{R}
		\left\lfloor\sqrt{R^2-a^2}\right\rfloor,
		\)
		and define the complete-shell reference stencil
		\begin{equation}
			S_{M(R)}
			=
			\left\{
			(a,b)\in\mathbb Z^2:
			a^2+b^2\le R^2
			\right\}.
			\label{eq:completeStencil}
		\end{equation}
		Here, $(a,b)$ denotes an integer offset
		relative to a stencil center.
		In particular, $M(3)=29$, so the principal
		reference stencil is $S_{29}$.
		Assume that the same reference weights
		$\widehat w_{ab}$,
		$(a,b)\in S_{M(R)}$, are used at every
		interior grid node and satisfy
		\begin{equation}
			\widehat w_{a,b}
			=
			\widehat w_{-a,b}
			=
			\widehat w_{a,-b}
			=
			\widehat w_{-a,-b}.
			\label{eq:reflectedWeightSymmetry}
		\end{equation}
		Suppose homogeneous Dirichlet boundary
		conditions are imposed by odd reflection
		across all four sides of the square,
		with a consistent reflected extension.
		Then the assembled interior Laplacian
		$\LD$, defined by \eqref{eq:reflectedL},
		is symmetric:
		\[
		\LD=\LD^T.
		\]
	\end{proposition}

	\begin{proof}
		Let $n=N-1$, so that the grid indices
		range from $0$ to $n$, with interior
		indices $1,\ldots,n-1$ in each
		coordinate direction.
		Extend the reference weights to all of
		$\mathbb Z^2$ by setting
		\[
		K(a,b)=
		\begin{cases}
			w_{a,b}, &(a,b)\in\mathcal S_{M(R)},\\
			0,       &\text{otherwise}.
		\end{cases}
		\]
		Since the stencil is invariant under coordinate
		reflections, \eqref{eq:reflectedWeightSymmetry}
		implies
		\begin{equation}
			K(a,b)=K(-a,b)=K(a,-b)=K(-a,-b).
			\label{eq:kernelEvenSymmetry}
		\end{equation}
		For an interior grid function $U$, its odd
		extension is $2n$-periodic in each coordinate.
		Since the boundary indices are $0$ and $n$,
		the interior indices are $1,\ldots,n-1$.
		Successive odd reflections across both
		boundaries produce a $2n$-periodic extension
		in each coordinate direction.
		In particular, for $1\le p,q\le n-1$,
		the extended values satisfy
		\[
		\widetilde U_{p+2mn,q+2rn}=U_{p,q},
		\]
		\[
		\widetilde U_{-p+2mn,q+2rn}=-U_{p,q},
		\]
		\[
		\widetilde U_{p+2mn,-q+2rn}=-U_{p,q},
		\]
		and
		\[
		\widetilde U_{-p+2mn,-q+2rn}=U_{p,q},
		\]
		for all $m,r\in\mathbb Z$.
		Consequently, the coefficient multiplying
		$U_{p,q}$ in $(\LD U)_{i,j}$ is
		\begin{align}
			(\LD)_{(i,j),(p,q)}
			={}&
			\sum_{m,r\in\mathbb Z}
			\Big[
			K(p-i+2mn,q-j+2rn)
			-K(-p-i+2mn,q-j+2rn)
			\nonumber\\
			&\quad
			-K(p-i+2mn,-q-j+2rn)
			+K(-p-i+2mn,-q-j+2rn)
			\Big].
			\label{eq:reflectedMatrixEntries}
		\end{align}
		All sums are finite because $K$ has
		finite support.
		Define the periodized kernel
		\begin{equation}
			K_n(a,b)
			=
			\sum_{m,r\in\mathbb Z}
			K(a+2mn,b+2rn).
			\label{eq:periodizedKernel}
		\end{equation}
		By \eqref{eq:kernelEvenSymmetry},
		$K_n$ is even in each coordinate:
		\begin{equation}
			K_n(a,b)
			=
			K_n(-a,b)
			=
			K_n(a,-b).
			\label{eq:periodizedSymmetry}
		\end{equation}
		Using \eqref{eq:periodizedKernel},
		equation \eqref{eq:reflectedMatrixEntries}
		becomes
		\begin{align}
			(\LD)_{(i,j),(p,q)}
			={}&
			K_n(p-i,q-j)
			-K_n(p+i,q-j)
			-
			K_n(p-i,q+j)
			+K_n(p+i,q+j).
			\label{eq:reflectedEntryCompact}
		\end{align}
		Now interchange the two interior nodes.
		The reverse matrix entry is
		\begin{align*}
			(\LD)_{(p,q),(i,j)}
			={}&
			K_n(i-p,j-q)
			-K_n(i+p,j-q)
			-
			K_n(i-p,j+q)
			+K_n(i+p,j+q).
		\end{align*}
		By the coordinatewise evenness
		\eqref{eq:periodizedSymmetry},
		each term in this expression equals
		the corresponding term in
		\eqref{eq:reflectedEntryCompact}.
		Therefore,
		\[
		(\LD)_{(i,j),(p,q)}
		=
		(\LD)_{(p,q),(i,j)}
		\]
		for all interior nodes $(i,j)$ and $(p,q)$.
		Hence
		\(\LD=\LD^T,
		\)
		as claimed.
	\end{proof}

	\begin{remark}[Symmetry, definiteness, and accuracy]
		The global symmetry established above has
		two important limitations:
		\begin{enumerate}[label=(\roman*),leftmargin=*]
			
			\item \textbf{Negative definiteness.}
			Symmetry alone does not guarantee
			\(
			\bm z^T\LD\bm z<0
			\text{for all }\bm z\ne0.
			\)
			For the principal 29-point reflected
			operator, negative definiteness is
			established separately through the
			discrete Fourier-symbol analysis in Section\,\ref{sec:fourier}.
			
			\item \textbf{Approximation accuracy.}
			Global symmetry does not necessarily
			imply smaller discretization errors
			than those obtained with standard
			KNN RBF-FD stencils. Accuracy and
			dispersion are assessed independently
			in the numerical experiments.
			
		\end{enumerate}
	\end{remark}

	
	\subsubsection{Reflected-boundary PHS RBF-FD algorithm}
	
	For the principal radius-three stencil and
	the grid sizes considered here, each exterior
	stencil index requires at most one reflection
	per coordinate. Larger stencils on very coarse
	grids would require repeated reflection. Algorithm~\ref{alg:reflected} summarizes the construction
	of the reflected PHS RBF-FD Laplacian for homogeneous
	Dirichlet boundary conditions. Unlike the standard KNN
	method, the reflected construction uses a single
	reference stencil at every interior node. Boundary
	interactions are incorporated through odd reflection,
	rather than by constructing one-sided local stencils.
	
	\begin{algorithm}[h]
		\caption{Reflected-boundary PHS RBF-FD method}
		\label{alg:reflected}
		\begin{algorithmic}[1]

			\Require Grid size $N$,
			reference stencil $S_{29}$ from \eqref{eq:S29},
			PHS kernel $\phi(r)=r^{11}$,
			polynomial space $\mathbb P_5$.

			\Ensure Reflected interior Laplacian matrix
			$\LD\in\mathbb R^{(N-2)^2\times(N-2)^2}$.
			
			\State Set $h=\pi/n, n=N-1$ and generate the
			uniform grid on $[0,\pi]^2$.
			
			\State Set $\mathcal S_R=S_{29}$,
			corresponding to $R=3$ and 29 lattice points.

			\State Construct the augmented PHS system
			\eqref{eq:blocksys} on the reference offsets
			$(a,b)\in\mathcal S_R$.
			
			\State Solve for the dimensionless reference
			Laplacian weights $\widehat w_{ab}$.
			
			\State Scale the weights according to
			\(
			w_{ab}=h^{-2}\widehat w_{ab}.
			\)
			
			\State Verify the reference-weight symmetries
			\(
			w_{a,b}=w_{-a,b}=w_{a,-b}=w_{-a,-b}.
			\)
			
			\State Initialize the sparse interior matrix
			$\LD$ to zero.
			
			\For{each interior grid node $(i,j)$}
			
			\For{each offset $(a,b)\in\mathcal S_R$}

			\State Set
			\(
			(p,q)=(i+a,j+b),
			\quad \sigma=1.
			\)
			
			\If{$p<0$}
			\State $p\gets-p$;
			$\sigma\gets-\sigma$.
			\ElsIf{$p>n$}
			\State $p\gets2n-p$;
			$\sigma\gets-\sigma$.
			\EndIf
			
			\If{$q<0$}
			\State $q\gets-q$;
			$\sigma\gets-\sigma$.
			\ElsIf{$q>n$}
			\State $q\gets2n-q$;
			$\sigma\gets-\sigma$.
			\EndIf
			
			\If{$p=0$ or $p=n$ or
				$q=0$ or $q=n$}
			\State \textbf{continue}
			\EndIf
			
			\State Accumulate
			\(
			(\LD)_{(i,j),(p,q)}
			\gets
			(\LD)_{(i,j),(p,q)}
			+\sigma w_{ab}.
			\)
			\EndFor
			\EndFor
			
			\State Compute the symmetry defect
			\(
			\eta_{\mathrm{sym}}(\LD)
			=
			\frac{\|\LD-\LD^T\|_F}{\|\LD\|_F}.
			\)
			
			\State Perform numerical spectral diagnostics
			to check the sign of the computed eigenvalues.
			
			\State \Return $\LD$.
			
		\end{algorithmic}
	\end{algorithm}
	
	For comparison, we also consider the standard
	fourth-order centered finite-difference (FD4)
	approximation of the Laplacian on the same
	uniform Cartesian grid. Homogeneous Dirichlet
	boundary conditions are imposed using odd
	reflection, yielding a symmetric discrete
	Laplacian.
	
	\section{Structure preservation and spectral analysis}
	\label{sec:semidiscrete}
	
	The assembled spatial operator determines important structural
	properties of the semi-discrete Schr\"odinger system.
	Matrix symmetry leads to conservation of discrete mass and energy,
	while the spectrum controls modal oscillations and numerical
	dispersion. We first establish the conservation laws and
	Hamiltonian formulation, then derive the exact discrete sine
	eigenbasis of the reflected Laplacian. The resulting Fourier
	symbol provides a criterion for negative definiteness and a
	basis for analyzing eigenvalue consistency and phase errors.
	Finally, a modal residual decomposition distinguishes frequency
	error from nonmodal error in the comparison with KNN RBF-FD.
	Throughout, semi-discrete invariants are distinguished from
	their preservation by a time integrator
	\cite{HairerLubichWanner2006,SanzSerna1984,SanzSernaVerwer1986}.
	Conservation alone does not guarantee superior numerical accuracy.
	
	\subsection{Semi-discrete conservation laws and Hamiltonian structure}
	\label{subsec:discreteConservation}
	
	For either spatial operator
	$L\in\{\LK,\LD\}$, the interior
	semi-discrete equations are
	\begin{align}
		\bm u'&=-L\bm v,
		\label{eq:semiu}\\
		\bm v'&=\phantom{-}L\bm u.
		\label{eq:semiv}
	\end{align}
	Introducing
	$\bm y=(\bm u^T,\bm v^T)^T$,
	we obtain
	\begin{equation}
		\bm y'=\mathcal J_h\bm y,
		\qquad
		\mathcal J_h=
		\begin{pmatrix}
			0&-L\\
			L&0
		\end{pmatrix}.
		\label{eq:blockJh}
	\end{equation}
	The algebraic structure of $\mathcal J_h$
	depends directly on the symmetry of $L$.
	Indeed,
	\[
	\mathcal J_h^T=-\mathcal J_h
	\quad\Longleftrightarrow\quad
	L=L^T.
	\]
	Thus, a symmetric spatial Laplacian
	generates a skew-symmetric semi-discrete
	evolution operator, which is the key
	structural condition for the conservation
	results below.
	
	On the uniform Cartesian grid, we use
	the discrete $L^2$ inner product
	\[
	\langle\bm a,\bm b\rangle_h
	=h^2\bm a^T\bm b.
	\]
	The factor $h^2$ accounts for the grid-cell
	area and provides the appropriate scaling
	for the discrete counterparts of the
	continuous mass and Dirichlet energy.
	
	\begin{proposition}[Conservation of discrete mass]
		\label{prop:massConservation}
		Suppose $L=L^T$. Define the discrete mass by
		\begin{equation}
			M_h(t)
			=
			h^2\left(
			\bm u(t)^T\bm u(t)
			+
			\bm v(t)^T\bm v(t)
			\right).
			\label{eq:Mh}
		\end{equation}
		Then $M_h(t)$ is conserved along solutions
		of \eqref{eq:semiu}--\eqref{eq:semiv}.
	\end{proposition}
	
	\begin{proof}
		Differentiating $M_h$ and substituting
		the semi-discrete equations gives
		\begin{align}
			\frac{dM_h}{dt}
			&=
			2h^2\left(
			\bm u^T\bm u'
			+
			\bm v^T\bm v'
			\right)
			=
			2h^2\left(
			-\bm u^TL\bm v
			+
			\bm v^TL\bm u
			\right)
			=
			2h^2\bm u^T(L^T-L)\bm v.
			\label{eq:massDerivative}
		\end{align}
		Since $L=L^T$, it follows that
		$dM_h/dt=0$.
	\end{proof}
	
	\begin{proposition}[Conservation of discrete energy]
		\label{prop:energyConservation}
		Suppose $L=L^T$. Define the discrete
		Dirichlet energy by
		\begin{equation}
			E_h(t)
			=
			\frac{h^2}{2}
			\left(
			\bm u(t)^T(-L)\bm u(t)
			+
			\bm v(t)^T(-L)\bm v(t)
			\right).
			\label{eq:Eh}
		\end{equation}
		Then $E_h(t)$ is conserved along solutions
		of \eqref{eq:semiu}--\eqref{eq:semiv}.
		If, in addition, $L\preceq0$, then
		$E_h(t)\ge0$ for all $t$.
		If $L\prec0$, the energy is strictly
		positive for every nonzero state.
	\end{proposition}
	
	\begin{proof}
		Because $L$ is symmetric and independent
		of time,
		\begin{align}
			\frac{dE_h}{dt}
			&=
			-h^2\left(
			\bm u^TL\bm u'
			+
			\bm v^TL\bm v'
			\right)
			=
			h^2\left(
			\bm u^TL^2\bm v
			-
			\bm v^TL^2\bm u
			\right).
			\label{eq:energyDerivative}
		\end{align}
		Since $L^2$ is symmetric,
		\(
		\bm u^TL^2\bm v
		=
		\bm v^TL^2\bm u.
		\)
		Consequently,
		\(
		\frac{dE_h}{dt}=0.
		\)
		Furthermore, if $L\preceq0$, then $-L$
		is positive semidefinite, so both
		quadratic terms defining $E_h$ are
		nonnegative. If $L\prec0$, then $-L$
		is positive definite, and hence
		$E_h>0$ whenever
		$(\bm u,\bm v)\ne(0,0)$.
	\end{proof}
	
	\begin{remark}[Energy conservation and Hamiltonian structure]
		Conservation and positivity of the discrete
		energy are distinct properties. Symmetry of
		$L$ guarantees conservation of $E_h$,
		whereas nonnegativity additionally requires
		$L\preceq0$, meaning that
		$\bm z^TL\bm z\le0$ for every
		$\bm z\in\mathbb R^n$.
		If $L\prec0$, so that the inequality is
		strict for every nonzero $\bm z$, then
		$E_h$ is positive for every nonzero state.
		
		When $L=L^T$, the semi-discrete system
		also admits the canonical Hamiltonian
		representation
		\begin{equation}
			\bm y'=J\nabla H_h(\bm y),
			\qquad
			J=
			\begin{pmatrix}
				0&I\\
				-I&0
			\end{pmatrix},
			\qquad
			H_h(\bm y)=h^{-2}E_h(\bm y).
			\label{eq:canonicalHamiltonian}
		\end{equation}
		Indeed, by \eqref{eq:Eh},
		\[
		\nabla H_h(\bm y)
		=
		\begin{pmatrix}
			-L&0\\
			0&-L
		\end{pmatrix}\bm y,
		\]
		so that
		$J\nabla H_h(\bm y)=\mathcal J_h\bm y$.
		Since $H_h$ differs from $E_h$ only
		by the constant factor $h^{-2}$,
		their conservation is equivalent.
	\end{remark}
	
	The discrete mass $M_h$ and energy $E_h$
	are counterparts of the continuous
	invariants defined in
	\eqref{eq:continuousMass} and
	\eqref{eq:energy}, respectively.
	For the reflected construction,
	Proposition~\ref{prop:reflectedSymmetry}
	guarantees symmetry of the assembled
	Laplacian and hence exact semi-discrete
	conservation of both quantities.
	In contrast, the standard KNN operator
	is generally nonsymmetric, so these
	conservation properties are not guaranteed.
	In particular, \eqref{eq:massDerivative}
	shows how the antisymmetric part of $L$
	contributes to variations in discrete mass.

	\subsection{Discrete sine diagonalization of the reflected Laplacian}
	The reflection symmetry of the reference stencil
	gives a stronger conclusion than matrix symmetry:
	the entire reflected operator admits an explicit
	discrete sine eigenbasis. Recall that $n=N-1$
	denotes the number of grid intervals in each
	coordinate direction,
	so that $h=\pi/n$ and the interior
	grid indices are $1,\ldots,n-1$. For a complete-shell stencil of integer radius $R$, write
	$S_{M(R)}=\{(a,b)\in\mathbb Z^2:a^2+b^2\le R^2\}$,
	where $M(R)=|S_{M(R)}|$; in particular,
	$S_{M(3)}=S_{29}$. Let $\widehat w_{ab}$
	denote the dimensionless reference weights in
	\eqref{eq:weightScaling}. The odd extension is
	understood as the $2n$-periodic continuation
	in each coordinate, odd about the boundary
	grid indices $0$ and $n$.
	
	\begin{theorem}[Discrete sine diagonalization of the reflected Laplacian]
		\label{thm:dst}
		Assume that the same finite offset set $S_{M(R)}$ and reference
		weights are used at every interior grid point, that $S_{M(R)}$
		is invariant under independent coordinate sign changes, and that
		$\widehat w_{ab}=\widehat w_{-a,b}=\widehat w_{a,-b}$.
		Suppose boundary-crossing stencil values are evaluated by the
		$2n$-periodic odd extension described above.  Then, for
		$1\le k,\ell\le n-1$, the interior grid vectors
		\[
		\Phi_{k\ell}(i,j)=
		\sin\!\left(\frac{k\pi i}{n}\right)
		\sin\!\left(\frac{\ell\pi j}{n}\right),
		\qquad 1\le i,j\le n-1,
		\]
		are eigenvectors of $\LD$:
		\begin{equation}
			\LD\Phi_{k\ell}=-\omega^h_{k\ell}\Phi_{k\ell},
			\qquad
			\omega^h_{k\ell}=
			-\frac1{h^2}\sum_{(a,b)\in S_{M(R)}}
			\widehat w_{ab}
			\cos(kah)\cos(\ell bh).
			\label{eq:discreteSymbol}
		\end{equation}
		In particular, $\LD$ is real symmetric and orthogonally diagonalizable
		by the two-dimensional discrete sine transform.
	\end{theorem}
	
	\begin{proof}
		The odd $2n$-periodic extension of $\Phi_{k\ell}$ coincides at
		every integer pair with the unrestricted product of sine functions:
		the identities $\sin(k(-i)h)=-\sin(kih)$ and
		$\sin(k(2n-i)h)=-\sin(kih)$ hold for integer $k$,
		and similarly in the second coordinate.  Thus, for any interior
		$(i,j)$, the reflected convolution formula becomes
		\[
		(\LD\Phi_{k\ell})_{ij}
		=\frac1{h^2}\sum_{(a,b)\in S_{M(R)}}
		\widehat w_{ab}\,
		\sin(k(i+a)h)\sin(\ell(j+b)h).
		\]
		Expand both shifted sines using the addition formulas.
		In the resulting four-term expansion, terms containing
		$\sin(kah)$ vanish when summed over $a$ by the symmetry
		$\widehat w_{ab}=\widehat w_{-a,b}$; terms containing
		$\sin(\ell bh)$ vanish by the corresponding symmetry in $b$.
		The surviving term is
		\[
		(\LD\Phi_{k\ell})_{ij}
		=\Phi_{k\ell}(i,j)\frac1{h^2}
		\sum_{(a,b)\in S_{M(R)}}
		\widehat w_{ab}\cos(kah)\cos(\ell bh),
		\]
		which is \eqref{eq:discreteSymbol}.
		The $(n-1)^2$ sine-product vectors form an orthogonal basis
		of the interior grid space.  Therefore, if $Q$ contains their
		normalized versions as columns, then
		$\LD=-Q\operatorname{diag}(\omega^h_{k\ell})Q^T$.
		The diagonal entries are real, so $\LD=\LD^T$.
	\end{proof}

	\paragraph{Spectral implications.}
	If $L\bm\varphi=\lambda\bm\varphi$, then the corresponding eigenvalues of
	$\mathcal J_h$ are $\pm i\lambda$.  In particular, a real spectrum for $L$ gives a purely
	imaginary spectrum for the block generator.  Symmetry of $L$ is a sufficient condition for
	reality of its spectrum; negative semidefiniteness additionally gives the correct sign for
	a discrete Dirichlet Laplacian.
	
	For the KNN operator, which may be nonsymmetric, we therefore examine both the matrix
	symmetry defect \eqref{eq:symdefect} and the imaginary parts of the eigenvalues of $L$.
	For the reflected operator, the symmetry diagnostic serves as an implementation check and
	should be near floating-point roundoff when the reflected symmetry assumptions
	are realized numerically.

	\subsection{Fourier-symbol analysis and negative definiteness}
	\label{sec:fourier}
	
	Fourier-symbol analysis provides a convenient way to study
	the spectral properties of a discrete differential operator.
	The Fourier symbol describes how the operator acts on
	individual Fourier modes and determines the corresponding
	discrete eigenvalues. For the reflected RBF-FD Laplacian,
	the discrete sine functions form an exact eigenbasis.
	Consequently, its Fourier symbol can be used to establish
	negative definiteness of the assembled operator and to
	analyze numerical dispersion, which measures the difference
	between the discrete and continuous oscillation frequencies.
	
	We now specialize the discrete sine eigenvalue formula
	established in Theorem~\ref{thm:dst} to the
	29-point stencil $S_{29}$ defined in \eqref{eq:S29}.
	Throughout this subsection, $w_{ab}$ denotes the
	dimensionless reference RBF-FD weight $\widehat w_{ab}$
	associated with the stencil offset $(a,b)$.
	The discrete Fourier symbol is the weighted sum
	\begin{equation}
		\mathscr D(\theta,\varphi)
		=
		-\sum_{(a,b)\in S_{29}}
		w_{ab}\cos(a\theta)\cos(b\varphi),
		\label{eq:symbolDefinition}
	\end{equation}
	where $\theta$ and $\varphi$ are dimensionless
	frequencies in the two coordinate directions.
	Evaluating the symbol at
	$\theta=kh$ and $\varphi=\ell h$ gives
	the discrete sine eigenvalue
	\[
	L_{\mathrm{ref}}\Phi_{k\ell}
	=
	-\omega_{k\ell}^h\Phi_{k\ell},
	\qquad
	\omega_{k\ell}^h
	=
	h^{-2}\mathscr D(kh,\ell h).
	\]
	Thus, the sign of the Fourier symbol determines
	the sign of the discrete Laplacian eigenvalues,
	while its dependence on the frequencies
	determines the numerical dispersion.
	
	\subsubsection{Exact symbol representation and spectral sign}
	Set $t=\sin^2(\theta/2)$, $s=\sin^2(\varphi/2)$. For either of the two analyzed reflected kernels, the complete-shell geometry, dihedral symmetry, and degree-four polynomial reproduction imply the exact symbol
	\begin{equation}\label{eq:symbol}
		\mathscr D(\theta,\varphi)=-\widehat L(\theta,\varphi)=4(t+s)+\tfrac43(t^2+s^2)+a(t^3+s^3)-b(t^2s+ts^2)+ct^2s^2.
	\end{equation}
	\begin{proposition}[Exact symbol structure]
		The dihedrally symmetric $S_{29}$ weights reproducing $\Pi_4$
		have the exact symbol \eqref{eq:symbol}. The reason is that the offset support
		$a^2+b^2\le9$ permits $t^3,s^3,t^2s,ts^2,t^2s^2$, but excludes
		$t^3s,ts^3$ and higher mixed terms. Polynomial reproduction fixes the
		coefficients of the linear and quadratic terms. Hence no unverified
		numerical cancellation is needed to remove other terms.
	\end{proposition}
	\begin{proposition}[Spectral sign criterion]
		\label{prop:symbolSignCriterion}
		Under the hypotheses of Theorem~\ref{thm:dst}, the reflected
		operator is strictly negative definite if and only if
		\(\mathscr D(kh,\ell h)>0\) for all
		\(1\le k,\ell\le N-2\). In particular, positivity of
		\(\mathscr D(\theta,\varphi)\) for all
		\((\theta,\varphi)\in(0,\pi)^2\) is sufficient on every grid.
	\end{proposition}
	\begin{proof}
		Theorem~\ref{thm:dst} yields the orthogonal decomposition
		\(L_{\mathrm{ref}}=-Q\operatorname{diag}(h^{-2}
		\mathscr D(kh,\ell h))Q^T\). The assertion follows from the
		characterization of negative definiteness by eigenvalue signs.
	\end{proof}
	
	\begin{proposition}[Coefficient reduction for the 29-point stencil]
		\label{prop:principalCoefficientReduction}
		Let $w_{30}=w_{3,0}$, $w_{21}=w_{2,1}$, and
		$w_{22}=w_{2,2}$. For dihedrally symmetric weights,
		the coefficients of \eqref{eq:symbol} satisfy
		\begin{equation}\label{eq:principalCoefficientWeights}
			a=64w_{30},\qquad b=-64w_{21}-256w_{22},
			\qquad c=-256w_{22}.
		\end{equation}
	\end{proposition}
	\begin{proof}
		Expand the trigonometric symbol using
		$\cos(2\theta)=1-8t+8t^2$ and
		$\cos(3\theta)=1-18t+48t^2-32t^3$, and the
		analogous formulas in $s$. The $(\pm3,0)$ and $(0,\pm3)$
		offsets determine the pure cubic coefficient, while the
		$(\pm2,\pm1)$, $(\pm1,\pm2)$, and $(\pm2,\pm2)$ offsets
		determine the mixed cubic and quartic coefficients.
	\end{proof}
	
	\subsubsection{All-grid negative definiteness and coefficient verification}
	\begin{theorem}[A sufficient all-grid stability criterion]
		Suppose the 29-point weights are uniquely determined, dihedrally symmetric, and reproduce $\Pi_4$, and that the coefficients of \eqref{eq:symbol} satisfy
		\[
		a\ge0,\qquad b\ge0,\qquad c\ge0,\qquad \tfrac43+a-b\ge0.
		\]
		Then the odd-reflected Dirichlet matrix is symmetric and strictly negative definite for every grid with a nonempty interior. The symmetry and discrete sine eigenbasis follow from Theorem~\ref{thm:dst}; the coefficient conditions establish the required eigenvalue signs.
	\end{theorem}
	\begin{proof}
		For $0\le t,s\le1$, $t^2s+ts^2\le t^3+s^3$ and $t^3+s^3\le t^2+s^2$. Therefore
		\begin{align*}
			-\widehat L&\ge4(t+s)+\tfrac43(t^2+s^2)+(a-b)(t^3+s^3)+ct^2s^2\\
			&\ge4(t+s)+\bigl[\tfrac43+\min\{0,a-b\}\bigr](t^2+s^2)+ct^2s^2>0
		\end{align*}
		for every nonzero sine frequency. The exact sine eigenbasis yields strict negative definiteness. The matrix is symmetric because the dihedrally symmetric convolution and odd reflection preserve reciprocity.
	\end{proof}
	\begin{remark}
		The coefficient criterion is sufficient, not necessary. A numerical basis sweep alone does not prove all-grid stability for the other PHS kernels.
	\end{remark}
	\begin{table}[H]\centering\small
		\caption{Fourier-symbol coefficients for the two analyzed reflected kernels (rounded values). Interval-arithmetic scripts and execution logs are retained as reproducibility materials.}
		\begin{tabular}{lrrrr}\toprule Kernel & $a$ & $b$ & $c$ & $\frac43+a-b$ \\\midrule
			$r^7+\Pi_4$ & 1.9378688733 & 2.4610075589 & 9.6239719934 & 0.8101946477\\
			$r^{11}+\Pi_5$ & 0.8723150311 & 0.3333003630 & 4.7717416029 & 1.8723480015\\
			\bottomrule\end{tabular}\end{table}
	For $r^7$ the positive margin is approximately $0.81019465$; for $r^{11}$, $a-b\approx0.53901467>0$. Both satisfy the theorem's sufficient conditions. The reported coefficient values agree with the accompanying interval-arithmetic
	scripts for the two kernels. Those scripts use outward-rounded interval
	operations and pivot checks; the resulting enclosures are conditional on the
	correctness of the interval library and implementation. Both supplied interval-arithmetic programs were rerun successfully on
	November 14, 2026, and their certified sign checks returned true.
	These runs reproduce the supplied computational certificates, but do not
	constitute an independent implementation or formal verification of the
	interval-arithmetic library. The coefficient-sign theorem itself does
	not depend on numerical rounding.

	\paragraph{Coefficient verification protocol.}
	\label{subsec:symbolVerification}
	
	The Fourier-symbol analysis can be verified directly from
	the dimensionless reference weights $w_{ab}$.
	First, verify the polynomial reproduction conditions
	and the dihedral symmetry identities
	\[
	w_{a,b}=w_{-a,b}=w_{a,-b}=w_{b,a}.
	\]
	Next, construct the trigonometric symbol
	\begin{equation}
		\widehat L(\theta,\varphi)
		=
		\sum_{(a,b)\in S_{29}}
		w_{ab}\cos(a\theta)\cos(b\varphi).
		\label{eq:symbolVerification}
	\end{equation}
	Introduce the variables
	\[
	t=\sin^2\left(\frac{\theta}{2}\right),
	\qquad
	s=\sin^2\left(\frac{\varphi}{2}\right).
	\]
	The resulting expression is a polynomial in $t$ and $s$.
	The support geometry of $S_{29}$ determines which
	coefficients vanish identically, while polynomial
	reproduction determines the lower-order coefficients.
	For numerical verification of the coefficient inequalities, the remaining coefficients
	must be enclosed using validated arithmetic with documented rounding
	and verified nonzero pivots. High-precision numerical
	residuals alone do not constitute certified interval
	bounds.
	
	\noindent Finally, verify positivity of
	$-\widehat L(\theta,\varphi)$ for all nonzero
	frequencies by applying the sufficient coefficient
	conditions established in the stability theorem
	of Section~\ref{sec:fourier}.
	This procedure separates the algebraic derivation
	of the symbol from the numerical verification  of
	its coefficients and provides a reproducible
	verification framework for the reflected operator.

	The preceding sign conditions establish spectral stability for
	the specified reflected stencils. We next quantify the accuracy
	of their eigenvalues and the resulting phase errors.
	
	\subsection{Eigenvalue consistency and numerical dispersion}
	\label{subsec:eigenvalueConvergence}
	
	We now quantify the difference between the discrete frequencies
	$\omega^h_{k\ell}$ in \eqref{eq:discreteSymbol} and the exact Dirichlet
	frequencies $\omega_{k\ell}=k^2+\ell^2$.  Throughout this subsection,
	$k$ and $\ell$ are fixed positive integers, the reference offset set
	$S_{M(R)}$ is finite and independent of $h$, and the dimensionless
	weights $\widehat w_{ab}$ are computed on that fixed reference stencil.
	In particular, these weights are independent of $h$ under the natural
	$h^{-2}$ scaling of the Laplacian.  Define the stencil moments
	\begin{equation}
		M_{pq}=\sum_{(a,b)\in S_{M(R)}}\widehat w_{ab}a^pb^q,
		\qquad p,q\geq0.
		\label{eq:stencilMoments}
	\end{equation}

	\begin{lemma}[Moment identities from polynomial reproduction]
		\label{lem:momentIdentities}
		Suppose the reference PHS RBF-FD weights
		reproduce the Laplacian at the origin
		for every polynomial in $\Pi_d(\mathbb R^2)$,
		where $d\ge 2$.
		
		With the moments defined in \eqref{eq:stencilMoments},
		polynomial reproduction implies the following identities.
		Then, for $p+q\le d$,
		\[
		M_{pq}
		=
		\begin{cases}
			2, & (p,q)=(2,0)\text{ or }(0,2),\\
			0, & \text{otherwise}.
		\end{cases}
		\]
		In particular,
		\begin{equation}
			M_{00}=0,\qquad
			M_{20}=M_{02}=2,\qquad
			M_{11}=0.
			\label{eq:momentIdentities}
		\end{equation}
	\end{lemma}
	
	\begin{proof}
		The discrete moment $M_{pq}$ is the
		weighted sum of the monomial $a^pb^q$
		over all stencil offsets $(a,b)\in S_{M(R)}$.
		Thus, each moment describes how the
		RBF-FD differentiation weights act on
		a particular monomial. The term
		``moment'' refers to this weighted sum
		of powers of the stencil coordinates.
		
		By polynomial reproduction, the RBF-FD
		approximation must give the exact
		Laplacian at the stencil center for
		every polynomial of total degree at
		most $d$. In dimensionless reference
		coordinates, this means
		\[
		\sum_{(a,b)\in S_{M(R)}}
		\widehat w_{ab}P(a,b)
		=
		\Delta P(0,0),
		\qquad P\in\Pi_d(\mathbb R^2).
		\]
		To determine the moments, choose the
		monomial $P(x,y)=x^py^q$, where
		$p+q\le d$. Substituting this polynomial
		into the reproduction identity gives
		\[
		\underbrace{
			\sum_{(a,b)\in S_{M(R)}}
			\widehat w_{ab}a^pb^q
		}_{M_{pq}}
		=
		\left.
		\Delta(x^py^q)
		\right|_{(0,0)}.
		\]
		Therefore, each discrete moment equals
		the exact Laplacian of its corresponding
		monomial evaluated at the stencil center.
		Applying the Laplacian yields
		\[
		\Delta(x^py^q)
		=
		p(p-1)x^{p-2}y^q
		+
		q(q-1)x^py^{q-2},
		\]
		where terms with zero coefficients
		are understood to vanish.
		When evaluated at $(0,0)$, both terms
		vanish except in two cases:
		\[
		\Delta(x^2)\big|_{(0,0)}=2,
		\qquad
		\Delta(y^2)\big|_{(0,0)}=2.
		\]
		Consequently,
		\[
		M_{20}=M_{02}=2,
		\]
		while all other moments with
		$p+q\le d$ are zero.
		In particular, $M_{00}=0$ expresses
		that the discrete Laplacian annihilates
		constant functions, and $M_{11}=0$
		expresses that it reproduces the
		zero Laplacian of the mixed polynomial
		$xy$. This proves the stated identities.
	\end{proof}
	
	For the two formulations considered here, polynomial reproduction holds through degrees four and five, respectively. Combined with reflection symmetry, these identities eliminate all Taylor-expansion error terms of total degree below six, yielding fourth-order spatial consistency.

	The moment identities provide a direct
	connection between polynomial reproduction
	and the accuracy of the discrete eigenvalues.
	For the continuous Laplacian, the sine mode
	$\sin(kx)\sin(\ell y)$ has eigenvalue
	$-(k^2+\ell^2)$. The reflected discrete
	Laplacian has the same sine eigenvectors,
	but its eigenvalues generally differ from
	their continuous counterparts. We now
	use a Taylor expansion of the Fourier
	symbol to quantify this difference.
	The resulting estimate identifies the
	moments responsible for the leading
	frequency error and establishes
	fixed-mode eigenvalue consistency.
	
	\begin{theorem}[Discrete eigenvalue expansion and consistency]
		\label{thm:eigenvalueConsistency}
		Assume the hypotheses of Theorem~\ref{thm:dst} and
		Lemma~\ref{lem:momentIdentities}.  For each fixed mode $(k,\ell)$,
		\begin{equation}
			\displaystyle
			\omega^h_{k\ell}
			=k^2+\ell^2
			-\frac{h^2}{24}\bigl(k^4M_{40}+6k^2\ell^2M_{22}
			+\ell^4M_{04}\bigr)+\mathcal R_{k\ell}(h),
			\label{eq:eigenvalueExpansion}
		\end{equation}
		where the remainder satisfies
		\begin{equation}
			|\mathcal R_{k\ell}(h)|\leq
			\frac{h^4}{720}\sum_{(a,b)\in S_{M(R)}}
			|\widehat w_{ab}|\bigl(|k a|+|\ell b|\bigr)^6.
			\label{eq:eigenvalueRemainder}
		\end{equation}
		Consequently there is a constant $C_{k\ell}$ independent of $h$
		such that
		\begin{equation}
			|\omega^h_{k\ell}-(k^2+\ell^2)|\leq C_{k\ell}h^2
			\qquad (0<h\leq h_0),
			\label{eq:eigenvalueConvergence}
		\end{equation}
		for any fixed $h_0>0$.
	\end{theorem}
	\begin{proof}
		By independent sign-reflection symmetry of the weights,
		\[
		\sum_{a,b}\widehat w_{ab}\cos(kah)\cos(\ell bh)
		=\sum_{a,b}\widehat w_{ab}\cos\bigl(h(ka+\ell b)\bigr).
		\]
		Indeed, the sine--sine term in the cosine addition formula
		cancels when summed over either sign-reflected coordinate.
		Taylor's theorem with the degree-five polynomial for cosine gives,
		for every real $z$,
		\[
		\cos z=1-\frac{z^2}{2}+\frac{z^4}{24}+r_6(z),
		\qquad |r_6(z)|\leq\frac{|z|^6}{720}.
		\]
		Insert $z=h(ka+\ell b)$ into the symbol formula
		\eqref{eq:discreteSymbol}.  The constant term vanishes by $M_{00}=0$.
		This reflects the fact that the discrete
		Laplacian reproduces the zero Laplacian
		of constant functions. The quadratic
		terms determine the leading behavior
		of the discrete eigenvalue and must
		recover the continuous frequency
		$k^2+\ell^2$.
		The quadratic contribution is
		\[
		\frac12\sum_{a,b}\widehat w_{ab}(ka+\ell b)^2
		=\frac12(k^2M_{20}+2k\ell M_{11}+\ell^2M_{02})
		=k^2+\ell^2.
		\]
		Expanding the fourth power, all terms odd in either $a$ or $b$
		vanish by reflection symmetry.  Thus
		\[
		\sum_{a,b}\widehat w_{ab}(ka+\ell b)^4
		=k^4M_{40}+6k^2\ell^2M_{22}+\ell^4M_{04}.
		\]
		Thus, the fourth-order moments determine
		the leading possible discrepancy between
		the discrete and continuous eigenvalues.
		If these moments vanish through polynomial
		reproduction, the $O(h^2)$ correction
		disappears, leaving an error of order
		$h^4$ or smaller.
		The remainder after division by $h^2$ is bounded by
		\eqref{eq:eigenvalueRemainder}.  This proves
		\eqref{eq:eigenvalueExpansion}.  Because the reference stencil is
		finite and its dimensionless weights do not depend on $h$, the
		fourth- and sixth-order weighted sums are finite constants.
		For $h\leq h_0$, the $h^4$ remainder is bounded by a constant
		times $h^2$, establishing \eqref{eq:eigenvalueConvergence}.
	\end{proof}

	\begin{remark}[Leading-order dispersion and possible cancellation]
		\label{rem:leadingDispersion}
		The eigenvalue expansion separates the
		continuous frequency $k^2+\ell^2$ from
		the discretization error. In general,
		the fourth-order moments
		\(M_{40},  M_{22}, M_{04}\)
		determine the coefficient of the
		$O(h^2)$ error term.
		If
		\[
		D_{k\ell}
		=
		k^4M_{40}
		+6k^2\ell^2M_{22}
		+\ell^4M_{04}
		\neq 0,
		\]
		then
		\[
		\omega^h_{k\ell}-(k^2+\ell^2)
		=
		-\frac{h^2}{24}D_{k\ell}
		+O(h^4).
		\]
		If $D_{k\ell}=0$, the leading
		$O(h^2)$ term cancels and the error
		is $O(h^4)$. Such cancellation may
		occur for an individual mode or,
		as in the present polynomially
		augmented formulations, for all
		fixed modes through polynomial
		reproduction.
		The estimates concern fixed
		wavenumbers $k,\ell$ as $h\to0$.
		They do not imply uniform accuracy
		for high-frequency modes whose
		wavenumbers increase with $h^{-1}$.
	\end{remark}

	\begin{corollary}[Fourth-order fixed-mode eigenvalue consistency]
		\label{cor:p5EigenvalueConsistency}
		Under the hypotheses of
		Theorem~\ref{thm:eigenvalueConsistency},
		suppose in addition that the reference
		weights reproduce every polynomial in
		$\Pi_4(\mathbb R^2)$. Then, for each
		fixed sine mode $(k,\ell)$,
		\[
		\bigl|
		\omega^h_{k\ell}-(k^2+\ell^2)
		\bigr|
		\le
		\frac{h^4}{720}
		\sum_{(a,b)\in S_{M(R)}}
		|\widehat w_{ab}|
		\bigl(|ka|+|\ell b|\bigr)^6.
		\]
		In particular, the eigenvalue error
		is $O(h^4)$ as $h\to0$.
	\end{corollary}
	
	\begin{proof}
		Polynomial reproduction through degree
		four implies
		\[
		M_{40}=M_{22}=M_{04}=0,
		\]
		because the Laplacians of $x^4$,
		$x^2y^2$, and $y^4$ vanish at the
		stencil center.
		
		Therefore, the fourth-order moment
		contribution in
		\eqref{eq:eigenvalueExpansion}
		vanishes identically. The eigenvalue
		error is determined by the remainder,
		which satisfies
		\eqref{eq:eigenvalueRemainder}.
		Since the reference stencil is finite
		and its dimensionless weights are
		independent of $h$, the right-hand side
		is bounded by $C_{k\ell}h^4$ for
		each fixed mode $(k,\ell)$.
	\end{proof}

	The corollary establishes fourth-order
	convergence of the discrete eigenvalues
	for each fixed sine mode. This is a
	spatial consistency result and does not,
	by itself, establish fourth-order
	convergence of the fully discrete
	solution. The latter also depends on
	the temporal discretization, the time
	step, and the regularity and spectral
	content of the initial data.
	
	
	\subsubsection{Sixth-order moments and directional dispersion}
	
	The following expansion applies to both analyzed PHS
	kernels, although the numerical moments and favorable
	directions depend on the selected kernel.
	For comparison, we also include the standard
	fourth-order centered finite-difference (FD4)
	discretization of the Laplacian, with homogeneous
	Dirichlet conditions imposed by odd reflection.
	
	\paragraph{Sixth-order moments and favorable directions.}
	Define the sixth moments of the dimensionless weights by
	\[M_{60}=\sum_{(a,b)\in S_{29}}w_{ab}a^6,\qquad
	M_{42}=\sum_{(a,b)\in S_{29}}w_{ab}a^4b^2.\]
	For the $r^7+\Pi_4$ weights on this stencil,
	the numerically computed sixth-order moments are
	\begin{equation}
		M_{60}\approx 13.801025,
		\qquad
		M_{42}\approx -1.845756.
		\label{eq:sixthMoments}
	\end{equation}
	
	\begin{proposition}[Fixed-mode dispersion expansion]\label{prop:dispersion}
		For fixed positive integers $k,\ell$ as $h\to0$,
		\begin{align}
			\omega_{29}^h-(k^2+\ell^2)
			&=\frac{h^4}{720}\left[M_{60}(k^6+\ell^6)
			+15M_{42}(k^4\ell^2+k^2\ell^4)\right]+O(h^6),\label{eq:phsexp}\\
			\omega_{\rm FD4}^h-(k^2+\ell^2)
			&=-\frac{h^4}{90}(k^6+\ell^6)+O(h^6).\label{eq:fdexp}
		\end{align}
	\end{proposition}
	\begin{proof}
		Expand the cosine symbol about $(0,0)$. Polynomial reproduction
		through degree four cancels the fourth-order symbol term;
		the first possible frequency error arises from sixth-order
		terms in the cosine expansion. The sixth-order terms are expressed in terms
		of the stencil moments $M_{60}$ and $M_{42}$
		defined above. The FD4 expansion follows from
		$4\sin^2(\theta/2)+\frac43\sin^4(\theta/2)
		=\theta^2-\theta^6/90+O(\theta^8)$.
	\end{proof}
	
	Let $z=k^2\ell^2/(k^2+\ell^2)^2$. For the $r^7+\Pi_4$ method, the PHS29 leading absolute
	frequency error is smaller than that of FD4 exactly when
	\begin{equation}\label{eq:directioncondition}
		\left|M_{60}+(15M_{42}-3M_{60})z\right|<8(1-3z).
	\end{equation}
	For the $r^7+\Pi_4$ method, the endpoints of this strict interval are
	\[
	z_- =\frac{8-M_{60}}{15M_{42}-3M_{60}+24},\qquad
	z_+=\frac{-8-M_{60}}{15M_{42}-3M_{60}-24},
	\]
	Numerically, $z_-\approx 0.128656$ and
	$z_+\approx 0.234195$.
	Equivalently, writing
	$r=\min(k,\ell)/\max(k,\ell)$,
	the favorable directions are approximately
	characterized by
	\begin{equation}\label{eq:directionratio}
		0.422808 < r < 0.773407.
	\end{equation}
	The comparison is asymptotic for \emph{fixed} modes; it does not
	apply uniformly to modes whose indices grow as $h^{-1}$.
	For a finite grid, set
	\[
	\Delta=a(t^3+s^3)-b(t^2s+ts^2)+ct^2s^2,\qquad
	E=(kh)^2+(\ell h)^2-\left[4(t+s)+\frac43(t^2+s^2)\right]>0.
	\]
	Then the exact necessary-and-sufficient condition for a strictly
	smaller PHS29 absolute frequency error is
	\begin{equation}\label{eq:finitecriterion}
		\quad0<\Delta<2E.
	\end{equation}
	Indeed $h^2(\omega_{29}^h-\omega)=\Delta-E$ and
	$h^2(\omega_{\rm FD4}^h-\omega)=-E$; compare absolute values.

	\paragraph{Dispersion for the principal $r^{11}+\Pi_5$ formulation.}
	The same sixth-moment expansion applies with kernel-dependent moments
	$M_{60}^{(11)}=\sum_{(a,b)\in S_{29}}w_{ab}^{(11)}a^6$ and
	$M_{42}^{(11)}=\sum_{(a,b)\in S_{29}}w_{ab}^{(11)}a^4b^2$:
	\begin{equation}\label{eq:disp11}
		\omega_{11}^h-(k^2+\ell^2)=\frac{h^4}{720}
		\left[M_{60}^{(11)}(k^6+\ell^6)+15M_{42}^{(11)}(k^4\ell^2+k^2\ell^4)\right]+O(h^6).
	\end{equation}
	For the principal $r^{11}+\Pi_5$ weights, the corresponding
	numerically reconstructed moments are
	\[
	M_{60}^{(11)}\approx1.813544100277468,\qquad
	M_{42}^{(11)}\approx-0.2499752722529451.
	\]
	The directional interval \eqref{eq:directionratio} is specific to $r^7$;
	it must not be attributed to $r^{11}$. For the latter, the exact finite-grid
	symbol and the reported numerical tests provide direct evidence of improved
	accuracy for the selected modes. The theorem proves definiteness, not
	universal dispersion superiority.
	
	\subsection{Modal evolution and phase-error estimates}
	
	\begin{theorem}[Exact semidiscrete modal evolution and phase error]
		\label{thm:phase}
		Under the assumptions of Theorem~\ref{thm:dst}, let the initial
		interior vectors be $\bm u(0)=\Phi_{k\ell}$ and $\bm v(0)=\bm0$.
		Then the exact solution of the \emph{semidiscrete} system is
		\begin{equation}
			\bm u_h(t)=\cos(\omega^h_{k\ell}t)\Phi_{k\ell},
			\qquad
			\bm v_h(t)=-\sin(\omega^h_{k\ell}t)\Phi_{k\ell}.
			\label{eq:semidiscreteMode}
		\end{equation}
		Let $\bm u_{\rm ex}(t),\bm v_{\rm ex}(t)$ denote the samples of
		the continuous exact solution with frequency
		$\omega_{k\ell}=k^2+\ell^2$.  Then, for $t\ge0$,
		\begin{equation}
			\left(
			\|\bm u_h(t)-\bm u_{\rm ex}(t)\|_2^2+
			\|\bm v_h(t)-\bm v_{\rm ex}(t)\|_2^2
			\right)^{1/2}
			\le t\,|\omega^h_{k\ell}-\omega_{k\ell}|\,
			\|\Phi_{k\ell}\|_2.
			\label{eq:phaseBound}
		\end{equation}
	\end{theorem}
	\begin{proof}
		Substitution of $\LD\Phi_{k\ell}=-\omega^h_{k\ell}\Phi_{k\ell}$
		into the semidiscrete equations yields the oscillator system
		$a'=\omega^h_{k\ell}b$, $b'=-\omega^h_{k\ell}a$,
		with $a(0)=1$, $b(0)=0$.  This proves
		\eqref{eq:semidiscreteMode}.
		The squared combined error equals
		\[
		\|\Phi_{k\ell}\|_2^2
		\left[
		\bigl(\cos(\omega^h_{k\ell}t)-\cos(\omega_{k\ell}t)\bigr)^2+
		\bigl(\sin(\omega^h_{k\ell}t)-\sin(\omega_{k\ell}t)\bigr)^2
		\right].
		\]
		The bracket equals
		$4\sin^2((\omega^h_{k\ell}-\omega_{k\ell})t/2)$.
		Apply $2|\sin(z/2)|\le |z|$ to obtain \eqref{eq:phaseBound}.
	\end{proof}
	
	\begin{corollary}[Second-order semidiscrete phase convergence]
		\label{cor:phaseConvergence}
		Under the assumptions of Theorem~\ref{thm:eigenvalueConsistency},
		consider the initial sine mode in Theorem~\ref{thm:phase}.
		For each fixed $T>0$, the semidiscrete solution satisfies
		\begin{equation}
			\sup_{0\leq t\leq T}
			\left(
			\|\bm u_h(t)-\bm u_{\rm ex}(t)\|_{2,h}^2+
			\|\bm v_h(t)-\bm v_{\rm ex}(t)\|_{2,h}^2
			\right)^{1/2}
			\leq C_{k\ell}T h^2\|\Phi_{k\ell}\|_{2,h},
			\label{eq:phaseConvergence}
		\end{equation}
		where $\|\bm z\|_{2,h}=h\|\bm z\|_2$ is the
		mesh-weighted Euclidean norm on the two-dimensional interior grid.
		In particular, the combined error is $O(h^2)$ for fixed $k,\ell,T$.
	\end{corollary}
	\begin{proof}
		Multiply the estimate \eqref{eq:phaseBound} by $h$ and use
		\eqref{eq:eigenvalueConvergence}.  The discrete sine orthogonality
		identity gives
		\[
		\|\Phi_{k\ell}\|_{2,h}^2
		=h^2\left(\frac{n}{2}\right)^2=\frac{\pi^2}{4}
		\qquad(1\leq k,\ell\leq n-1),
		\]
		so this norm is bounded independently of grid refinement.
	\end{proof}
	
	\begin{corollary}[Fourth-order fixed-mode frequency and phase convergence]
		\label{cor:fourthOrderPhase}
		Assume the hypotheses of Theorem~\ref{thm:eigenvalueConsistency}
		and suppose the reference weights reproduce the Laplacian
		for every polynomial in $\Pi_4$. In particular, these
		hypotheses hold for the principal 29-point reflected
		construction whenever its augmented local weight system
		is nonsingular. Then for every fixed pair of positive integers
		$(k,\ell)$, the discrete frequency satisfies
		\begin{equation}
			\omega^h_{k\ell}=k^2+\ell^2
			+\frac{h^4}{720}\bigl(k^6M_{60}
			+15k^4\ell^2M_{42}+15k^2\ell^4M_{24}
			+\ell^6M_{06}\bigr)+O(h^6).
			\label{eq:fourthOrderFrequencyExpansion}
		\end{equation}
		In particular,
		$|\omega^h_{k\ell}-(k^2+\ell^2)|\le C^{(4)}_{k\ell}h^4$
		for $0<h\le h_0$, with a constant independent of $h$.
		For the initial sine mode of Theorem~\ref{thm:phase},
		\begin{equation}
			\sup_{0\le t\le T}
			\left(\|\bm u_h(t)-\bm u_{\rm ex}(t)\|_{2,h}^2
			+\|\bm v_h(t)-\bm v_{\rm ex}(t)\|_{2,h}^2\right)^{1/2}
			\le C^{(4)}_{k\ell}T h^4\|\Phi_{k\ell}\|_{2,h}
			\label{eq:fourthOrderPhaseBound}
		\end{equation}
		for every fixed $T>0$.
	\end{corollary}
	\begin{proof}
		Polynomial reproduction through degree four gives
		$M_{40}=M_{22}=M_{04}=0$ in addition to the lower-order
		moment identities. Independent coordinate reflection
		symmetry makes every moment with an odd exponent vanish.
		Expand the cosine in \eqref{eq:discreteSymbol} through degree
		seven, using the eighth-order Taylor remainder. The
		quadratic term yields $k^2+\ell^2$, the fourth-order
		term vanishes, and the sixth-order term yields
		\eqref{eq:fourthOrderFrequencyExpansion}. Since the stencil
		is finite and the dimensionless weights are independent
		of $h$, the remainder after division by $h^2$ is $O(h^6)$
		for fixed $(k,\ell)$. The frequency bound follows.
		Finally, apply \eqref{eq:phaseBound} in the mesh-weighted
		norm and take the supremum over $0\le t\le T$.
	\end{proof}
	
	\begin{remark}[Fixed-mode scope of the fourth-order estimate]
		The estimate in Corollary~\ref{cor:fourthOrderPhase} is
		not uniform for modes whose wavenumbers increase with
		$h^{-1}$. Its sixth-order coefficient is consistent with
		the more detailed Fourier-symbol analysis in
		Section~\ref{sec:fourier}.
	\end{remark}
	
	\begin{remark}[Semi-discrete versus fully discrete conservation and accuracy]
		The conservation and phase estimates above concern the exact
		continuous-time evolution of the semi-discrete system. They do not
		imply exact conservation under an arbitrary time integrator.
		Classical RK4 generally introduces temporal mass and energy errors.
		For a fixed sine mode on a bounded time interval, its temporal error
		is $O(\Delta t^4)$ under the usual smoothness and stability
		conditions. Thus, the second-order fixed-mode estimate gives
		$O(h^2+\Delta t^4)$ for the fully discrete solution; when the
		fourth-order spatial consistency applies, the corresponding estimate
		is $O(h^4+\Delta t^4)$. These are fixed-mode conclusions, not
		uniform-in-mode stability theorems for the full grid operator.
		By contrast, the implicit midpoint method preserves the quadratic
		invariants $M_h$ and $E_h$ when $L=L^T$.
	\end{remark}
	
	\subsection{Modal residuals and comparison with KNN}
	The exact sine eigenvectors allow the frequency error of the reflected
	operator to be separated from the nonmodal residual of a general
	RBF-FD matrix. The following identity is useful when comparing
	reflected and KNN discretizations.
	
	\begin{proposition}[Modal residual decomposition]
		\label{prop:modalResidual}
		For any real matrix $L$ and nonzero sampled
		mode $\bm\phi$, let $\omega=k^2+\ell^2$
		and define
		\[
		\omega_h^*
		=
		-\frac{\bm\phi^T L\bm\phi}
		{\bm\phi^T\bm\phi},
		\qquad
		\bm\eta
		=
		L\bm\phi+\omega_h^*\bm\phi.
		\]
		Then $\bm\phi^T\bm\eta=0$, and
		\begin{equation}
			\|L\bm\phi+\omega\bm\phi\|_2^2
			=
			(\omega-\omega_h^*)^2
			\|\bm\phi\|_2^2
			+
			\|\bm\eta\|_2^2.
			\label{eq:orthogonalResidual}
		\end{equation}
		For $L=\LD$ under
		Theorem~\ref{thm:dst}, and
		$\bm\phi=\Phi_{k\ell}$, one has
		$\omega_h^*=\omega^h_{k\ell}$ and
		$\bm\eta=\bm0$ in exact arithmetic.
	\end{proposition}
	
	\begin{proof}
		By the definition of $\omega_h^*$,
		\[
		\begin{aligned}
			\bm\phi^T\bm\eta
			&=
			\bm\phi^T
			(L\bm\phi+\omega_h^*\bm\phi)=
			\bm\phi^TL\bm\phi
			+\omega_h^*\bm\phi^T\bm\phi=0.
		\end{aligned}
		\]
		Thus, $\bm\eta$ is orthogonal to
		$\bm\phi$.
		
		Next, adding and subtracting
		$\omega_h^*\bm\phi$ gives
		\[
		\begin{aligned}
			L\bm\phi+\omega\bm\phi
			&=
			L\bm\phi+\omega_h^*\bm\phi
			+(\omega-\omega_h^*)\bm\phi=
			\bm\eta
			+(\omega-\omega_h^*)\bm\phi.
		\end{aligned}
		\]
		Since $\bm\eta$ and $\bm\phi$ are
		orthogonal, taking squared Euclidean
		norms yields
		\[
		\begin{aligned}
			\|L\bm\phi+\omega\bm\phi\|_2^2
			&=
			\|\bm\eta\|_2^2
			+(\omega-\omega_h^*)^2
			\|\bm\phi\|_2^2,
		\end{aligned}
		\]
		which proves
		\eqref{eq:orthogonalResidual}.
		
		Finally, Theorem~\ref{thm:dst}
		establishes that
		\[
		\LD\Phi_{k\ell}
		=
		-\omega^h_{k\ell}\Phi_{k\ell}.
		\]
		Consequently,
		\[
		\omega_h^*
		=
		-\frac{
			\Phi_{k\ell}^T\LD\Phi_{k\ell}
		}{
			\Phi_{k\ell}^T\Phi_{k\ell}
		}
		=
		\omega^h_{k\ell},
		\]
		and hence
		\(
		\bm\eta
		=
		\LD\Phi_{k\ell}
		+\omega^h_{k\ell}\Phi_{k\ell}
		=
		\bm0.
		\)
	\end{proof}

	The decomposition \eqref{eq:orthogonalResidual}
	separates the modal residual into a frequency
	error and an eigenvector residual. For the
	reflected operator, the discrete sine modes
	are exact eigenvectors, so the eigenvector
	residual vanishes. The remaining modal
	residual is therefore determined entirely
	by the discrete frequency error
	$\omega-\omega^h_{k\ell}$.
	For a general KNN operator, the sampled
	sine modes need not be eigenvectors,
	and both terms may contribute.
	\section{Numerical experiments}\label{sec:results}
	The experiments address three questions: (i) how the reflected PHS
	formulations compare in accuracy and dispersion with alternative spatial
	operators; (ii) whether the structural properties established in
	Section~\ref{sec:semidiscrete} are reflected in assembled-matrix and
	time-integration diagnostics; and (iii) what accuracy--cost tradeoffs
	arise relative to rank-safe KNN RBF-FD and fourth-order finite differences.
	We distinguish modal propagation, spatial-operator tests, and full
	sparse-matrix time integration throughout. Unless stated otherwise,
	reported errors are numerical observations rather than general convergence
	theorems.

	We consider two numerical examples: the free Schr\"odinger
	equation, for which an exact solution is available, and
	the Schr\"odinger equation with a spatially varying
	potential, which provides a more general test of the
	numerical methods.
	
	The experiments address three questions:
	(i) how the reflected PHS formulations compare in accuracy
	and dispersion with alternative spatial operators;
	(ii) whether the structural properties established in
	Section~\ref{sec:semidiscrete} are reflected in
	assembled-matrix and time-integration diagnostics; and
	(iii) what accuracy--cost tradeoffs arise relative to
	rank-safe KNN RBF-FD and fourth-order finite differences.
	
	We distinguish modal propagation, spatial-operator tests,
	and full sparse-matrix time integration throughout.
	Unless stated otherwise, reported errors are numerical
	observations rather than general convergence theorems.
	
	\subsection{Experimental setup and numerical consistency}
	\label{subsec:experimentalSetup}
	
	\paragraph{Time integration and conservation diagnostics.}
	We employ two time-integration methods with complementary
	objectives. Classical fourth-order Runge--Kutta (RK4) is used
	as the common accuracy baseline for the KNN and reflected
	RBF-FD spatial discretizations. The implicit midpoint method
	is used to investigate the preservation of discrete mass and
	energy. Both integrators are applied to the same spatial
	operators, allowing spatial discretization effects to be
	distinguished from temporal integration effects.
	
	For the semi-discrete Schr\"odinger system
	\[
	\bm y'(t)=\mathcal J_h\bm y(t),
	\]
	the classical RK4 method is fourth-order accurate in time
	but does not generally preserve quadratic invariants exactly.
	Consequently, the observed mass and energy drift under RK4
	may contain contributions from both spatial discretization
	and temporal integration.
	
	The implicit midpoint method is given by
	\begin{equation}
		\left(I-\frac{\Delta t}{2}\mathcal J_h\right)
		\bm y^{n+1}
		=
		\left(I+\frac{\Delta t}{2}\mathcal J_h\right)
		\bm y^n.
		\label{eq:midpointScheme}
	\end{equation}
	When $L=L^T$, the matrix $\mathcal J_h$ is skew-symmetric,
	and the midpoint update is the Cayley transform of a
	skew-symmetric matrix. It is therefore orthogonal and
	preserves discrete mass exactly in exact arithmetic.
	Because the discrete energy is also a quadratic invariant
	of the symmetric semi-discrete system, implicit midpoint
	preserves it as well.
	
	For the variable-potential problem, the corresponding
	discrete Hamiltonian includes the potential matrix.
	When this Hamiltonian is symmetric, the same quadratic
	invariant-preservation property holds.
	
	The midpoint method is second-order accurate in time,
	whereas RK4 is fourth-order. Accordingly, RK4 is used
	primarily for solution-error and convergence comparisons,
	while midpoint provides a complementary test of the
	structure-preserving properties of the spatial operators.
	
	For both integrators, the default time step is selected
	using a bound based on the maximum absolute row sum of
	the spatial operator, with an additional contribution
	from the potential when present. The same time step is
	used when comparing the two integrators at a given
	spatial resolution. This restriction is required for
	the explicit RK4 stability considerations but is not
	an unconditional stability requirement for implicit
	midpoint.
	
	\paragraph{Consistency and boundary treatment.}
	For standard KNN RBF-FD, local consistency depends on
	polynomial reproduction, stencil geometry, and the
	behavior of the differentiation weights. The augmented
	PHS systems enforce exact reproduction of the
	Laplacian on the prescribed polynomial space
	$\Pi_m(\mathbb R^2)$, provided the corresponding
	local systems are nonsingular. In our implementation,
	the solvability of these systems is checked numerically.
	
	For the reflected construction, the same reference
	weights are used throughout the domain, while odd
	extension supplies the values required at
	boundary-crossing stencil positions. This construction
	preserves the symmetry of the assembled discrete
	Laplacian.
	
	Under sufficient smoothness and uniform bounds on
	the appropriately scaled differentiation weights,
	a local Taylor expansion gives an
	$O(h^{m-1})$ truncation-error estimate for
	second derivatives when polynomials through degree
	$m$ are reproduced. Additional symmetry and
	moment cancellation may improve this order.
	
	For the reflected method, the corresponding
	consistency argument also requires sufficient
	smoothness of the odd continuation across the
	boundary. The sine eigenmodes considered in the
	verification experiments admit such a continuation.
	
	These consistency considerations provide a basis for
	the expected convergence behavior but do not constitute
	a general RBF-FD convergence theorem. For the reflected
	method, the local consistency argument additionally
	requires sufficient smoothness of the odd extension
	across the boundary, a condition satisfied by the sine
	eigenmodes used in our verification experiments.
	The convergence rates are examined numerically below.

	\paragraph{Computational environment.}
	The numerical experiments were implemented in C++
	using the Eigen linear algebra library. Classical RK4
	uses explicit evaluations of the semi-discrete right-hand
	side, whereas implicit midpoint requires the solution
	of a sparse linear system at each time step. Since the
	discrete Hamiltonian is time-independent in the experiments,
	the midpoint coefficient matrix is factorized once and
	the factorization is reused throughout the simulation.
	
	Python was used to generate figures and visualize the
	numerical results. All computations were performed on a
	computer equipped with a 13th-generation Intel Core
	i7-1355U processor (1.70 GHz) and 16 GB of RAM, running
	a 64-bit Windows operating system.
	
	\subsection{Kernel and polynomial screening}
	The first diagnostic experiment, \cref{tab2}, evaluates ten kernel--polynomial combinations on the same complete-shell 29-point stencil. Each local system is solved with 65-decimal-digit arithmetic, and the weights are converted to double precision for global computations. The combinations are chosen to examine odd-power and logarithmic PHS families; because polynomial order changes, the table is a screening study rather than a controlled estimate of the effect of the kernel alone. Screening error is measured for the mixed-mode initial condition
	\[
	z_0=\sin x\sin y+0.8\sin(2x)\sin(3y)+0.5\sin x\sin(5y)+0.4\sin(3x)\sin(4y)+0.2\sin(5x)\sin(5y).
	\]
	The exact solution multiplies each sine product by $e^{-i(k^2+\ell^2)T}$. 
	\begin{table}[H]\centering\small
		\caption{Exploratory basis and polynomial-order screening at $N=31$, $T=1$ for the mixed-mode test. The listed polynomial spaces vary with the kernel; the screening is not a controlled one-factor comparison.}
		\label{tab2}
		\begin{tabular}{lrr}
			\toprule
			Kernel and augmentation
			& Relative $L^2$ error
			& RK4 temporal error \\
			\midrule
			$r^3,\mathbb P_2$             & 0.3477   & $5.09\times10^{-8}$ \\
			$r^5,\mathbb P_3$             & 0.1399   & $8.50\times10^{-8}$ \\
			$r^7,\mathbb P_4$             & 0.01264  & $8.65\times10^{-8}$ \\
			$r^9,\mathbb P_4$             & 0.006111 & $8.67\times10^{-8}$ \\
			$r^{11},\mathbb P_5$          & 0.001152 & $8.69\times10^{-8}$ \\
			$r^2\log r,\mathbb P_1$       & 0.2635   & $7.64\times10^{-8}$ \\
			$r^4\log r,\mathbb P_2$       & 0.2629   & $8.20\times10^{-8}$ \\
			$r^6\log r,\mathbb P_3$       & 0.05811  & $8.65\times10^{-8}$ \\
			$r^8\log r,\mathbb P_4$       & 0.009299 & $8.66\times10^{-8}$ \\
			$r^{10}\log r,\mathbb P_5$    & 0.003352 & $8.68\times10^{-8}$ \\
			\bottomrule
		\end{tabular}
	\end{table}
	The screening identifies $r^{11}+\mathbb P_5$ as the leading tested choice on the specified stencil, domain, initial data, and resolution. This does not establish global optimality over all kernels, polynomial degrees, stencil sizes, or PDEs. We retain $r^7+\mathbb P_4$ as a second analytically certified example and FD4 as a computationally economical baseline.
	
	
	\subsection{Modal accuracy and stencil sensitivity}
	
	\paragraph{Modal convergence.}
	For the modal convergence experiments, we use
	$N=21,31,41,61,81,101$ and final times
	$T=1$ and $T=20$.
	Because each reflected matrix is exactly
	diagonalized by the discrete sine basis,
	applying the RK4 stability polynomial to each
	eigenvalue is algebraically equivalent to
	sparse-matrix RK4 for the prescribed sine-mode
	initial data, apart from floating-point roundoff.
	This modal RK4 implementation permits fine
	temporal refinement and separates temporal error
	from spatial dispersion. It is not used to infer
	sparse-matrix wall-clock cost.
	Time steps are chosen below a conservative
	imaginary-axis stability threshold, with a cap
	of $0.0007$.
	The relative $L^2$ error is computed at interior
	nodes against the exact continuous solution.

	The mixed-mode initial condition is the five-mode
	combination introduced in the preceding subsection.
	To examine the directional dispersion behavior
	predicted in Section~3, we also consider the
	two-mode initial condition
	\[
	z_0(x,y)
	=0.8\sin(2x)\sin(3y)
	+0.4\sin(3x)\sin(4y).
	\]
	Both modes satisfy the favorable-direction
	criterion \eqref{eq:directionratio} for
	$r^7+\mathbb P_4$ relative to FD4.
	We refer to this test as the
	\emph{favorable-mode} experiment.

\begin{figure}[pos=htbp]
		\centering
		
		\begin{subfigure}[t]{0.48\textwidth}
			\centering
			\includegraphics[width=\linewidth]{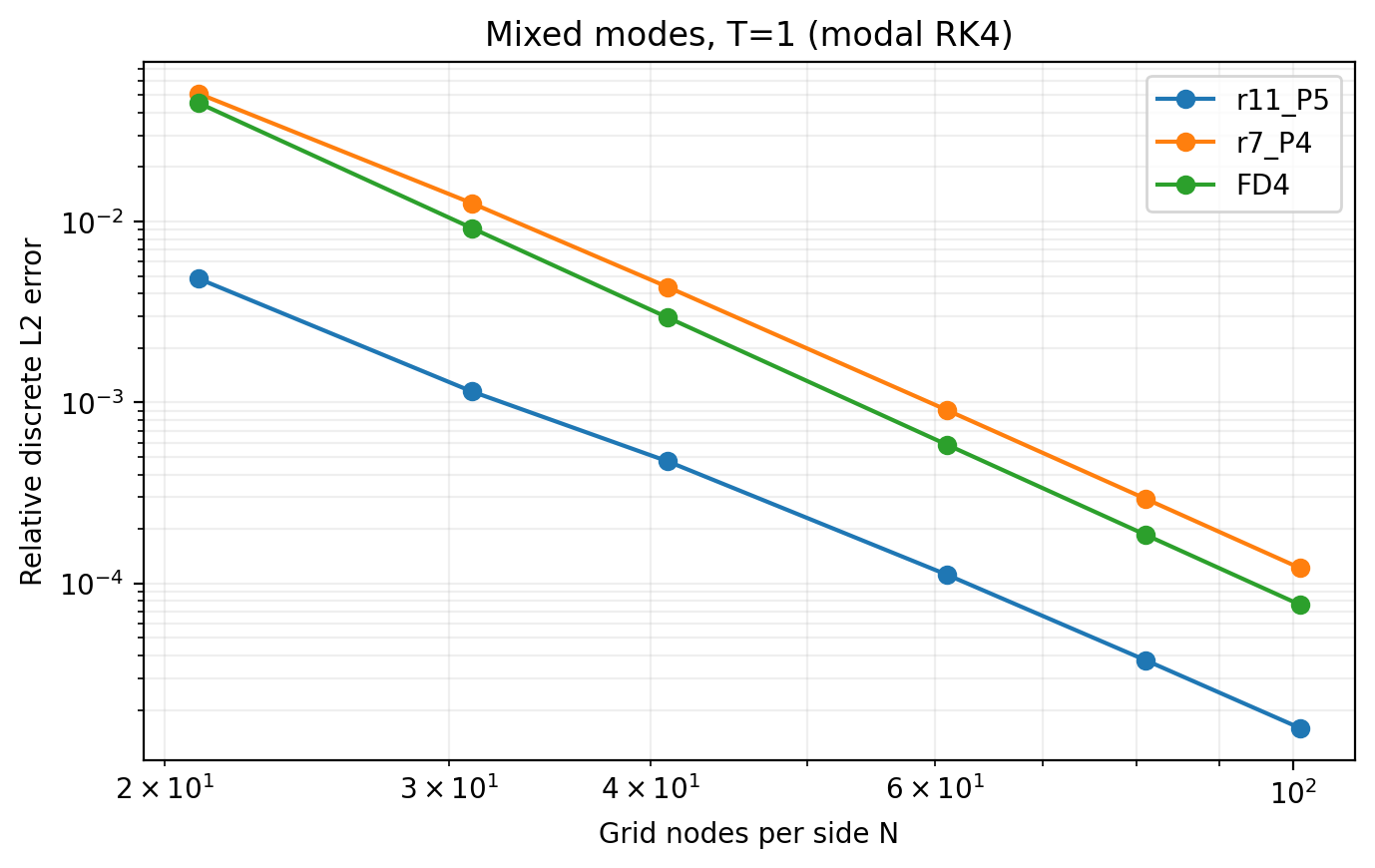}
			\caption{Mixed modes, $T=1$.}
			\label{fig:mixed-T1}
		\end{subfigure}
		\begin{subfigure}[t]{0.48\textwidth}
			\centering
			\includegraphics[width=\linewidth]{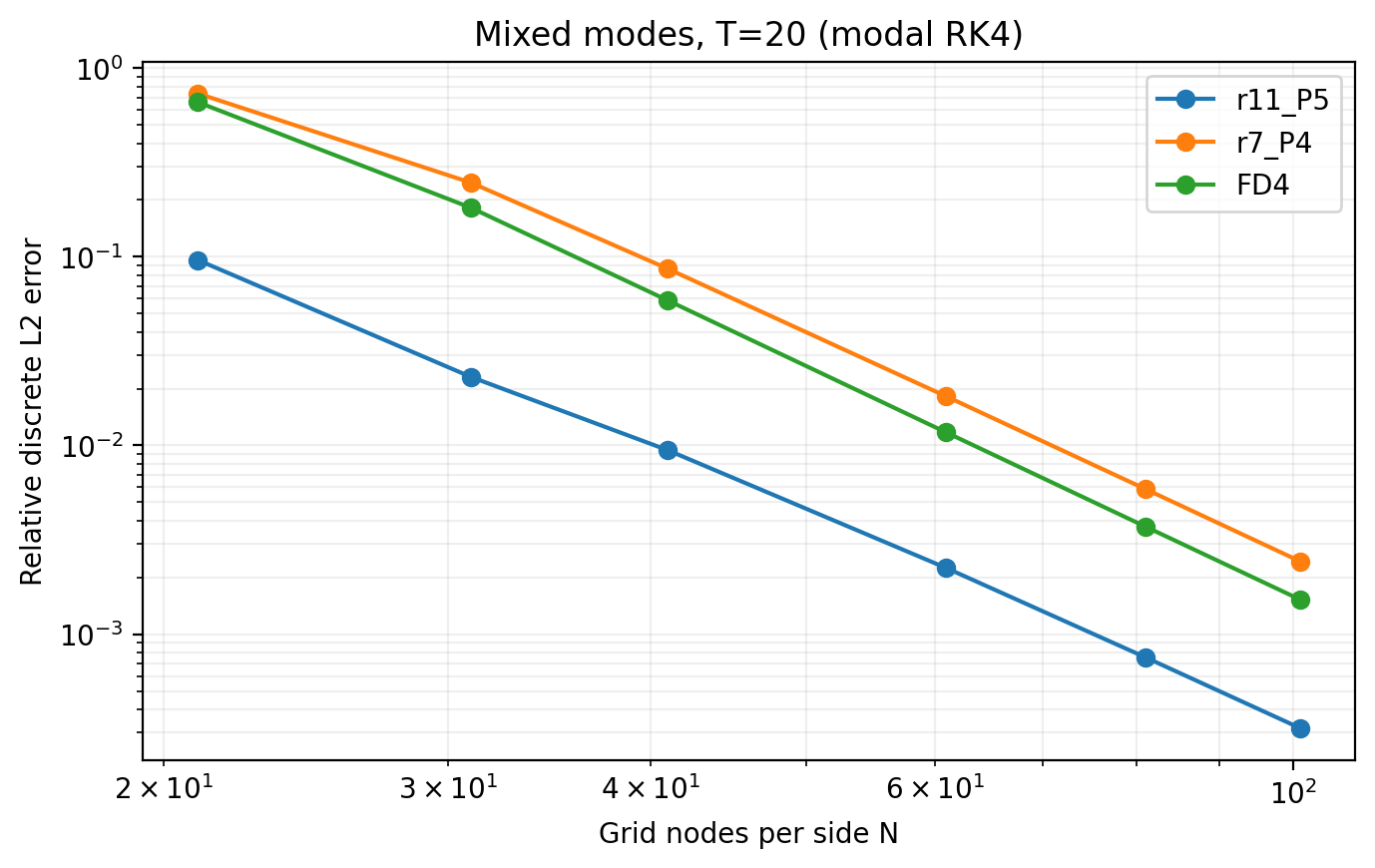}
			\caption{Mixed modes, $T=20$.}
			\label{fig:mixed-T20}
		\end{subfigure}
		\begin{subfigure}[t]{0.48\textwidth}
			\centering
			\includegraphics[width=\linewidth]{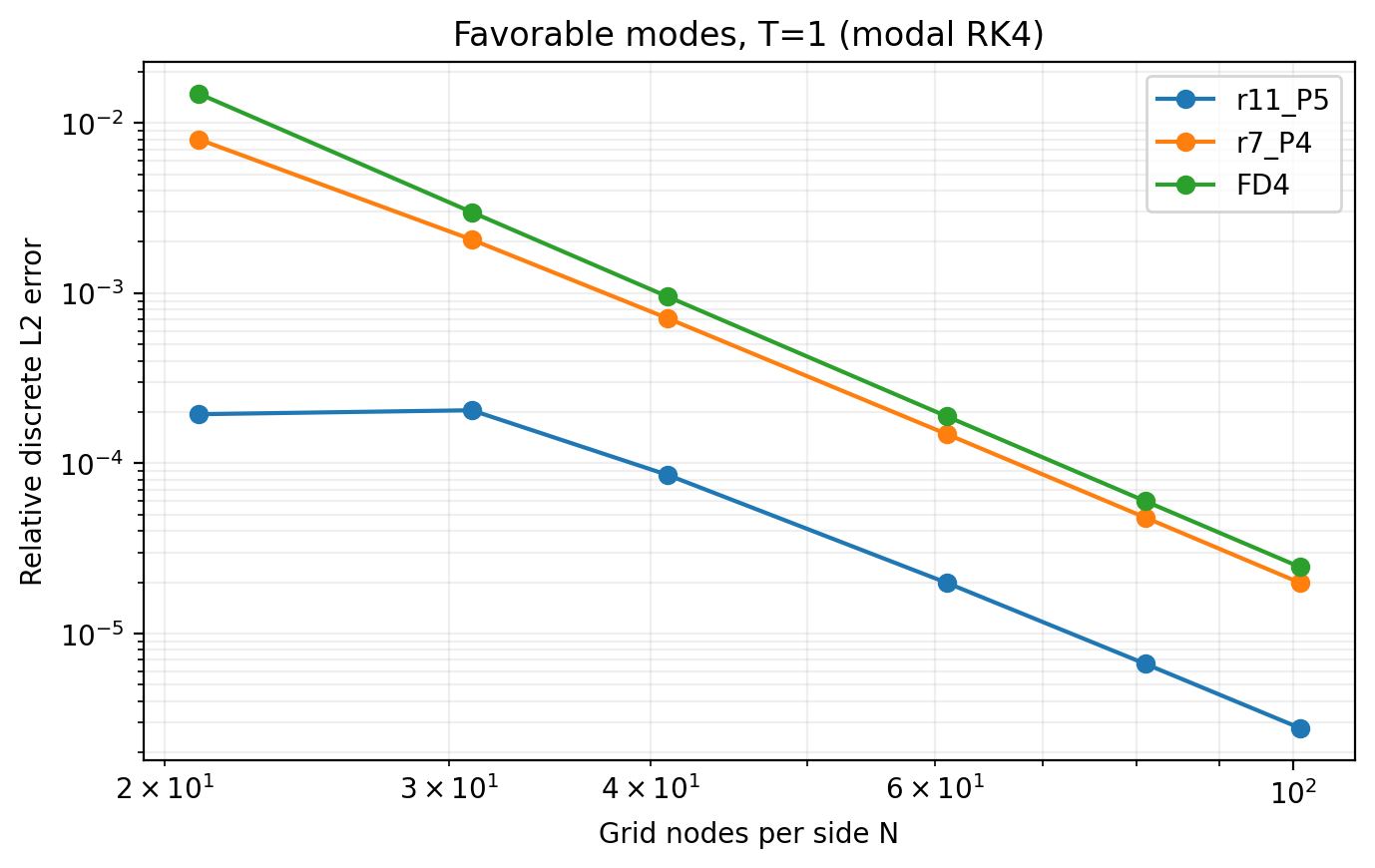}
			\caption{Favorable modes, $T=1$.}
			\label{fig:favorable-T1}
		\end{subfigure}
		\begin{subfigure}[t]{0.48\textwidth}
			\centering
			\includegraphics[width=\linewidth]{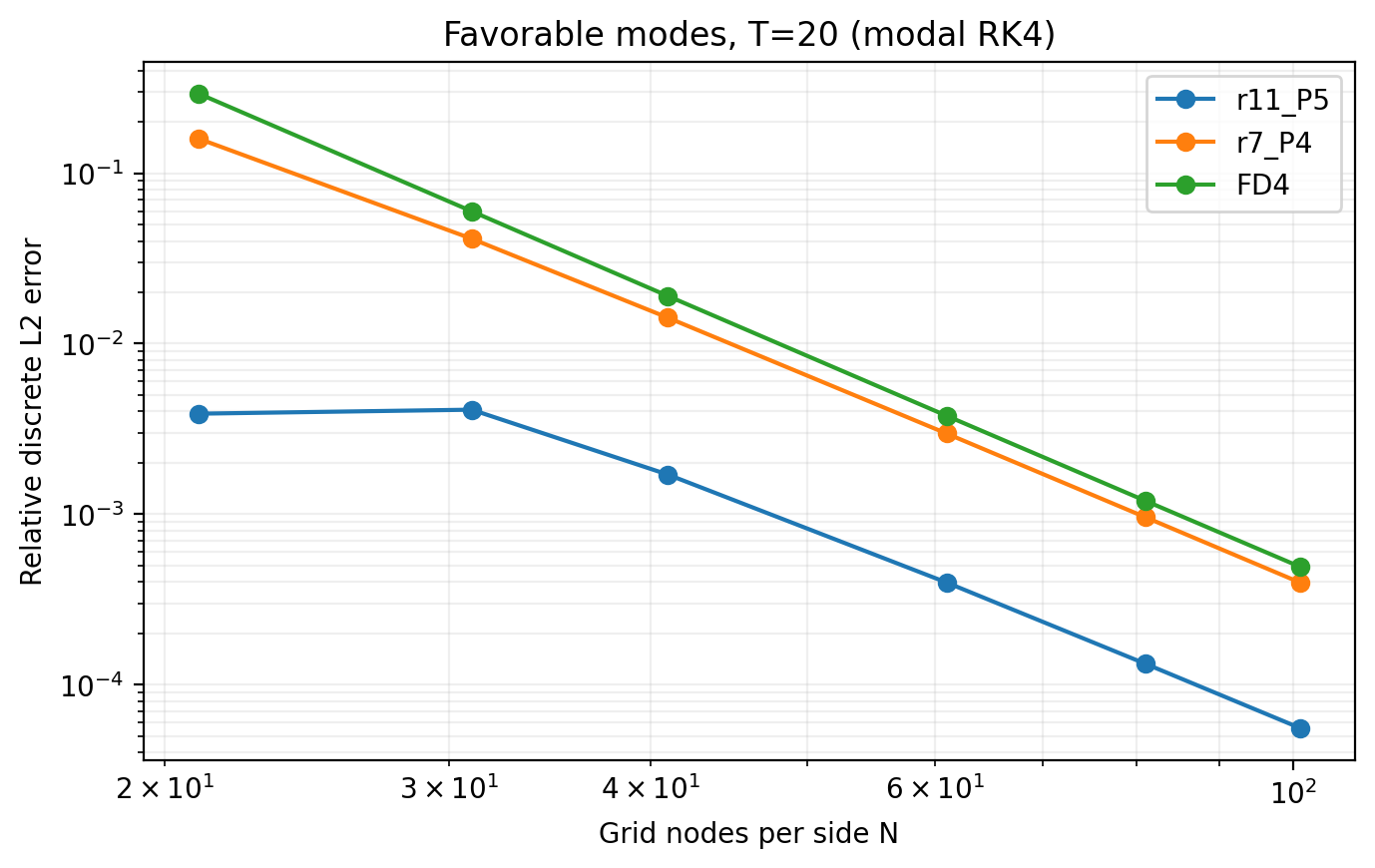}
			\caption{Favorable modes, $T=20$.}
			\label{fig:favorable-T20}
		\end{subfigure}
		
		\caption{Modal convergence for the free
			Schr\"odinger equation. Relative discrete
			$L^2$ errors for the reflected
			$r^{11}+\mathbb P_5$, $r^7+\mathbb P_4$,
			and fourth-order finite-difference (FD4)
			discretizations using modal RK4 on six
			spatial grids. Results are shown for
			(a) mixed modes at $T=1$,
			(b) mixed modes at $T=20$,
			(c) favorable modes at $T=1$, and
			(d) favorable modes at $T=20$.}
		\label{fig:modal-convergence}
	\end{figure}

	\cref{fig:modal-convergence}(a)--(d)
	presents the relative discrete $L^2$ errors
	for the reflected $r^{11}+\mathbb P_5$,
	$r^7+\mathbb P_4$, and FD4 discretizations
	under identical modal RK4 settings.
	The time step is chosen sufficiently small
	to make the RK4 temporal error negligible
	relative to the spatial discretization error,
	so that the reported convergence behavior
	primarily reflects spatial dispersion.
	
	All three methods exhibit decreasing errors
	with grid refinement overall, although a
	slight nonmonotonic variation is observed
	for the reflected $r^{11}+\mathbb P_5$
	formulation in the favorable-mode test
	between the two coarsest grids.
	Among the tested methods, the reflected
	$r^{11}+\mathbb P_5$ formulation consistently
	produces the smallest errors for both
	mixed and favorable modes across all
	six spatial grids.
	Comparing \cref{fig:modal-convergence}(a)
	and \cref{fig:modal-convergence}(b),
	as well as \cref{fig:modal-convergence}(c)
	and \cref{fig:modal-convergence}(d),
	the errors generally increase when the
	final time is extended from $T=1$ to $T=20$.
	This behavior is consistent with the
	accumulation of numerical phase error
	over longer integration intervals.
	The relative accuracy of $r^7+\mathbb P_4$
	and FD4 depends on the selected modes.
	FD4 performs better for the mixed-mode tests,
	whereas $r^7+\mathbb P_4$ is more accurate
	for the favorable-mode tests at both
	$T=1$ and $T=20$.
	This observation is consistent with the
	directional dispersion analysis in Section~3.

	These results demonstrate the favorable
	accuracy of the $r^{11}+\mathbb P_5$
	formulation for the selected benchmarks,
	without establishing its universal superiority
	for arbitrary initial data or modes.

	Having examined the modal accuracy of the
	reflected discretizations, we next investigate
	the sensitivity of conventional KNN RBF-FD
	operators to stencil selection. This experiment
	focuses on frequency accuracy, matrix symmetry,
	and boundary consistency rather than direct
	convergence comparisons with the reflected methods.
	
\begin{figure}[pos=htbp]
		\centering
		
		\begin{subfigure}[t]{0.32\textwidth}
			\centering
			\includegraphics[width=\linewidth]{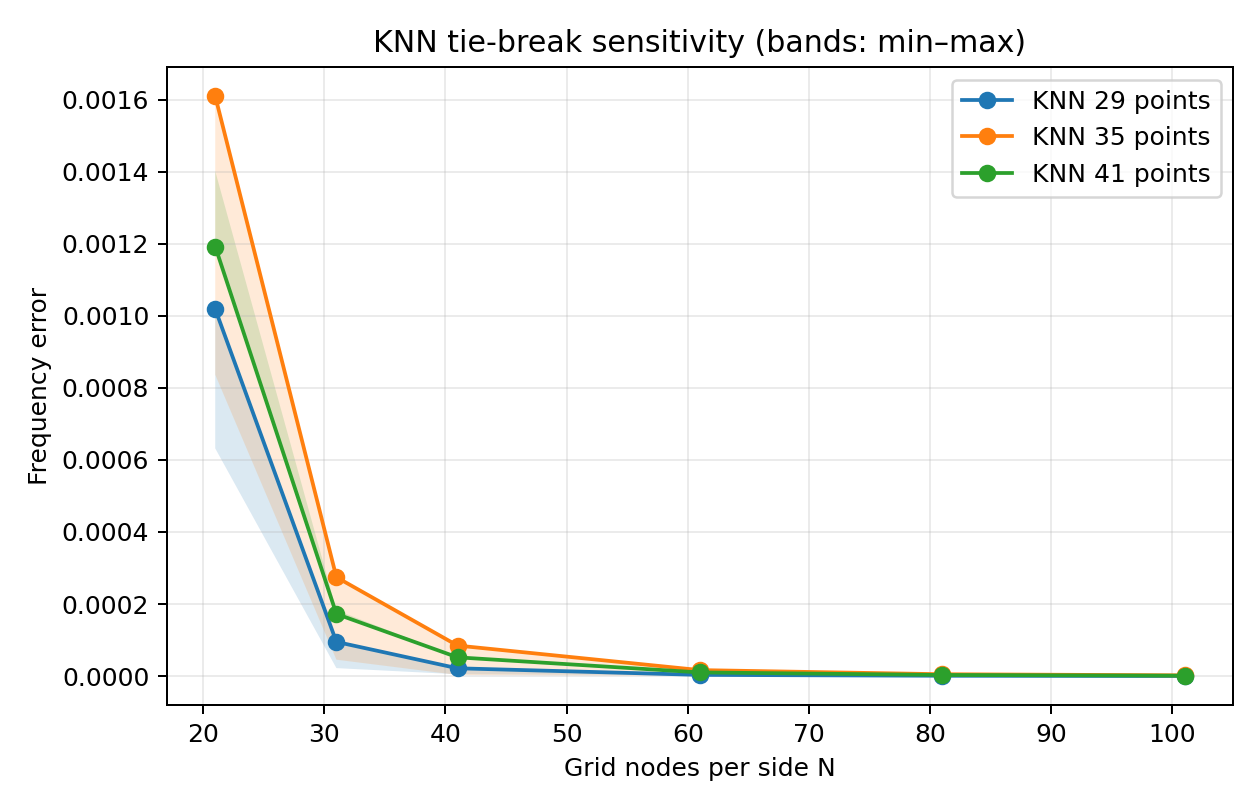}
			\caption{KNN frequency sensitivity.}
			\label{fig:knn-ties-frequency}
		\end{subfigure}
		\hfill
		\begin{subfigure}[t]{0.32\textwidth}
			\centering
			\includegraphics[width=\linewidth]{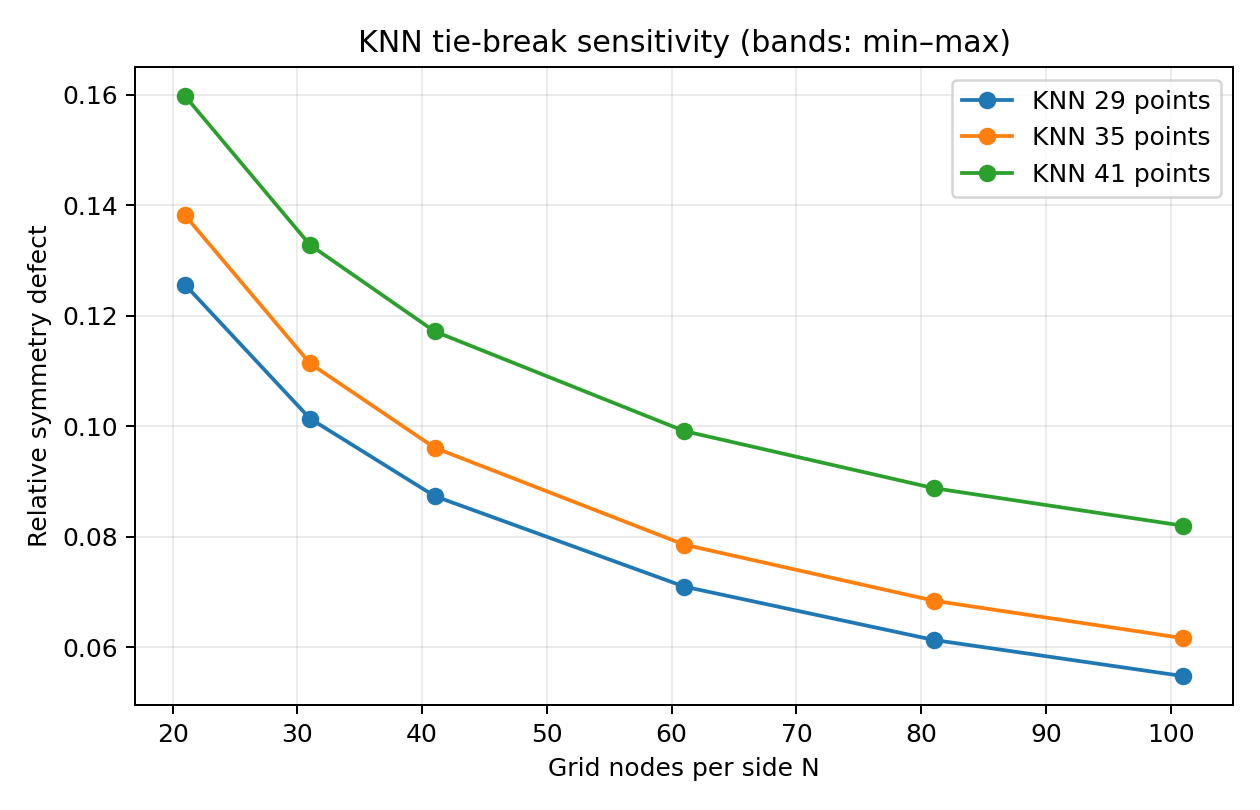}
			\caption{KNN symmetry sensitivity.}
			\label{fig:knn-ties-symmetry}
		\end{subfigure}
		\hfill
		\begin{subfigure}[t]{0.32\textwidth}
			\centering
			\includegraphics[width=\linewidth]{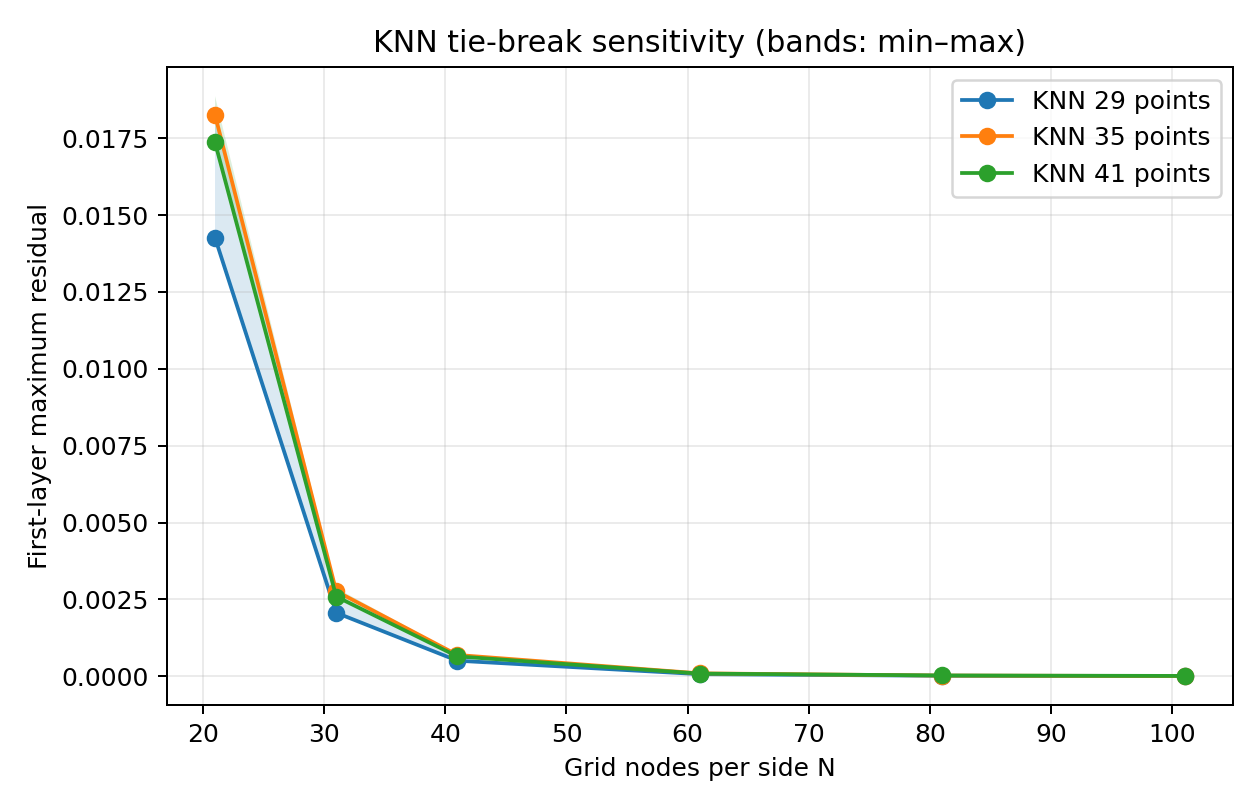}
			\caption{KNN boundary residuals.}
			\label{fig:knn-ties-boundary}
		\end{subfigure}
		
		\caption{Sensitivity of the rank-safe KNN
			RBF-FD discretization to stencil selection.
			Results compare 29-, 35-, and 41-point
			stencils using the $r^{11}+\mathbb P_5$
			formulation with deterministic distance-tie
			breaking. The panels show
			(a) Rayleigh-frequency errors,
			(b) relative global matrix symmetry defects,
			and (c) boundary-layer operator residuals
			under grid refinement.}
		\label{fig:knn-sensitivity}
	\end{figure}

	\paragraph{KNN stencil-selection sensitivity.}
	\cref{fig:knn-sensitivity} and \cref{tab:knn-ties} examine the
	sensitivity of the rank-safe KNN RBF-FD
	discretization to stencil size and
	deterministic distance-tie-breaking policies.
	We consider 29-, 35-, and 41-point stencils
	with three tie-breaking policies:
	\texttt{xy\_ascending},
	\texttt{yx\_ascending}, and
	\texttt{xy\_descending}.
	The kernel--polynomial pair
	$r^{11}+\mathbb P_5$ is fixed throughout,
	so the observed differences arise from
	stencil selection rather than kernel choice.
	
	The Rayleigh-frequency errors, global matrix
	symmetry defects, and boundary-layer operator
	residuals vary with stencil configuration
	and grid resolution. These diagnostics
	measure different properties of the spatial
	operator and need not improve simultaneously
	as stencil size increases. In particular,
	for the 35-point stencil, changing only
	the tie-breaking policy can materially
	affect the Rayleigh-frequency error.
	This demonstrates sensitivity to stencil
	selection without establishing the universal
	superiority of any particular policy.
	Unlike the reflected complete-shell
	construction, which yields an exactly
	symmetric discrete Laplacian on the
	uniform tensor-product grid, the KNN
	construction does not generally enforce
	global matrix symmetry. Increasing the
	number of KNN neighbors therefore does
	not guarantee the structure-preserving
	properties of the reflected method.
	This experiment characterizes the
	sensitivity of conventional KNN stencil
	selection rather than providing a direct
	accuracy ranking between KNN and the
	reflected discretization.

	\begin{table}[H]\centering\small
		\caption{Observed range across the three KNN tie policies for selected grids.
			The range is not a statistical confidence interval.}\label{tab:knn-ties}
		\begin{tabular}{rrrrr}\toprule
			$N$ & KNN points & Minimum frequency error & Maximum frequency error & Symmetry-defect range\\\midrule
			31 & 29 & 2.231e-05 & 9.505e-05 & 1.014e-01--1.014e-01 \\
			101 & 29 & 4.012e-07 & 4.301e-07 & 5.474e-02--5.474e-02 \\
			31 & 35 & 4.567e-05 & 2.746e-04 & 1.114e-01--1.114e-01 \\
			101 & 35 & 1.202e-07 & 2.200e-06 & 6.161e-02--6.161e-02 \\
			31 & 41 & 1.725e-04 & 1.832e-04 & 1.328e-01--1.328e-01 \\
			101 & 41 & 1.200e-06 & 1.358e-06 & 8.196e-02--8.196e-02 \\
			\bottomrule\end{tabular}\end{table}

	\subsection{Matched KNN and reflected operator diagnostics}
	\label{subsec:matched-higher-grid}
	
	To isolate the effect of stencil construction,
	we compare the reflected and rank-safe KNN
	RBF-FD discretizations using the same
	$r^{11}+\mathbb P_5$ kernel--polynomial pair,
	Cartesian domain, and homogeneous Dirichlet
	boundary conditions.
	
	The reflected method uses the 29-point
	complete shell $S_{29}$. The KNN method
	selects physical nodes in increasing distance
	order while enforcing polynomial unisolvency.
	A polynomially independent subset is selected
	first, after which the remaining stencil
	positions are filled according to the
	prescribed distance and tie-breaking rules.
	Consequently, near the boundary, the selected
	stencil need not consist of the literal
	nearest physical neighbors.
	All 21 degree-five polynomial columns
	are required to be linearly independent,
	and the local augmented systems are solved
	using high-precision arithmetic.
	We consider six spatial grids,
	$N=21,31,41,61,81,101$. At each grid,
	we evaluate the reflected method and nine
	KNN configurations obtained from three
	stencil sizes (29, 35, and 41 points)
	and three deterministic tie-breaking policies.

	For the sampled sine mode
	\(
	\phi_{23}(x,y)=\sin(2x)\sin(3y),
	\)
	the continuous Laplacian satisfies
	\(
	-\Delta\phi_{23}=13\phi_{23}.
	\)
	Consequently, the exact solution of the
	free Schr\"odinger equation evolves as
	$e^{-13it}\phi_{23}$, with angular
	frequency $13$.
	For a discrete Laplacian $L$, we define
	the Rayleigh-frequency estimate
	\begin{equation}
		\omega_h^*
		=-\frac{\phi_{23}^T L\phi_{23}}
		{\|\phi_{23}\|_2^2},
	\end{equation}
	and its error
	$|\omega_h^*-13|$.
	This quantity measures how accurately
	the discrete operator reproduces the
	continuous modal frequency.
	We also define the mode-shape residual
	\begin{equation}
		\eta=L\phi_{23}+\omega_h^*\phi_{23}.
	\end{equation}
	While the Rayleigh-frequency error
	measures the discrepancy in the
	estimated frequency, the normalized
	residual $\|\eta\|_2/\|\phi_{23}\|_2$
	measures how closely the sampled sine
	mode behaves as an eigenvector of $L$.
	Together with the global matrix symmetry
	defect, these quantities characterize
	the spatial operator rather than the
	fully discrete time-dependent solution.
	
	\begin{table}[H]\centering\small
		\caption{Six-grid matched $r^{11}+\mathbb P_5$ comparison. KNN uses 29 points and the
			\texttt{xy\_ascending} tie policy; values are from the C++ spatial-operator study.}
		\label{tab:matched-r11}
		\resizebox{\linewidth}{!}{\begin{tabular}{rrrrrr}\toprule
				$N$ & Reflected $|\omega_h^*-13|$ & KNN $|\omega_h^*-13|$ & KNN $\eta_{\rm sym}$ & Reflected shape residual & KNN shape residual\\\midrule
				21 & 1.023e-04 & 1.018e-03 & 1.256e-01 & 3.160e-14 & 7.044e-03 \\
				31 & 3.776e-05 & 9.505e-05 & 1.014e-01 & 7.092e-14 & 8.310e-04 \\
				41 & 1.401e-05 & 2.180e-05 & 8.733e-02 & 1.633e-13 & 1.763e-04 \\
				61 & 3.064e-06 & 3.526e-06 & 7.094e-02 & 4.864e-13 & 1.941e-05 \\
				81 & 1.003e-06 & 1.065e-06 & 6.129e-02 & 8.356e-13 & 4.033e-06 \\
				101 & 4.171e-07 & 4.301e-07 & 5.474e-02 & 1.355e-12 & 1.191e-06 \\
				\bottomrule\end{tabular}}\end{table}
	
	\cref{tab:matched-r11} reports the
	six-grid comparison between the
	29-point reflected operator and
	the rank-safe KNN operator with
	the \texttt{xy\_ascending} tie policy.
	The reflected construction produces
	smaller Rayleigh-frequency errors
	than this KNN configuration on all
	six grids, while its mode-shape
	residual remains near floating-point
	precision. This is consistent with
	the exact discrete sine eigenvector
	property of the reflected operator.
	The KNN frequency and mode-shape
	errors decrease under grid refinement,
	but the KNN matrix remains
	nonsymmetric. At $N=101$, the
	frequency errors are
	$4.171\times10^{-7}$ and
	$4.301\times10^{-7}$ for the
	reflected and KNN methods,
	respectively, whereas the KNN
	symmetry defect as defined in \eqref{eq:symdefect} remains
	$5.474\times10^{-2}$.
	Thus, accurate modal frequencies
	do not necessarily imply exact
	sine-mode invariance or matrix
	symmetry.

	\begin{figure}[pos=htbp]
		\centering
		
		\begin{subfigure}[t]{0.45\textwidth}
			\centering
			\includegraphics[width=\linewidth]{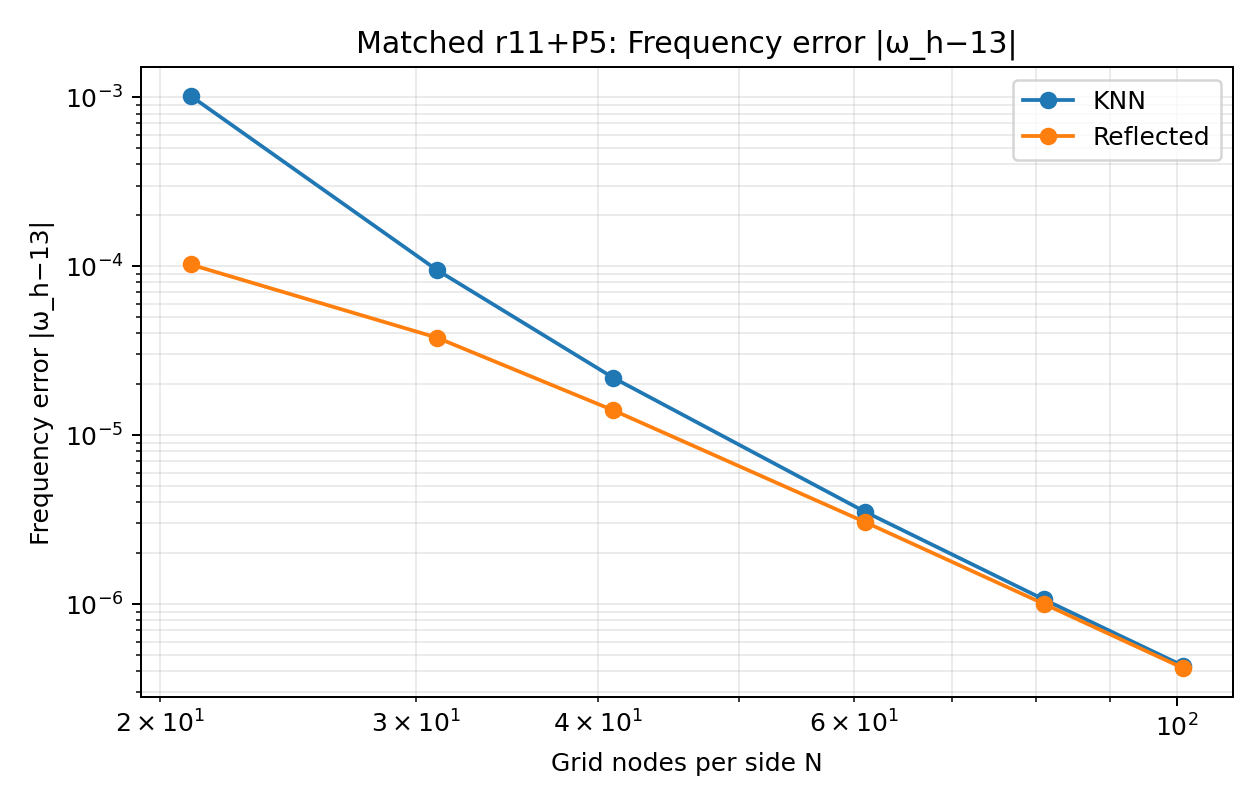}
			\caption{Rayleigh-frequency error.}
			\label{fig:higher-frequency}
		\end{subfigure}
		\begin{subfigure}[t]{0.45\textwidth}
			\centering
			\includegraphics[width=\linewidth]{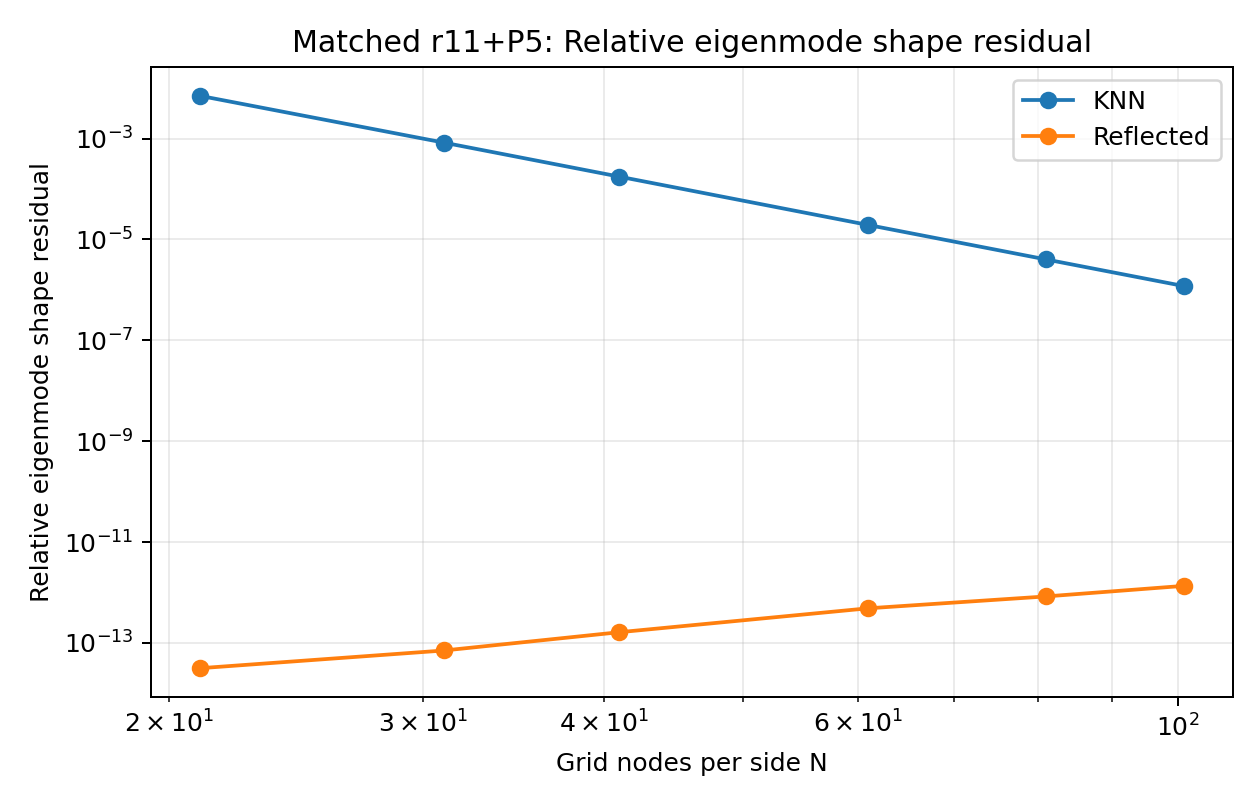}
			\caption{Mode-shape residual.}
			\label{fig:higher-shape}
		\end{subfigure}
		\begin{subfigure}[t]{0.45\textwidth}
			\centering
			\includegraphics[width=\linewidth]{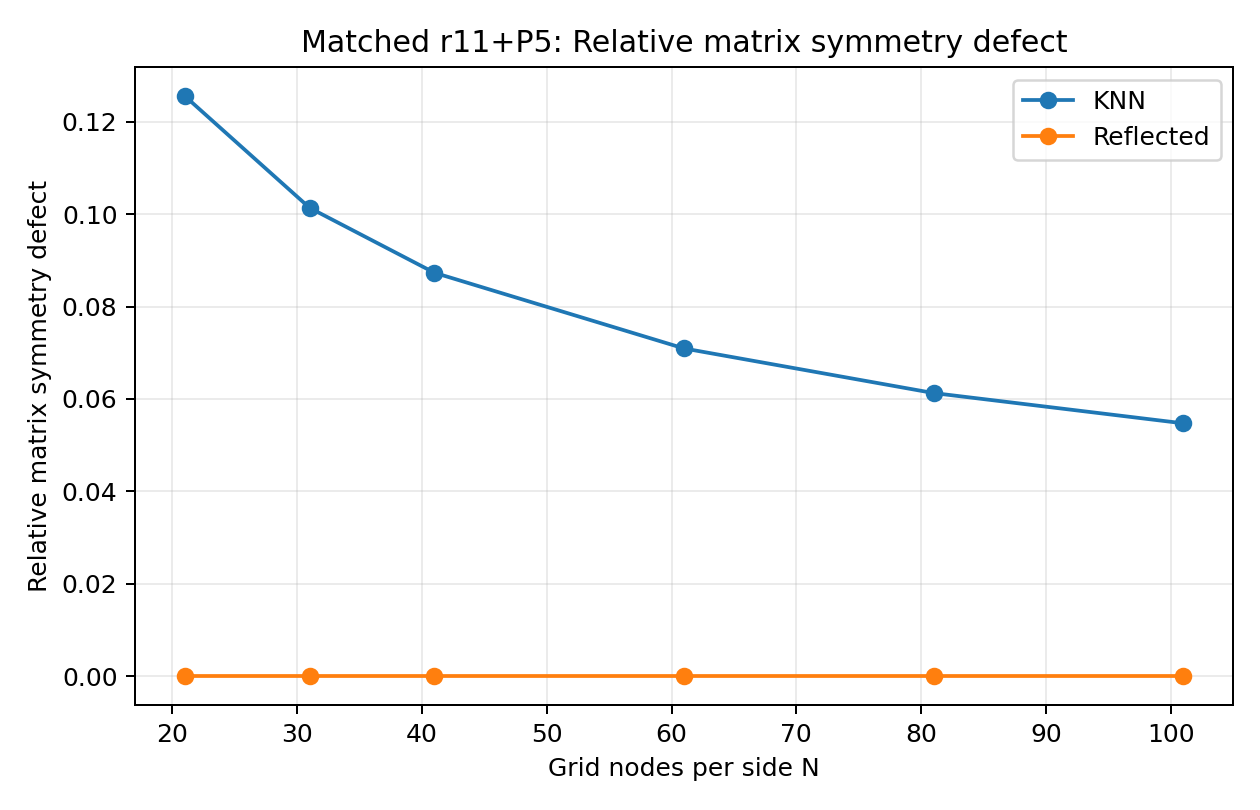}
			\caption{Global matrix symmetry defect.}
			\label{fig:higher-symmetry}
		\end{subfigure}
		\begin{subfigure}[t]{0.45\textwidth}
			\centering
			\includegraphics[width=\linewidth]{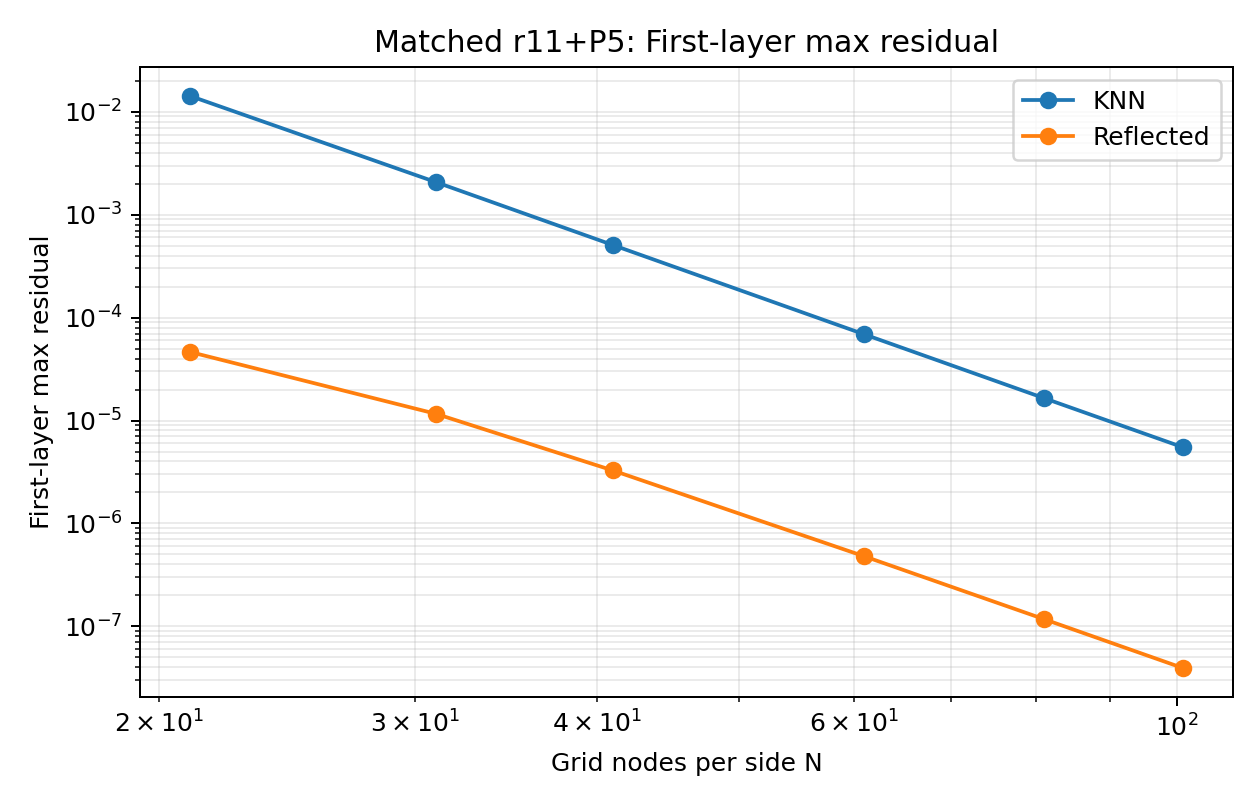}
			\caption{First- and second-layer residuals.}
			\label{fig:higher-boundary}
		\end{subfigure}
		
		\caption{Matched $r^{11}+\mathbb P_5$ KNN and reflected
			RBF-FD spatial-operator comparisons over six grid sizes
			($N=21,31,41,61,81,101$):
			(a) Rayleigh-frequency error;
			(b) mode-shape residual;
			(c) global matrix symmetry defect;
			and (d) first- and second-interior-layer
			sine-mode operator residuals.}
		\label{fig:higher-grid-comparison}
	\end{figure}
	\cref{fig:higher-grid-comparison}
	provides the corresponding six-grid
	operator diagnostics, including
	frequency error, mode-shape residual,
	matrix symmetry, and boundary-layer
	residuals. The boundary-layer quantities
	measure the residual of the sampled sine
	mode at the first two interior layers;
	they are not errors in the strongly
	imposed homogeneous Dirichlet values.
	These results illustrate an important
	distinction: a nonsymmetric KNN operator
	can achieve accurate modal frequencies
	without possessing the exact symmetry
	and sine-mode invariance of the
	reflected complete-shell construction.

	\subsection{Sparse time integration and conservation}
	\label{subsec:new-sparse}
	
	We examine the time-integration accuracy,
	conservation behavior, and spectral stability
	of the matched $r^{11}+\mathbb P_5$
	KNN and reflected RBF-FD operators.
	
	The matched-operator experiments include
	sparse RK4 and implicit-midpoint trajectories
	at $N=21,31$ through $T=2$, together with
	long-time stability diagnostics at $N=21$
	through $T=20$. For the $T=2$ experiments,
	the RK4 and midpoint time steps are
	$\Delta t=0.0005$ and $\Delta t=0.005$,
	respectively.
	For the symmetric reflected operator,
	implicit midpoint preserves the discrete
	quadratic mass and energy invariants
	in exact arithmetic. Observed numerical
	drift may also reflect errors in solving
	the implicit linear systems.
	For nonsymmetric KNN operators, the
	quadratic form
	$-(u^TLu+v^TLv)/2$ is a diagnostic
	quantity rather than a generally
	conserved discrete Hamiltonian.
	
	\begin{figure}[pos=htbp]
		\centering
		
		\begin{subfigure}[t]{0.32\textwidth}
			\centering
			\includegraphics[width=\linewidth]
			{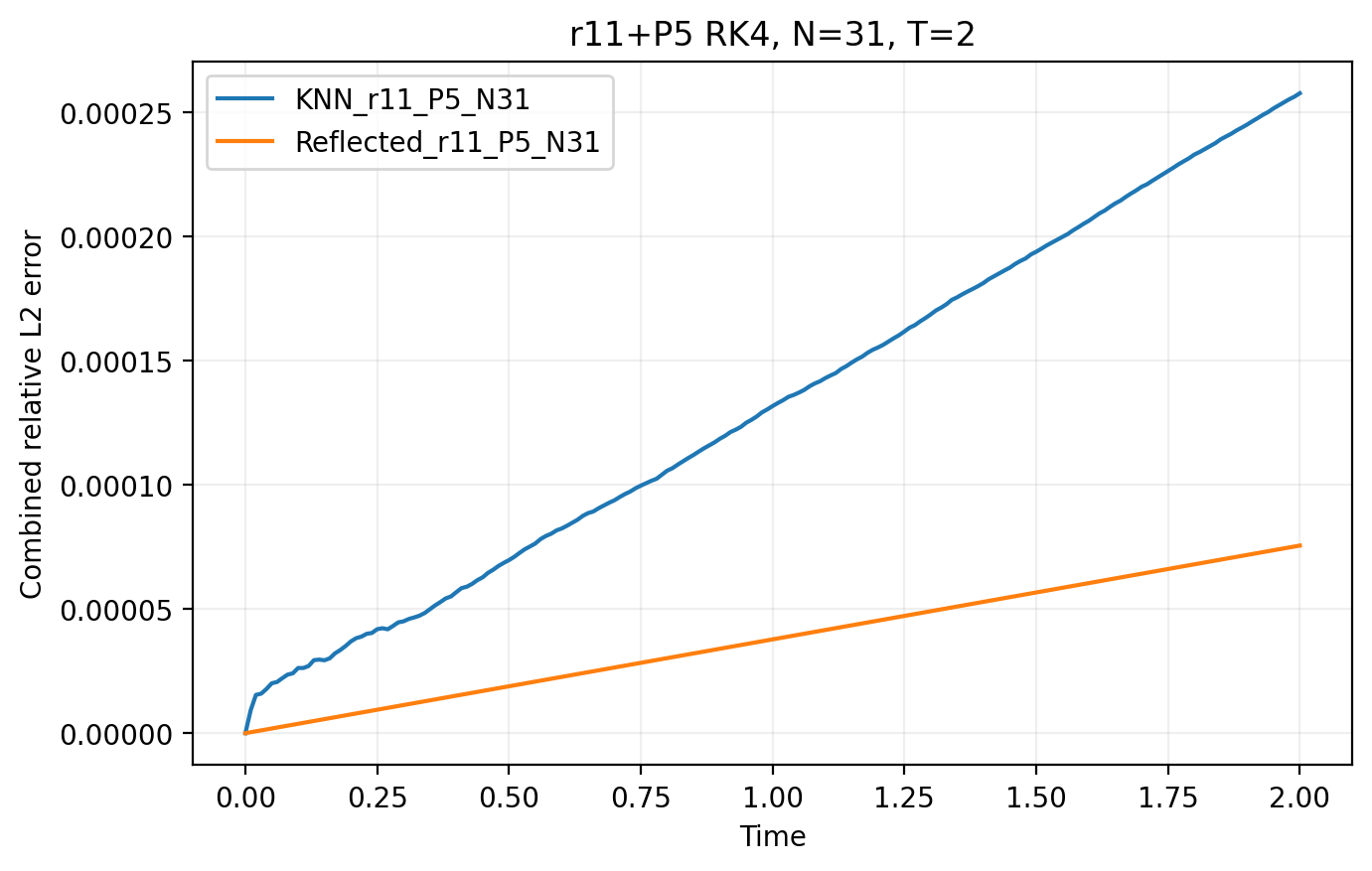}
			\caption{RK4 relative $L^2$ error, $N=31$.}
			\label{fig:sparse-rk4}
		\end{subfigure}
		\begin{subfigure}[t]{0.32\textwidth}
			\centering
			\includegraphics[width=\linewidth]
			{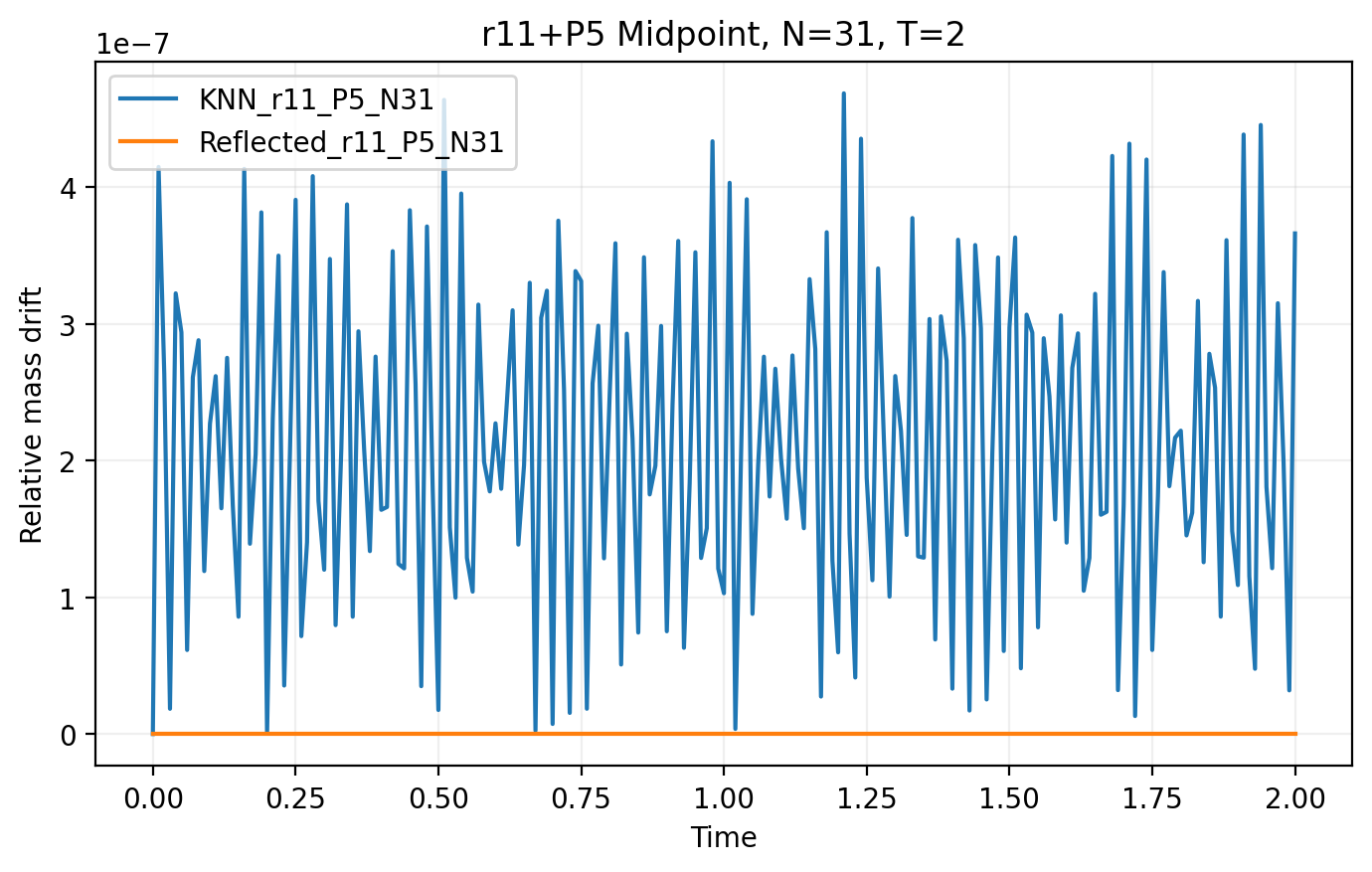}
			\caption{Midpoint mass drift, $N=31$.}
			\label{fig:midpoint-mass}
		\end{subfigure}
		\begin{subfigure}[t]{0.32\textwidth}
			\centering
			\includegraphics[width=\linewidth]
			{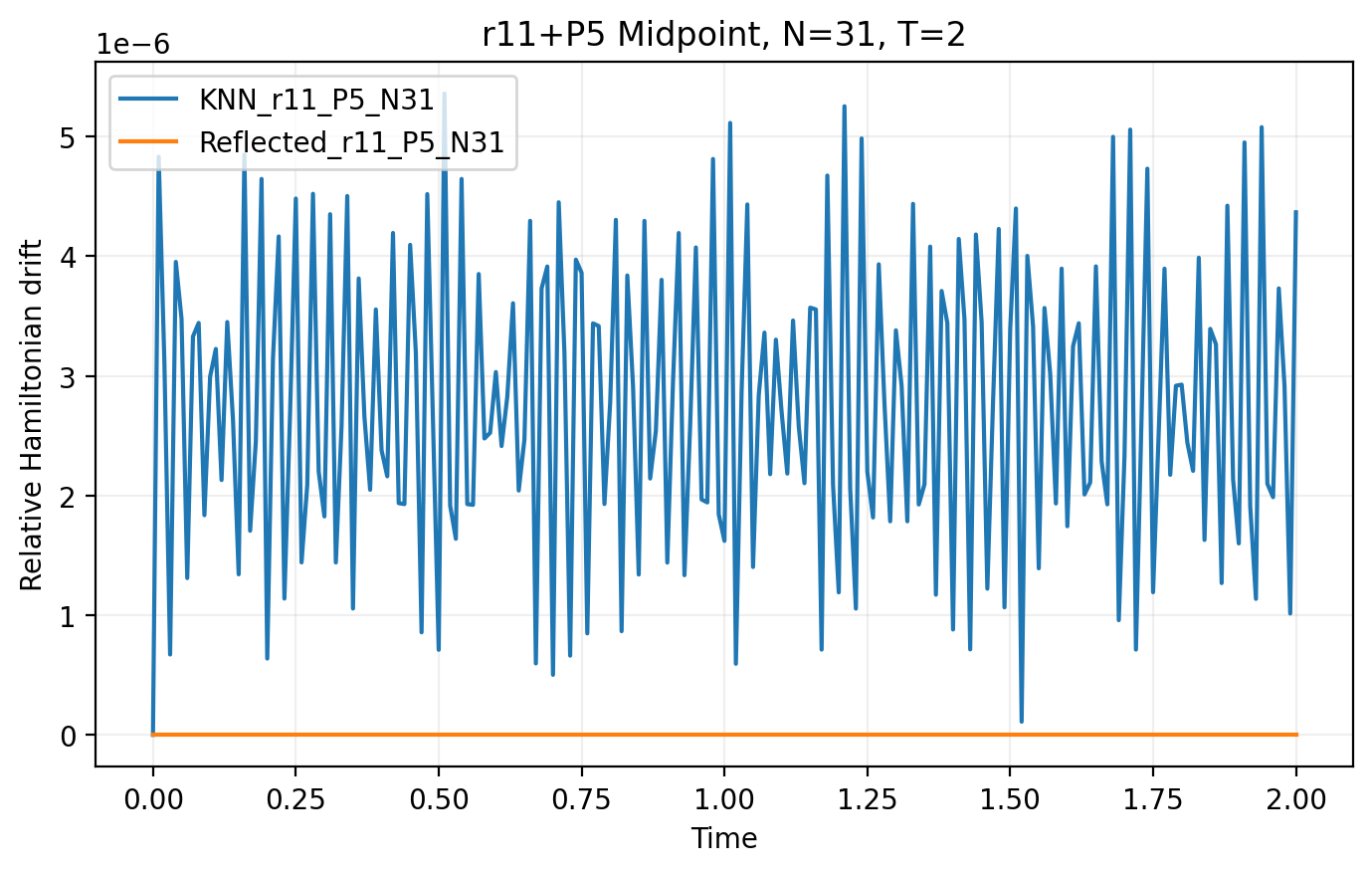}
			\caption{Midpoint energy drift, $N=31$.}
			\label{fig:midpoint-energy}
		\end{subfigure}
		
		\caption{Sparse time-integration and
			conservation diagnostics for the matched
			$r^{11}+\mathbb P_5$ KNN and reflected
			RBF-FD operators at $N=31$ through $T=2$.
			(a) Relative $L^2$ solution-error history
			under RK4; (b) relative mass drift under
			implicit midpoint; and (c) relative
			energy-diagnostic drift under implicit
			midpoint. For nonsymmetric KNN, the
			quadratic energy diagnostic is not
			generally conserved.}
		\label{fig:time-integration-comparison}
	\end{figure}
	
	\begin{figure}[pos=htbp]
		\centering
		\begin{subfigure}[t]{0.4\textwidth}
			\centering
			\includegraphics[width=\linewidth]
			{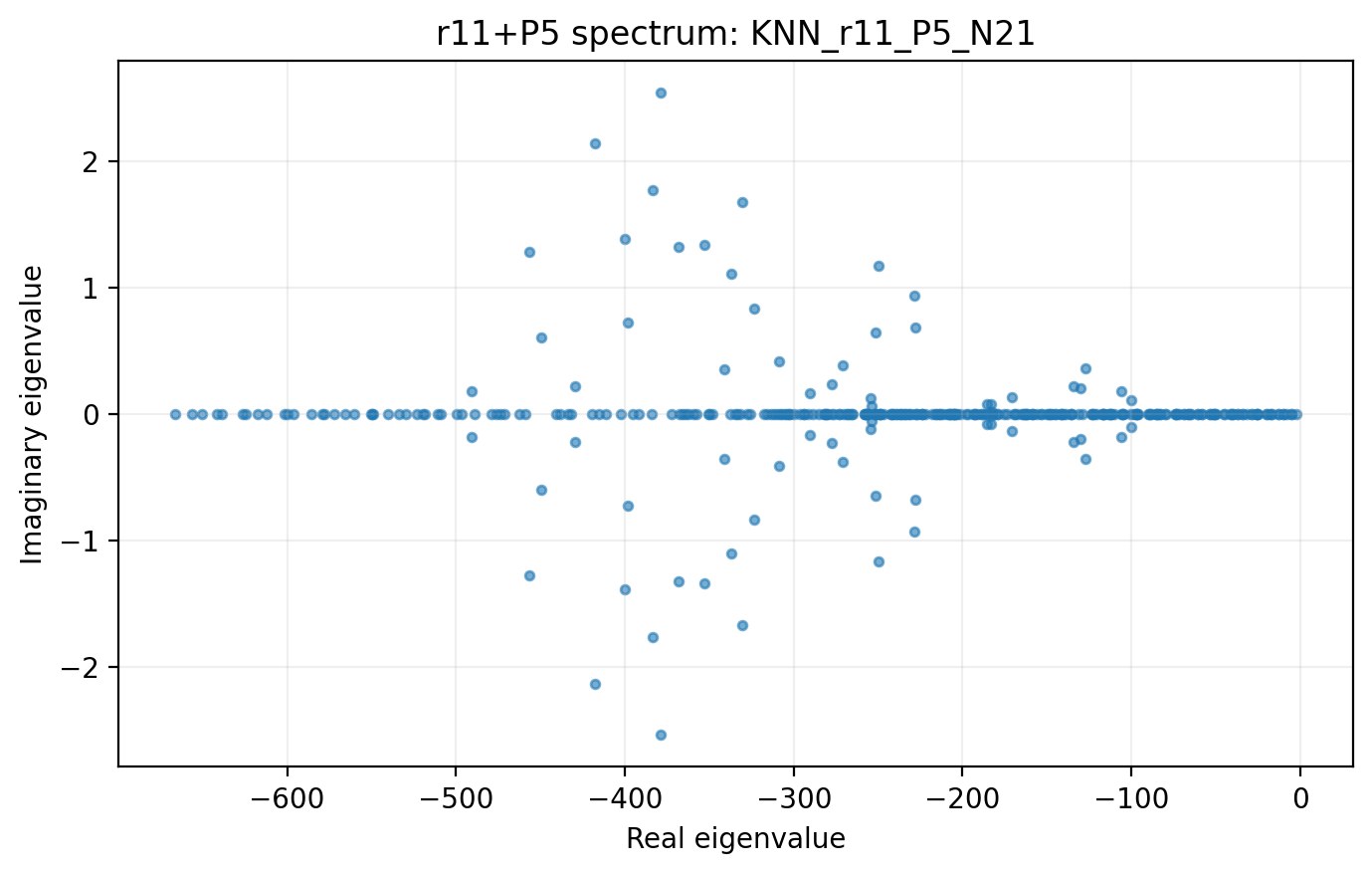}
			\caption{KNN spatial spectrum.}
			\label{fig:knn-r11-spectrum}
		\end{subfigure}
		\begin{subfigure}[t]{0.4\textwidth}
			\centering
			\includegraphics[width=\linewidth]
			{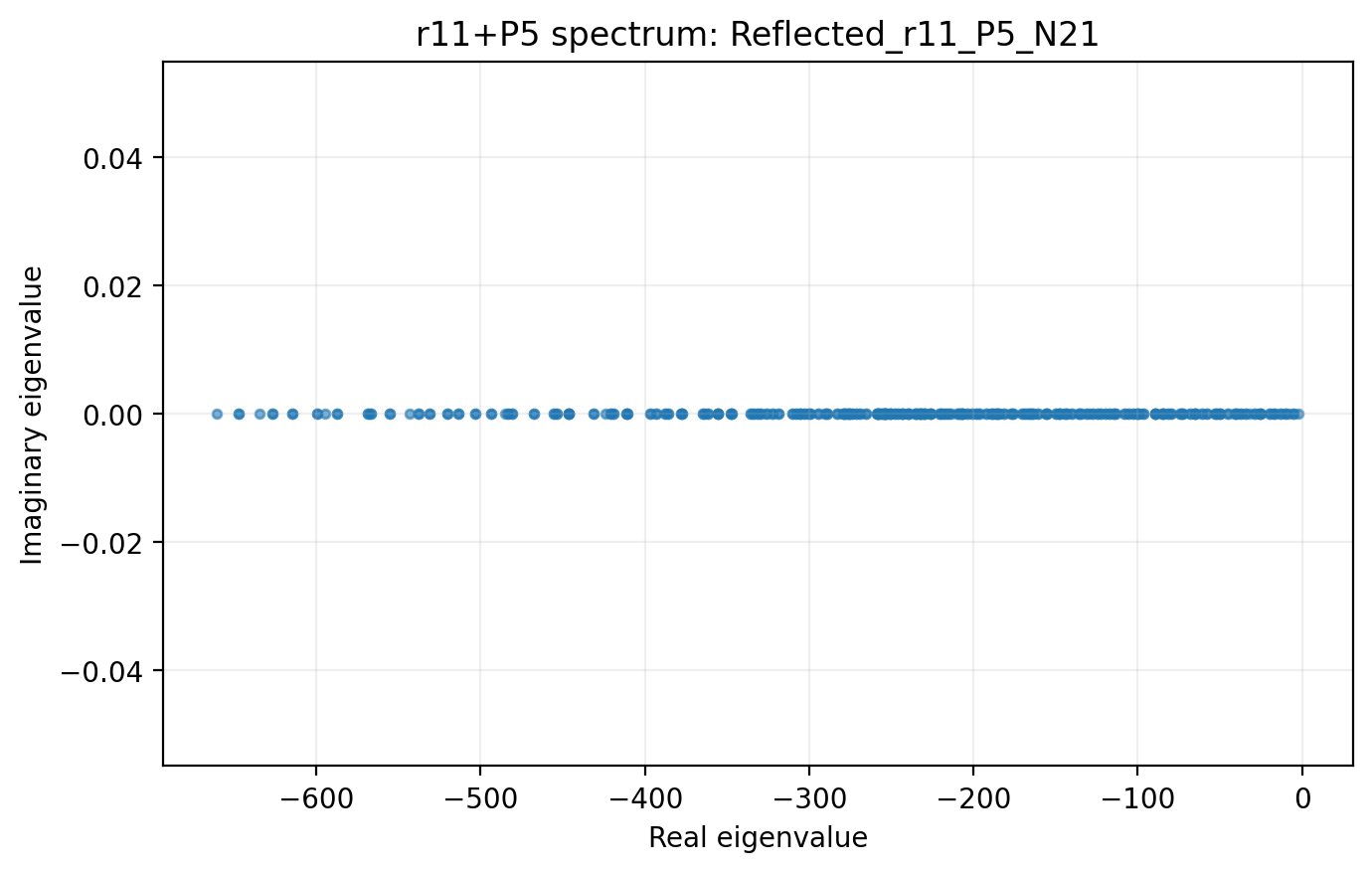}
			\caption{Reflected spatial spectrum.}
			\label{fig:ref-r11-spectrum}
		\end{subfigure}
		\begin{subfigure}[t]{0.4\textwidth}
			\centering
			\includegraphics[width=\linewidth]
			{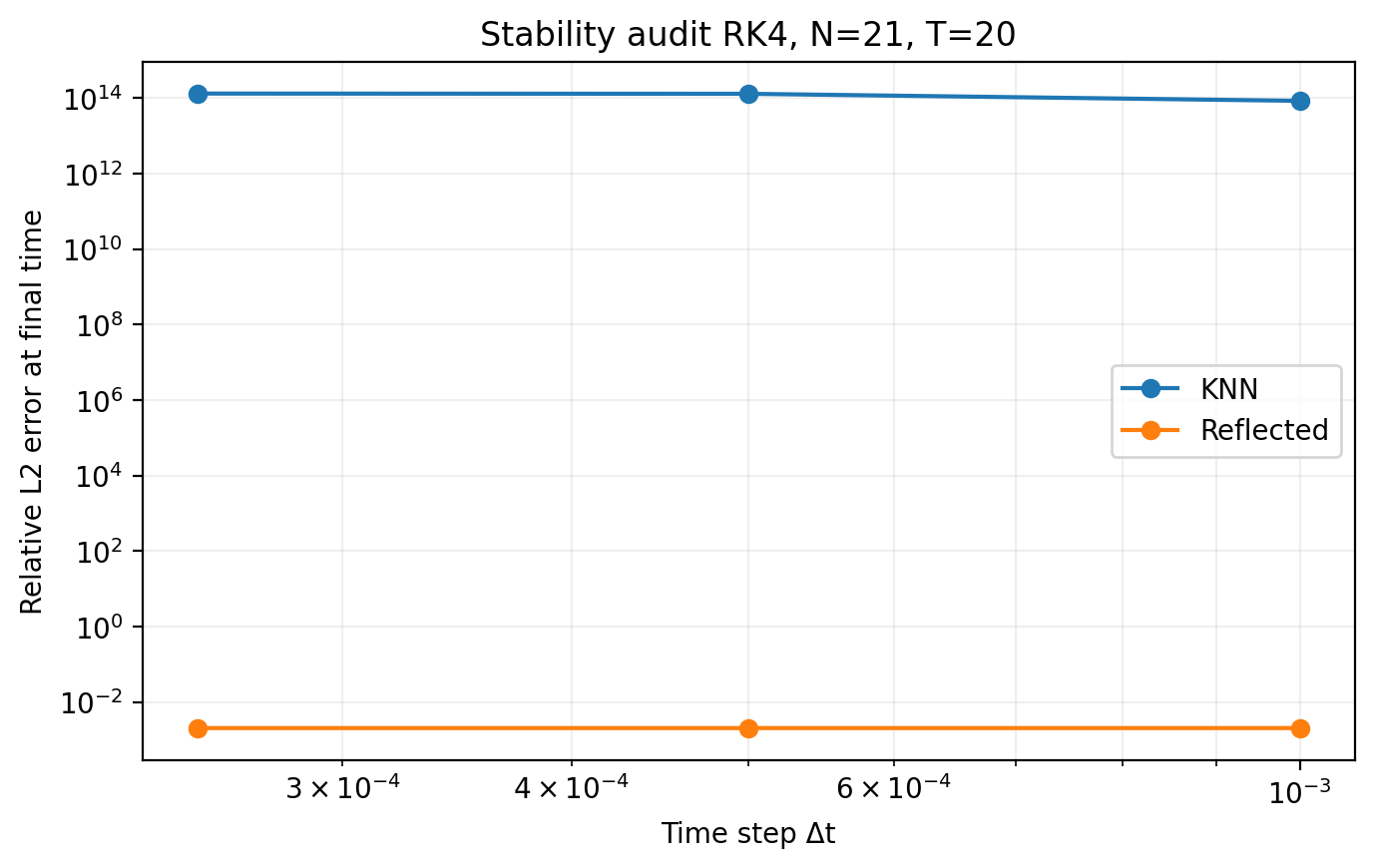}
			\caption{RK4 step-size audit.}
			\label{fig:r11-audit-rk4}
		\end{subfigure}
		\begin{subfigure}[t]{0.4\textwidth}
			\centering
			\includegraphics[width=\linewidth]
			{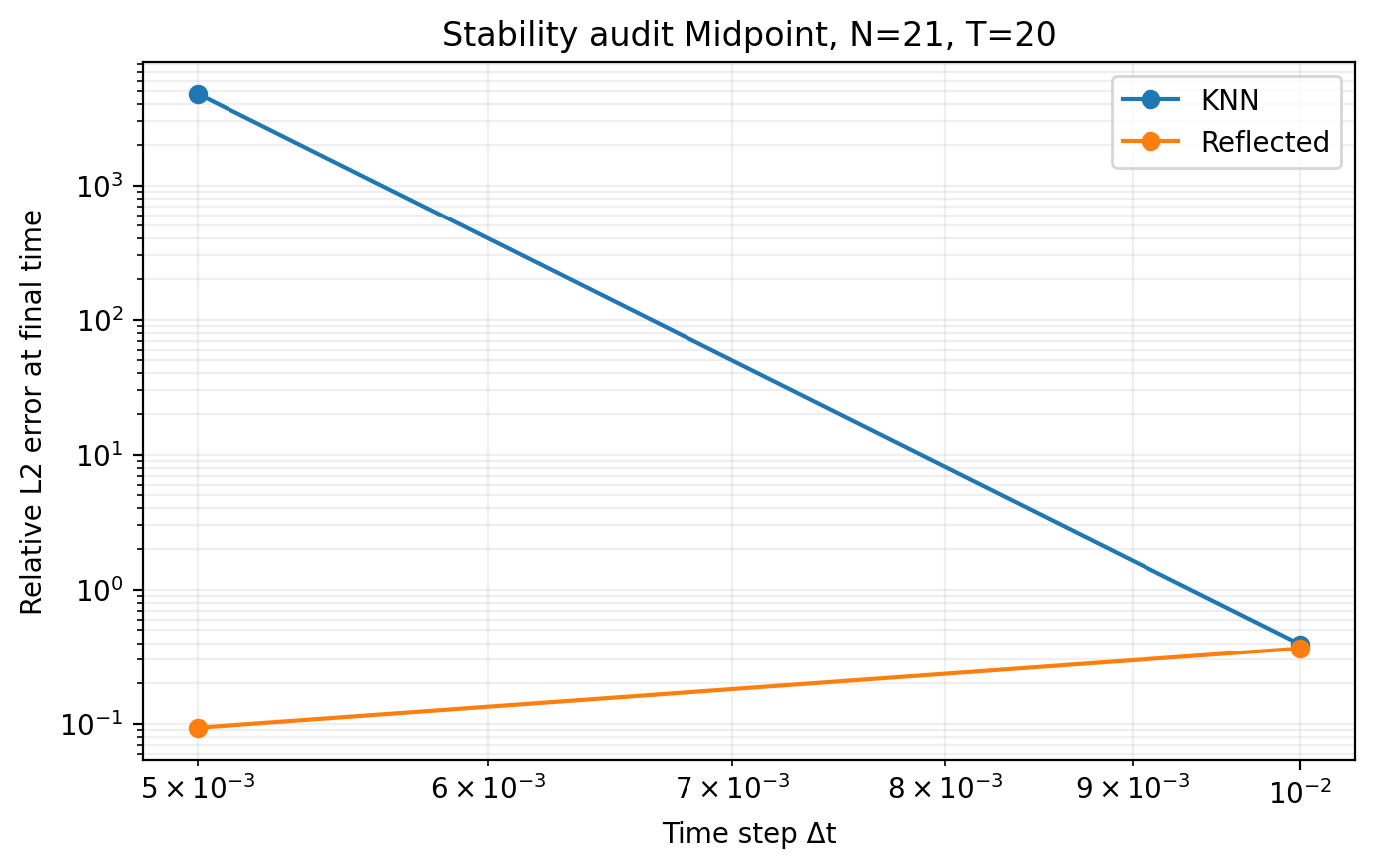}
			\caption{Midpoint step-size audit.}
			\label{fig:r11-audit-midpoint}
		\end{subfigure}
		
		\caption{Spatial spectra and long-time
			stability diagnostics for the matched
			$r^{11}+\mathbb P_5$ KNN and reflected
			RBF-FD operators at $N=21$.
			(a)--(b) Computed spatial spectra;
			(c)--(d) RK4 and implicit-midpoint
			solution-error dependence on time-step
			size at $T=20$.}
		\label{fig:spectral-stability-comparison}
	\end{figure}

	\cref{fig:time-integration-comparison}
	compares the solution accuracy and
	conservation behavior of the matched
	KNN and reflected operators at
	$N=31$ through $T=2$.
	Under RK4, the relative $L^2$ solution
	error increases approximately linearly
	with time for both discretizations,
	with the KNN error consistently larger
	than that of the reflected method.
	This behavior is consistent with the
	accumulation of spatial phase and
	mode-shape errors during propagation.
	The implicit-midpoint results reveal
	a different distinction. For the
	symmetric reflected operator, both
	mass and energy drift remain near
	numerical precision, consistent with
	the exact conservation of quadratic
	invariants by implicit midpoint.
	For KNN, the corresponding diagnostics
	exhibit bounded oscillations over
	the displayed interval, with mass
	and energy deviations of approximately
	$10^{-7}$ and $10^{-6}$, respectively.
	These results illustrate that small
	conservation defects do not necessarily
	imply small solution errors.
	The reflected construction combines
	improved solution accuracy with
	the structural conservation properties
	of its symmetric spatial operator.

	\cref{fig:spectral-stability-comparison}
	examines the spatial spectra and
	long-time integration behavior at
	$N=21$ through $T=20$.
	The reflected Laplacian has a real
	spectrum to numerical precision,
	consistent with its symmetric
	construction. In contrast, the KNN
	Laplacian possesses eigenvalues with
	nonzero imaginary parts, reflecting
	the loss of discrete self-adjointness.
	For the free Schr\"odinger equation,
	the semidiscrete evolution is governed
	by $\psi_t=iL\psi$. Consequently,
	nonreal eigenvalues of $L$ can
	introduce exponentially growing or
	decaying components, depending on
	the sign of their imaginary parts.
	This provides a mechanism by which
	spatial nonsymmetry may compromise
	long-time stability.
	The RK4 step-size audit in
	\cref{fig:r11-audit-rk4} shows severe
	error amplification for the KNN
	discretization, with relative errors
	of order $10^{14}$ over the tested
	time-step range. In contrast, the
	reflected discretization maintains
	substantially smaller errors.
	The implicit-midpoint audit in
	\cref{fig:r11-audit-midpoint}
	also shows pronounced time-step
	sensitivity for KNN, while the
	reflected method remains substantially
	more accurate over the displayed range.
	These results demonstrate that
	spatial-operator structure affects
	more than conservation diagnostics:
	it can also influence the stability
	of the semidiscrete evolution.
	However, the observed long-time
	errors depend on the initial data,
	the excited discrete modes, and the
	time-integration method.

	\subsection{Variable-potential Schr\"odinger benchmark}
	\label{subsec:variablePotential}
	
	To evaluate the reflected discretization beyond the
	sine-diagonalizable free Schr\"odinger equation, we
	consider the variable-potential problem
	\begin{equation}\label{eq:variablePotentialPDE}
		\begin{aligned}
			i\psi_t &=-\Delta\psi+V(x,y)\psi,
			&& (x,y)\in\Omega=(0,\pi)^2,\\
			\psi|_{\partial\Omega}&=0,
		\end{aligned}
	\end{equation}
	where
	\begin{equation}
		V(x,y)=5\sin^2x\sin^2y
		+2\sin x\sin(2y).
	\end{equation}
	The initial condition is a localized wave packet,
	\begin{equation}
		\psi(x,y,0)=
		\sin x\sin y\,
		\exp\!\left[
		-8(x-\pi/2)^2-8(y-\pi/2)^2
		\right]
		e^{i(3x+2y)}.
	\end{equation}
	The real, nonseparable potential couples the sine
	modes, so the full Hamiltonian is not diagonal in
	the Fourier-sine basis. This benchmark therefore
	tests the discretizations beyond the exact
	sine-mode evolution of the free problem.
	
	\paragraph{Spatial discretization and time integration.}
	We compare the 29-point odd-reflected
	$r^{11}+\mathbb P_5$ RBF-FD method, the matched
	29-point rank-safe KNN RBF-FD method, and the
	odd-reflected fourth-order finite-difference
	method (FD4).
	For a discrete Laplacian $L_h$, define
	\begin{equation}
		H_h=-L_h+\operatorname{diag}(V_j),
		\qquad V_j=V(x_j,y_j),
	\end{equation}
	so that the semidiscrete system is
	\begin{equation}
		i\dot{\boldsymbol\psi}=H_h\boldsymbol\psi.
	\end{equation}
	Since the potential is real, symmetry of $L_h$
	implies $H_h=H_h^*$. Consequently, the
	semidiscrete mass and energy,
	\begin{equation}
		M_h=\boldsymbol\psi^*\boldsymbol\psi,
		\qquad
		E_h=\boldsymbol\psi^*H_h\boldsymbol\psi,
	\end{equation}
	are conserved for the reflected and FD4
	discretizations. These conservation properties
	are not generally guaranteed for the
	nonsymmetric KNN operator.
	The semidiscrete systems are integrated using
	classical fourth-order Runge--Kutta (RK4).
	For automatic time-step selection, the
	implementation uses
	\begin{equation}
		\Delta t_{\mathrm{trial}}
		=
		\min\!\left\{
		0.002,\,
		\frac{1.5}{\max\{1,B_h\}}
		\right\},
		\qquad
		B_h=
		\max_i\sum_j |(L_h)_{ij}|
		+\max_j|V_j|.
	\end{equation}
	The number of steps and the actual step size
	are then determined by
	\begin{equation}
		n_t=
		\left\lceil
		\frac{T}{\Delta t_{\mathrm{trial}}}
		\right\rceil,
		\qquad
		\Delta t=\frac{T}{n_t}.
	\end{equation}
	This procedure reaches the prescribed final
	time exactly. The row-sum estimate provides
	a conservative step-selection heuristic,
	rather than a rigorous RK4 stability guarantee
	for arbitrary nonsymmetric operators.
	Although the symmetric spatial discretizations
	preserve mass and energy at the semidiscrete
	level, explicit RK4 does not preserve these
	quadratic invariants exactly. The measured
	conservation defects therefore also reflect
	temporal discretization errors.
	
	\paragraph{Reference solution.}
	In the absence of a known analytical solution,
	we construct an independent numerical reference
	using Fourier-sine pseudospectral kinetic
	evolution combined with a fourth-order
	symmetric composition of second-order
	kinetic--potential splitting steps
	\cite{Strang1968,Yoshida1990,
		BaoJinMarkowich2002,Lubich2008}.
	
	The reference calculation is checked by
	spatial and temporal refinement.
	Increasing the reference resolution from
	$N_{\mathrm{ref}}=241$ to $481$, with
	512 splitting steps, gives a relative
	difference of $8.95\times10^{-8}$.
	Increasing the number of splitting steps
	from 512 to 1024 on the $481\times481$
	grid gives a relative difference of
	$2.21\times10^{-6}$.
	The final reference uses
	$N_{\mathrm{ref}}=481$ and 1024 splitting
	steps. Numerical errors are evaluated
	against this reference at the coarse-grid
	nodes. The refinement comparisons provide
	empirical consistency checks rather than
	rigorous error bounds.
	
	\begin{table}[htbp]
		\centering
		\caption{Variable-potential Schr\"odinger benchmark
			at $T=0.2$: relative discrete $L^2$ errors
			and observed convergence orders.}
		\label{tab:variable_potential_accuracy}
		\small
		\setlength{\tabcolsep}{5pt}
		\begin{tabular}{r rr rr rr}
			\toprule
			& \multicolumn{2}{c}{Reflected PHS29}
			& \multicolumn{2}{c}{Matched KNN}
			& \multicolumn{2}{c}{FD4}\\
			$N$ & Error & Order & Error & Order
			& Error & Order\\
			\midrule
			21  & 7.569e-02 & ---
			& 1.823e-01 & ---
			& 2.352e-01 & ---\\
			31  & 7.401e-03 & 5.73
			& 3.049e-02 & 4.41
			& 6.799e-02 & 3.06\\
			41  & 1.539e-03 & 5.46
			& 1.121e-02 & 3.48
			& 2.369e-02 & 3.67\\
			61  & 5.764e-04 & 2.42
			& 1.647e-03 & 4.73
			& 4.915e-03 & 3.88\\
			81  & 2.289e-04 & 3.21
			& 4.015e-04 & 4.91
			& 1.579e-03 & 3.95\\
			101 & 1.034e-04 & 3.56
			& 1.424e-04 & 4.64
			& 6.512e-04 & 3.97\\
			\bottomrule
		\end{tabular}
	\end{table}

	\paragraph{Accuracy and convergence.}
	The relative discrete $L^2$ errors and
	observed convergence orders at $T=0.2$
	are reported in
	\cref{tab:variable_potential_accuracy}.
	All three methods exhibit decreasing errors
	under grid refinement. The reflected PHS29
	method achieves the smallest error at every
	tested resolution.
	At $N=101$, the relative errors are
	$1.034\times10^{-4}$,
	$1.424\times10^{-4}$, and
	$6.512\times10^{-4}$ for the reflected,
	KNN, and FD4 methods, respectively.
	Thus, the reflected method reduces the
	error by approximately $27\%$ relative
	to KNN and by a factor of $6.3$
	relative to FD4.
	The observed convergence orders vary
	across the grids, particularly for
	the two RBF-FD methods. On the finest
	grid pair, the reflected, KNN, and FD4
	methods have observed orders of
	$3.56$, $4.64$, and $3.97$, respectively.
	These values describe the measured
	fully discrete errors rather than
	establishing asymptotic spatial
	convergence rates.
	
	\begin{table}[htbp]
		\centering
		\caption{Matrix sparsity and measured C++
			runtime for the variable-potential benchmark.
			Runtime is the program-reported total;
			assembly and propagation are not timed
			separately.}
		\label{tab:variable_potential_cost}
		\begin{tabular}{r l r r}
			\toprule
			$N$ & Method & Nonzeros & Total runtime (s)\\
			\midrule
			21  & Reflected PHS29 & 9137   & 0.00578\\
			21  & Matched KNN     & 9539   & 0.01845\\
			21  & FD4             & 3021   & 0.00393\\
			31  & Reflected PHS29 & 22337  & 0.01655\\
			31  & Matched KNN     & 22959  & 0.14834\\
			31  & FD4             & 7221   & 0.00670\\
			41  & Reflected PHS29 & 41337  & 0.11089\\
			41  & Matched KNN     & 42179  & 0.13317\\
			41  & FD4             & 13221  & 0.05656\\
			61  & Reflected PHS29 & 96737  & 0.28060\\
			61  & Matched KNN     & 98019  & 0.38458\\
			61  & FD4             & 30621  & 0.18383\\
			81  & Reflected PHS29 & 175337 & 0.92850\\
			81  & Matched KNN     & 177059 & 1.10273\\
			81  & FD4             & 55221  & 0.41905\\
			101 & Reflected PHS29 & 277137 & 2.36032\\
			101 & Matched KNN     & 279299 & 2.32525\\
			101 & FD4             & 87021  & 0.78722\\
			\bottomrule
		\end{tabular}
	\end{table}

	\paragraph{Computational cost.}
	The matrix sparsity and measured
	execution times are summarized in
	\cref{tab:variable_potential_cost}.
	As indicated in \cref{tab:variable_potential_cost},
	the reflected and KNN matrices have
	comparable numbers of nonzero entries,
	whereas FD4 produces substantially
	sparser matrices.
	FD4 requires less execution time on
	the same grid, but also produces
	larger solution errors, as seen in
	\cref{tab:variable_potential_accuracy}.
	The reflected method is faster than
	KNN on most tested grids, although
	their measured runtimes are similar
	at $N=101$.
	Because the reported timings include
	both matrix construction and time
	propagation, they do not establish
	an equal-accuracy efficiency ranking.

	\begin{table}[htbp]
		\centering
		\caption{Matrix symmetry defect and relative
			mass and energy drift at $T=0.2$ for the
			variable-potential Schr\"odinger benchmark.
			Drift values are measured under explicit RK4.}
		\label{tab:variable_potential_structure}
		\small
		\setlength{\tabcolsep}{6pt}
		\begin{tabular}{r l r r r}
			\toprule
			$N$ & Method & Symmetry defect
			& Mass drift & Energy drift\\
			\midrule
			21  & Reflected PHS29 & 9.587e-14
			& 2.086e-05 & 9.336e-05\\
			21  & Matched KNN & 1.284e-01
			& 5.204e-03 & 5.205e-03\\
			21  & FD4 & 0.000e+00
			& 1.497e-05 & 5.990e-05\\
			31  & Reflected PHS29 & 9.724e-14
			& 6.669e-07 & 3.065e-06\\
			31  & Matched KNN & 1.038e-01
			& 3.799e-04 & 9.571e-04\\
			31  & FD4 & 0.000e+00
			& 5.256e-06 & 2.327e-05\\
			41  & Reflected PHS29 & 9.785e-14
			& 3.840e-08 & 1.771e-07\\
			41  & Matched KNN & 8.947e-02
			& 1.443e-05 & 2.904e-04\\
			41  & FD4 & 0.000e+00
			& 3.239e-07 & 1.473e-06\\
			61  & Reflected PHS29 & 9.852e-14
			& 6.684e-10 & 3.084e-09\\
			61  & Matched KNN & 7.273e-02
			& 5.309e-06 & 4.520e-05\\
			61  & FD4 & 0.000e+00
			& 5.896e-09 & 2.712e-08\\
			81  & Reflected PHS29 & 9.879e-14
			& 3.790e-11 & 1.749e-10\\
			81  & Matched KNN & 6.285e-02
			& 1.055e-06 & 8.240e-06\\
			81  & FD4 & 0.000e+00
			& 3.345e-10 & 1.542e-09\\
			101 & Reflected PHS29 & 9.897e-14
			& 4.055e-12 & 1.876e-11\\
			101 & Matched KNN & 5.614e-02
			& 2.556e-07 & 1.962e-06\\
			101 & FD4 & 0.000e+00
			& 3.619e-11 & 1.669e-10\\
			\bottomrule
		\end{tabular}
	\end{table}
	
	\paragraph{Structure preservation.}
	The matrix symmetry defects and
	relative mass and energy drifts
	are reported in
	\cref{tab:variable_potential_structure}.
	As shown in \cref{tab:variable_potential_structure},
	the reflected Laplacian remains symmetric
	to numerical precision, with defects of
	order $10^{-13}$, while FD4 is symmetric
	to the reported precision.
	The KNN matrix remains nonsymmetric,
	although its symmetry defect decreases
	under grid refinement.
	The reflected method also exhibits
	smaller mass and energy drift than KNN
	at every tested resolution. At $N=101$,
	the reflected mass and energy drifts
	are $4.055\times10^{-12}$ and
	$1.876\times10^{-11}$, compared with
	$2.556\times10^{-7}$ and
	$1.962\times10^{-6}$ for KNN.
	
	These findings are consistent with
	the semidiscrete conservation properties
	of the symmetric reflected operator.
	Nevertheless, because explicit RK4
	is not exactly conservative, the
	measured drift reflects both the
	spatial operator and time integration.

	\begin{figure}[pos=htbp]
		\centering
		
		\begin{subfigure}[t]{0.47\textwidth}
			\centering
			\includegraphics[width=\linewidth]
			{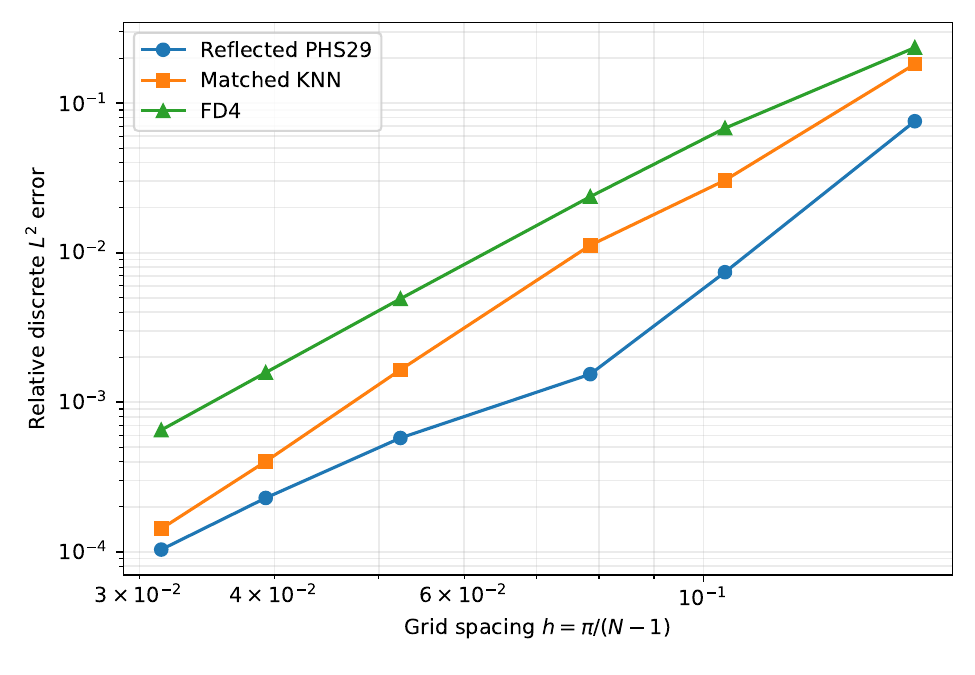}
			\caption{Relative solution-error convergence.}
			\label{fig:potential_conv}
		\end{subfigure}
		\begin{subfigure}[t]{0.44\textwidth}
			\centering
			\includegraphics[width=\linewidth]
			{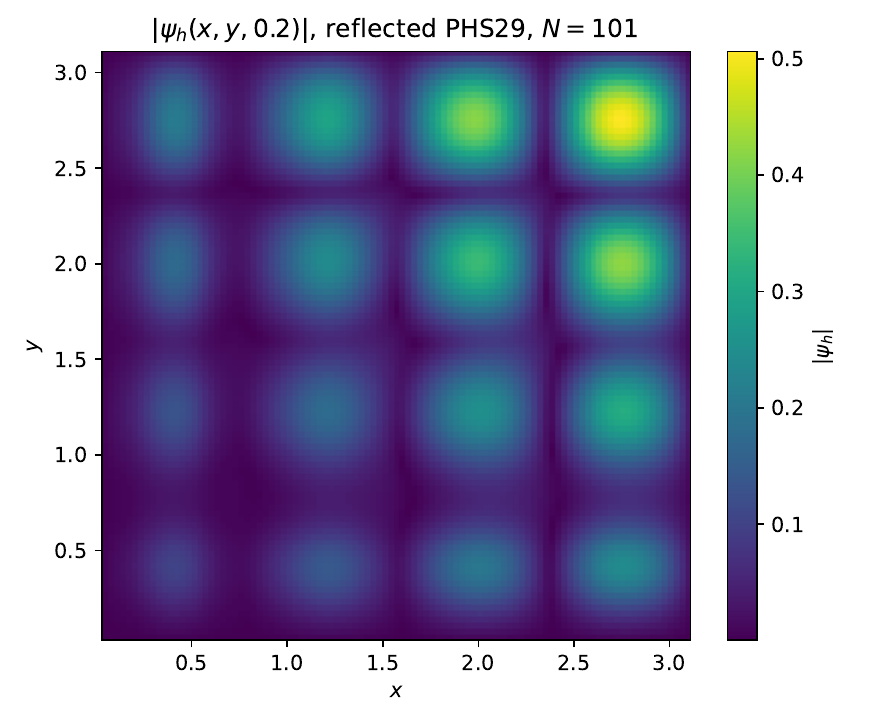}
			\caption{Solution magnitude at $T=0.2$.}
			\label{fig:potential_modulus}
		\end{subfigure}
		
		\caption{Variable-potential Schr\"odinger
			benchmark at $T=0.2$:
			(a) relative discrete $L^2$ errors
			for reflected PHS29, matched KNN,
			and FD4 across six grids;
			(b) computed solution magnitude
			$|\psi_h(x,y,T)|$ for the reflected
			method on the $N=101$ grid.}
		\label{fig:variable_potential}
	\end{figure}

	\paragraph{Solution visualization.}
	The convergence behavior and computed
	solution magnitude are presented in
	\cref{fig:variable_potential}. As illustrated in \cref{fig:potential_conv},
	all three methods exhibit decreasing
	solution errors with refinement,
	consistent with
	\cref{tab:variable_potential_accuracy}.
	The reflected method maintains the
	smallest error over the tested grids.
	\cref{fig:potential_modulus} displays
	the magnitude of the propagated
	wave packet for the reflected method
	at $N=101$ and $T=0.2$.
	Unlike the single-mode free problem,
	the evolution involves coupling
	among sine modes induced by the
	nonseparable potential.

	Taken together,
	\cref{tab:variable_potential_accuracy},
	\cref{tab:variable_potential_cost},
	\cref{tab:variable_potential_structure}
	and \cref{fig:variable_potential}
	illustrate the tradeoff between accuracy,
	computational cost, and discrete structure
	preservation. For this benchmark, the
	reflected RBF-FD method achieves the
	smallest solution errors and retains
	matrix symmetry, while requiring more
	nonzero interactions than FD4.

	\subsection{Interpretation, computational cost,
		and reproducibility}
	\label{subsec:interpretation}
	
	The numerical experiments distinguish spatial
	accuracy, discrete operator structure, and
	time-integration behavior. The reflected
	complete-shell construction has exact matrix
	symmetry and discrete sine diagonalization
	on the uniform grid, while negative
	definiteness follows separately from the
	Fourier-symbol analysis for the
	coefficient-certified kernels.
	
	The matched operator diagnostics quantify
	symmetry defects, frequency errors, and
	boundary-layer residuals across six grids.
	The free Schr\"odinger experiments examine
	modal accuracy and sparse time propagation,
	while the variable-potential benchmark
	extends the comparison to a nonseparable
	Hamiltonian. In particular,
	\cref{tab:variable_potential_accuracy},
	\cref{tab:variable_potential_structure},
	\cref{fig:variable_potential}
	show that the reflected method retains
	its accuracy and symmetry advantages
	for the variable-potential benchmark problem.
	
	\paragraph{Computational cost.}
	The 29-point reflected RBF-FD stencil
	requires more sparse matrix entries and
	matrix-vector operations than the
	standard FD4 stencil. Consequently,
	greater accuracy on the same grid
	does not necessarily imply greater
	efficiency at a fixed computational cost.
	
	The measured C++ runtimes in
	\cref{tab:variable_potential_cost}
	provide a direct comparison of the
	implemented methods. However, these
	timings include matrix construction
	and time propagation and do not
	separately measure their costs.
	Moreover, the automatic RK4 time-step
	selection may produce different numbers
	of time steps for different operators.
	Thus, the reported timings should not
	be interpreted as an equal-accuracy
	efficiency comparison.
	
	\paragraph{Scope of the numerical evidence.}
	The matched $r^{11}+\mathbb P_5$
	KNN/reflected comparison uses rank-safe
	stencil selection to maintain polynomial
	unisolvency. This procedure is not
	identical to unmodified nearest-neighbor
	selection, particularly near the boundary.
	Accordingly, the observed KNN behavior
	should not be generalized to all
	nearest-neighbor RBF-FD constructions.
	
	The parameter sweep was performed using
	the same model problems subsequently
	reported. The resulting numerical ranking
	is therefore exploratory rather than an
	independent out-of-sample validation.
	Similarly, the observed long-time KNN
	behavior does not establish a universal
	instability theorem, and the higher-grid
	spatial-operator diagnostics are not
	full sparse long-time convergence tests.
	
	Local augmented RBF-FD systems can
	be severely ill-conditioned for high
	PHS orders. High-precision weight
	generation is used to reduce the
	resulting numerical errors. The
	analytical stability proof is restricted
	to the two coefficient-certified kernels;
	the numerical experiments do not
	extend that proof to arbitrary PHS
	orders or stencil configurations.
	
	\paragraph{Reproducibility.}
	The reproducibility package includes
	weight-generation code, parameter-screening
	diagnostics, modal RK4 error tables,
	matched KNN/reflected operator diagnostics,
	and sparse RK4 and implicit-midpoint
	experiments on the stated grids.
	The principal results concern the
	$r^{11}+\mathbb P_5$ formulation.
	A historical $r^7+\mathbb P_3$
	experiment is archived separately
	and is not used to support claims
	about the principal method.
	
	\section{Discussion and Conclusions}
	\label{sec:discussionConclusions}
	
	The principal contribution of this work is a
	structure-preserving PHS RBF-FD discretization
	of the Laplacian for the two-dimensional
	Schr\"odinger equation on a square with
	homogeneous Dirichlet boundary conditions.
	The construction combines a common
	29-point complete-shell stencil with odd
	boundary reflection, allowing one set of
	reference differentiation weights to be
	used throughout the interior. An important
	distinction is that symmetry of the local
	polynomially augmented RBF system does not
	automatically imply symmetry of the
	assembled differentiation matrix.
	By enforcing reciprocal couplings through
	the stencil geometry and boundary reflection,
	the proposed method produces a symmetric
	global Laplacian without post-assembly
	symmetrization.
	
	The reflected operator admits an exact
	discrete sine eigenbasis. Its Fourier
	symbol provides a sufficient criterion
	for strict negative definiteness, with
	the required coefficient conditions
	verified for the $r^7+\mathbb P_4$
	and $r^{11}+\mathbb P_5$ formulations
	using validated arithmetic.
	These analytical results also clarify
	the semidiscrete dynamics: symmetry of
	the Laplacian yields a skew-Hermitian
	generator for the free Schr\"odinger
	equation and hence a unitary exact
	semidiscrete flow. For a real-valued
	potential represented by a diagonal
	matrix, the corresponding Hamiltonian
	remains Hermitian, although the
	complete variable-potential operator
	is generally no longer diagonal in
	the sine basis. Exact semidiscrete
	evolution conserves discrete mass
	and energy, while implicit midpoint
	preserves these quadratic invariants
	in exact arithmetic; classical RK4
	does not generally do so.
	
	The numerical experiments examine
	kernel selection, modal dispersion,
	spatial convergence, matrix symmetry,
	boundary-layer residuals, and
	time-integration behavior.
	The matched KNN/reflected comparisons
	use rank-safe polynomially unisolvent
	stencil selection and investigate the
	effects of stencil cardinality and
	distance-tie policies. Among the
	screened formulations, the
	$r^{11}+\mathbb P_5$ reflected method
	performs favorably in the reported
	modal experiments. The full sparse
	RK4 and implicit-midpoint calculations
	provide complementary evidence about
	solution accuracy, conservation,
	and long-time behavior. The
	variable-potential benchmark further
	demonstrates that the reflected method
	retains its symmetry and achieves
	the smallest relative $L^2$ errors
	among the three tested discretizations
	on the reported grids. Because the
	reference solution is numerical rather
	than analytical, these comparisons
	are subject to the stated reference
	refinement checks.
	
	The numerical results should be
	interpreted within several limitations.
	The screening and final modal tests
	use the same model problems, so the
	resulting ranking is exploratory rather
	than independent out-of-sample
	validation. The rank-safe KNN
	construction is not identical to
	unmodified nearest-neighbor selection,
	particularly near boundaries, and the
	observed behavior does not establish
	a universal instability result for
	KNN RBF-FD methods. Likewise, the
	analytical negative-definiteness
	result applies only to the two
	coefficient-certified PHS formulations.
	Modal RK4 evaluation is algebraically
	equivalent to applying RK4 to an exact
	discrete eigenmode, but its computational
	cost is not representative of full
	sparse matrix propagation.
	
	Computational efficiency also requires
	careful interpretation. The reflected
	29-point stencil has substantially
	more nonzero matrix entries than
	fourth-order finite differences.
	Consequently, improved accuracy on
	a fixed grid does not necessarily
	imply improved efficiency at a fixed
	computational cost. The reported
	wall-clock measurements are
	implementation- and hardware-dependent
	and do not establish an equal-accuracy
	efficiency ranking.
	
	The present analysis relies on a
	uniform Cartesian grid, rectangular
	geometry, and homogeneous Dirichlet
	conditions compatible with odd
	reflection. On scattered nodes or
	irregular domains, reflected or
	ghost-node stencils need not preserve
	reciprocal matrix couplings.
	Possible extensions include paired
	ghost-node constructions, weighted
	self-adjoint discretizations, mimetic
	constraints, and globally constrained
	assembly. Each would require a new
	analysis of the assembled operator
	and its spectral properties.
	Further directions include broader
	kernel classification, systematic
	time-to-accuracy comparisons, and
	extensions to nonlinear Schr\"odinger
	models with appropriate
	structure-preserving time integrators.
	
	Overall, the analysis demonstrates
	that enforcing symmetry at the
	assembled-operator level provides
	a mathematically justified route
	to structure-preserving RBF-FD
	discretizations. The resulting
	method combines exact discrete
	spectral structure with favorable
	accuracy in the tested problems,
	while making explicit the geometric
	assumptions and computational
	tradeoffs underlying these advantages.
	
	
	
	\section*{Funding}
	
	This research did not receive any specific grant from
	funding agencies in the public, commercial, or
	not-for-profit sectors.
	
	\section*{Declaration of competing interests}
	
	The authors declare that they have no known competing
	financial interests or personal relationships that
	could have appeared to influence the work reported
	in this paper.
	
	\section*{Data availability}
	
	The source code and numerical data supporting the
	findings of this study are available at
	\url{https://github.com/guangmingyao/structure-preserving-phs-rbffd}.
	
	\section*{Declaration of generative AI and
		AI-assisted technologies in the manuscript
		preparation process}
	
	During the preparation of this work, the authors
	used ChatGPT (OpenAI) to assist with manuscript
	organization, language editing, and LaTeX preparation.
	The authors reviewed and edited the resulting
	material, verified the mathematical and computational
	content, and take full responsibility for the
	content of the published article.
	\bibliographystyle{elsarticle-num}
	\bibliography{references}
	
\end{document}